\documentclass[11pt,a4paper,usenames,dvipsnames]{article}

\usepackage{url,amsmath,amssymb,latexsym,mathrsfs,comment,amsthm,enumerate,geometry,calc,ifthen,mathdots,textcomp,multicol, adjustbox}

\usepackage{cases}

\usepackage[noadjust]{cite}

\usepackage[colorlinks]{hyperref}

\usepackage[margin=10pt,font=small,labelfont=bf, labelsep=period]{caption}

\usepackage{makecell}

\usepackage[all,cmtip,2cell]{xy}

\usepackage[x11names, svgnames, rgb,table]{xcolor}
\usepackage{hhline}
\usepackage{multirow}

\usepackage{pst-node}
\usepackage{tikz-cd} 
\tikzcdset{scale cd/.style={every label/.append style={scale=#1},
    cells={nodes={scale=#1}}}}

\usepackage{enumitem}

\usepackage[utf8]{inputenc}
\usepackage[T1]{fontenc}

\usepackage{tikz}
\pgfdeclarelayer{background layer}
\pgfsetlayers{background layer,main}

\usetikzlibrary{decorations.markings}
\usetikzlibrary{arrows,matrix}
  
\tikzset{->-/.style={decoration={
  markings,
  mark=at position #1 with {\arrow{{latex}}}},postaction={decorate}}}

\newcommand\nc\newcommand
\nc\rnc\renewcommand

\nc{\custpartn}[3]{{\lower1.4 ex\hbox{
\begin{tikzpicture}[scale=.3]
\foreach \x in {#1}
{ \uvert{\x}  }
\foreach \x in {#2}
{ \lvert{\x}  }
#3 \end{tikzpicture}
}}}
\nc{\uvert}[1]{\fill (#1,2)circle(.2);}
\rnc{\lvert}[1]{\fill (#1,0)circle(.2);}

\nc\cplab[2]{\node[above] () at (#1,2) {\tiny $#2$};}

\nc\too\rightsquigarrow
\nc\comm{\textup{c}}

\nc\ca\circledast
\nc\wh\widehat

\nc\N{\mathbb N}
\nc\R{\operatorname{\mathscr R}}
\rnc\L{\operatorname{\mathscr L}}
\nc\J{\operatorname{\mathscr J}}
\nc\D{\operatorname{\mathscr D}}
\rnc\H{\operatorname{\mathscr H}}
\nc\K{\operatorname{\mathscr K}}
\nc\leqL{\leq_{\L}}
\nc\leqR{\leq_{\R}}
\nc\leqH{\leq_{\H}}
\nc\leqJ{\leq_{\J}}
\nc\leqK{\leq_{\K}}
\nc\Reg{\operatorname{Reg}}
\nc\Aut{\operatorname{Aut}}
\nc\PAut{\operatorname{PAut}}
\nc\End{\operatorname{End}}
\nc\PEnd{\operatorname{PEnd}}
\nc\Eq{\mathfrak{Eq}}

\nc\trans[1]{\left(\begin{smallmatrix}#1\end{smallmatrix}\right)}
\nc\sm\setminus

\nc{\mtlab}[1]{\mapstochar \xrightarrow {\ #1\ }}

\nc{\ssim}{\mathrel{\raise0.25 ex\hbox{\oalign{$\approx$\crcr\noalign{\kern-0.84 ex}$\approx$}}}}

\nc{\ldb}{[\hspace{-0.55truemm}[}
\nc{\rdb}{]\hspace{-0.55truemm}]}

\nc\larr{\mathrel{\begin{tikzpicture}\draw[>-](0,0)--(.5,0);\end{tikzpicture}}}
\nc\rarr{\mathrel{\begin{tikzpicture}\draw[->](0,0)--(.5,0);\end{tikzpicture}}}
\nc\llarr{\mathrel{\begin{tikzpicture}\draw[>-<](0,0)--(.5,0);\end{tikzpicture}}}
\nc\rrarr{\mathrel{\begin{tikzpicture}\draw[<->](0,0)--(.5,0);\end{tikzpicture}}}
\nc\lrarr{\mathrel{\begin{tikzpicture}\draw[>->](0,0)--(.5,0);\end{tikzpicture}}}
\nc\arrl{\mathrel{\begin{tikzpicture}\draw[-<](0,0)--(.5,0);\end{tikzpicture}}}
\nc\arrr{\mathrel{\begin{tikzpicture}\draw[<-](0,0)--(.5,0);\end{tikzpicture}}}

\nc\Z{\mathbb Z}

\nc\dda{{\scalebox{0.5}{\begin{tikzpicture}\draw[thick,>->](0,.5)--(0,0);\end{tikzpicture}}}}

\nc\pre\preceq

\newcommand{\uvx}[2]{\fill (#1,2)circle(#2);}
\newcommand{\lvx}[2]{\fill (#1,0)circle(#2);}
\newcommand{\uv}[1]{\fill (#1,2)circle(.17);}
\newcommand{\lv}[1]{\fill (#1,0)circle(.17);}

\newcommand{\uvs}[1]{{\foreach \x in {#1} { \uv{\x}}}}
\newcommand{\lvs}[1]{{\foreach \x in {#1} { \lv{\x}}}}
\newcommand{\darcx}[3]{\draw(#1,0)arc(180:90:#3) (#1+#3,#3)--(#2-#3,#3) (#2-#3,#3) arc(90:0:#3);}
\newcommand{\darc}[2]{\darcx{#1}{#2}{.4}}
\newcommand{\uarcx}[3]{\draw(#1,2)arc(180:270:#3) (#1+#3,2-#3)--(#2-#3,2-#3) (#2-#3,2-#3) arc(270:360:#3);}
\newcommand{\uarc}[2]{\uarcx{#1}{#2}{.4}}
\newcommand{\stline}[2]{\draw(#1,2)--(#2,0);}
\newcommand{\uline}[2]{\draw(#1,2)--(#2,2);}
\newcommand{\dline}[2]{\draw(#1,0)--(#2,0);}
\nc{\uarcs}[1]{{\foreach \x/\y in {#1}{ \uarc{\x}{\y} }}}
\nc{\darcs}[1]{{\foreach \x/\y in {#1}{ \darc{\x}{\y} }}}

\nc\udotted[2]{\draw[dotted](#1+.5,2)--(#2-.5,2);}
\nc\ldotted[2]{\draw[dotted](#1+.5,0)--(#2-.5,0);}

\newcommand{\uarcxcol}[4]{\draw[#4](#1,2)arc(180:270:#3) (#1+#3,2-#3)--(#2-#3,2-#3) (#2-#3,2-#3) arc(270:360:#3);}
\newcommand{\darcxcol}[4]{\draw[#4](#1,0)arc(180:90:#3) (#1+#3,#3)--(#2-#3,#3) (#2-#3,#3) arc(90:0:#3);}
\newcommand{\uvcol}[2]{\fill[#2] (#1,2)circle(.2);}

\newcommand{\stlinex}[3]{\draw[#3](#1,2)--(#2,0);}
\newcommand{\ulinex}[3]{\draw[#3](#1,2)--(#2,2);}
\newcommand{\dlinex}[3]{\draw[#3](#1,0)--(#2,0);}

\nc\Ptwo[3]{{\lower1.0 ex\hbox{
\begin{tikzpicture}[scale=.23]
\uvx1{.25}\uvx2{.25}\lvx1{.25}\lvx2{.25}
\foreach \x/\y in {#1} {\stlinex\x\y{thick}}
\foreach \x/\y in {#2} {\ulinex\x\y{thick}}
\foreach \x/\y in {#3} {\dlinex\x\y{thick}}
\end{tikzpicture}
}}}

\nc\RedPtwo[3]{{\lower1.0 ex\hbox{
\begin{tikzpicture}[red,scale=.23]
\uvx1{.25}\uvx2{.25}\lvx1{.25}\lvx2{.25}
\foreach \x/\y in {#1} {\stlinex\x\y{thick}}
\foreach \x/\y in {#2} {\ulinex\x\y{thick}}
\foreach \x/\y in {#3} {\dlinex\x\y{thick}}
\end{tikzpicture}
}}}

\nc\BluePtwo[3]{{\lower1.0 ex\hbox{
\begin{tikzpicture}[blue,scale=.23]
\uvx1{.25}\uvx2{.25}\lvx1{.25}\lvx2{.25}
\foreach \x/\y in {#1} {\stlinex\x\y{thick}}
\foreach \x/\y in {#2} {\ulinex\x\y{thick}}
\foreach \x/\y in {#3} {\dlinex\x\y{thick}}
\end{tikzpicture}
}}}

\nc\mot[2]{{\lower1.0 ex\hbox{
\begin{tikzpicture}[scale=.23]
\uv1\uv2\uv3\lv1\lv2\lv3
\foreach \x in {#1} {\stline\x\x}
\foreach \x/\y in {#2} {\uarc\x\y \darc\x\y}
\end{tikzpicture}
}}}

\nc\mott[3]{{\lower1.0 ex\hbox{
\begin{tikzpicture}[scale=.23]
\uv1\uv2\uv3\lv1\lv2\lv3
\foreach \x/\y in {#1} {\stline\x\y}
\foreach \x/\y in {#2} {\uarc\x\y}
\foreach \x/\y in {#3} {\darc\x\y}
\end{tikzpicture}
}}}

\nc\moot[2]{{\lower1.4 ex\hbox{
\begin{tikzpicture}[scale=.3]
\uv1\uv2\uv3\uv4\lv1\lv2\lv3\lv4
\foreach \x in {#1} {\stline\x\x}
\foreach \x/\y in {#2} {\uarc\x\y \darc\x\y}
\end{tikzpicture}
}}}

\newcounter{ncols}
\newcounter{incols}
\newenvironment{partn}[1]{
  \setcounter{ncols}{#1} \setcounter{incols}{\thencols - 1}\setlength{\arraycolsep}{1pt}
  \Bigl( \hspace{-1.5truemm}\scriptsize 
    \begin{array}{@{\hskip 3pt} c *{\theincols}{|c} @{\hskip 3pt}  }
}{
     \end{array}
     \normalsize \hspace{-1.5truemm}\Bigr)\setlength{\arraycolsep}{6pt}
}

\nc\la\langle
\nc\ra\rangle

\nc\ol\overline
\nc\ul\underline
\nc\da\downarrow
\nc\pr\bullet
\nc\ccirc\circledcirc
\nc\bc\odot

\nc\RSS{{\bf RSS}}
\nc\BS{{\bf BS}}
\nc\RBS{{\bf RBS}}
\nc\SBS{{\bf SBS}}
\nc\RSBS{{\bf RSBS}}
\nc\PA{{\bf PA}}
\nc\OG{{\bf OG}}
\nc\CPG{{\bf CPG}}
\nc\Set{{\bf Set}}
\nc\Sgp{{\bf Sgp}}
\nc\IS{{\bf IS}}
\nc\IG{\operatorname{\textup{\textsf{IG}}}}
\nc\RIG{\operatorname{\textup{\textsf{RIG}}}}
\nc\bG{{\bf G}}
\nc\bS{{\bf S}}
\nc\bE{{\bf E}}
\nc\bP{{\bf P}}
\nc\bEE{{\bf E}}
\nc\bPP{{\bf P}}
\nc\bB{{\bf B}}
\nc\bC{{\bf C}}
\nc\bD{{\bf D}}
\nc\bbG{\mathbb G}
\nc\bbS{\mathbb S}

\nc\PG{\operatorname{\textup{\textsf{PG}}}}
\nc\PGPz{\PG(P_0)}

\nc\ds\displaystyle
\nc\ts\textstyle

\nc\bn{{\bf n}}

\nc\al\alpha
\nc\be\beta
\nc\ga\gamma
\nc\de\delta
\nc\ve\varepsilon
\nc\lam\lambda
\nc\om\omega
\nc\ze\zeta
\nc\si\sigma
\nc\Si\Sigma
\nc\Ga\Gamma
\nc\De\Delta
\nc\Om\Omega
\nc\io\iota
\nc\ka\kappa

\nc\nab\nabla

\nc\vt\vartheta
\rnc\th\theta
\nc\Th\Theta

\nc\BY{\qquad\text{by}\qquad}
\nc\GIVENBY{\qquad\text{given by}\qquad}
\nc\ISGIVENBY{\qquad\text{is given by}\qquad}
\nc\OR{\qquad\text{or}\qquad}
\nc\Or{\quad\text{or}\quad}
\nc\AND{\qquad\text{and}\qquad}
\nc\ANDSIM{\qquad\text{and similarly}\qquad}
\nc\ANDSIm{\quad\text{and similarly}\quad}
\nc\ANDSO{\qquad\text{and so}\qquad}
\nc\ANd{\quad\text{and}\quad}
\nc\COMMA{,\qquad}
\nc\COMMa{,\quad}
\nc\WHERE{\qquad\text{where}\qquad}

\rnc\iff{\ \Leftrightarrow\ }
\nc\IFf{\quad \Leftrightarrow\quad }
\nc\Iff{\ \ \Leftrightarrow\ \ }
\nc\IFF{\qquad \Leftrightarrow\qquad }
\rnc\implies{\ \Rightarrow\ }
\nc\IMPLIES{\qquad \Rightarrow\qquad }

\nc\set[2]{\{#1:#2\}}
\nc\bigset[2]{\big\{#1:#2\big\}}
\nc\pres[2]{\la#1:#2\ra}

\nc\bit{\begin{itemize}[label=\textbullet, leftmargin=5mm]}
\nc\eit{\end{itemize}}
\nc\ben{\begin{enumerate}[label=\textup{(\roman*)},leftmargin=10mm]}
\nc\bena{\begin{enumerate}[label=\textup{(\alph*)},leftmargin=10mm]}
\nc\een{\end{enumerate}}
\nc\bmc{\begin{multicols}}
\nc\emc{\end{multicols}}

\nc\pf{\begin{proof}}
\nc\epf{\end{proof}}
\nc\pfclaim{\begin{quote}\begin{proof}}
\nc\epfclaim{\end{proof}\end{quote}}
\nc\epfres{\hfill\qed}
\nc\epfreseq{\tag*{\qed}}
\let\oldproofname=\proofname
\renewcommand{\proofname}{\rm\bf{\oldproofname}}
\nc{\pfitem}[1]{\medskip \noindent #1.}
\nc{\firstpfitem}[1]{#1.}
\nc{\pfcase}[1]{\medskip\noindent {\bf Case #1.}}
\nc\aftercases{\medskip\noindent}

\nc\rH{\mathrel{\H}}
\nc\rL{\mathrel{\L}}
\nc\rR{\mathrel{\R}}
\nc\rD{\mathrel{\D}}
\nc\rJ{\mathrel{\J}}
\nc\rK{\mathrel{\K}}
\nc\rsi{\mathrel{\si}}

\nc\leqF{\leq_{\F}}
\nc\geqF{\geq_{\F}}

\newcommand{\Sing}{\operatorname{Sing}}

\newcommand{\dom}{\operatorname{dom}}
\newcommand{\codom}{\operatorname{codom}}
\newcommand{\coker}{\operatorname{coker}}
\newcommand{\rank}{\operatorname{rank}}

\nc\mr\mathrel
\nc\pc[2]{(#1,#2)^\sharp}

\nc\U{{\bf U}}
\nc\bF{{\bf F}}
\nc\GG{{\bf G}}

\nc\V{\mathcal V}
\nc\G{\mathcal G}
\rnc\iff{\ \Leftrightarrow\ }
\rnc\implies{\ \Rightarrow\ }
\nc\Implies{\quad \Rightarrow\quad }
\nc\F{\mathrel{\mathscr F}}
\nc\C{\mathscr C}
\nc\M{\mathcal M}
\nc\CC{\mathcal C}
\nc\CP{\CC_P}
\nc\CPz{\CC_{P_0}}
\nc\EP{\bE(P)}
\nc\EPd{\bE(P')}
\nc\PE{\bP(E)}
\nc\PEd{\bP(E')}
\nc\DD{\mathcal D}
\nc\I{\mathcal I}
\rnc\O{\mathcal O}
\rnc\S{\mathcal S}
\rnc\P{\mathcal P}
\nc\TL{\mathcal T\!\mathcal L}
\nc\PP{\mathscr P\mathcal P}
\nc\T{\mathcal T}
\nc\p{\mathfrak p}
\nc\q{\mathfrak q}
\rnc\r{\mathfrak r}
\nc\s{\mathfrak s}
\rnc\t{\mathfrak t}
\nc\bd{{\bf d}}
\nc\br{{\bf r}}
\nc\op\oplus
\nc\lra{\mathrel\leftrightarrow}
\nc\es\varnothing
\nc\rev{\textup{rev}}
\nc\corestt{{\upharpoonleft}}
\nc\restt{{\upharpoonright}}
\nc\corest{{\downharpoonleft}}
\nc\rest{{\upharpoonright}}
\nc\WHERe{\quad\text{where}\quad}
\nc\ot\leftarrow
\rnc\a{\mathfrak a}
\rnc\b{\mathfrak b}
\rnc\c{\mathfrak c}
\rnc\d{\mathfrak d}
\nc\sub\subseteq
\nc\mt\mapsto
\nc\im{\operatorname{im}}
\nc\B{\mathcal B}
\nc\PB{\P\B}
\nc\PT{\P\T}
\nc\E{\mathbb E}
\nc\Ef{\E^\flat}
\nc\BP{\operatorname{\textup{\textsf{BP}}}}
\rnc\SS{\operatorname{\textup{\textsf{S}}}}
\nc\MM{\operatorname{\textup{\textsf{M}}}}

\numberwithin{equation}{section}

\newtheorem{thm}[equation]{Theorem}
\newtheorem{lemma}[equation]{Lemma}
\newtheorem{cor}[equation]{Corollary}
\newtheorem{prop}[equation]{Proposition}

\theoremstyle{definition}

\newtheorem{defn}[equation]{Definition}
\newtheorem{rem}[equation]{Remark}

\newcounter{caseco}

\newcounter{subcaseco}

\newcounter{stepco}

\newcounter{stageco}

\makeatletter
\def\blfootnote{\gdef\@thefnmark{}\@footnotetext}
\makeatother

\begin{document}

\title{\vspace{-10mm}Twisted products of monoids}

\date{}
\author{}

\maketitle

\vspace{-15mm}

\begin{center}
{\large 
James East,%
\hspace{-.25em}\footnote{\label{fn:JE}Centre for Research in Mathematics and Data Science, Western Sydney University, Locked Bag 1797, Penrith NSW 2751, Australia. {\it Emails:} {\tt J.East@westernsydney.edu.au}, {\tt A.ParayilAjmal@westernsydney.edu.au}.}
Robert D.~Gray,%
\hspace{-.25em}\footnote{School of Engineering, Mathematics and Physics, University of East Anglia, Norwich NR4 7TJ, England, UK. {\it Email:} {\tt Robert.D.Gray@uea.ac.uk}.}
P.A.~Azeef Muhammed,%
\hspace{-.25em}\textsuperscript{\ref{fn:JE}}
Nik Ru\v{s}kuc%
\footnote{Mathematical Institute, School of Mathematics and Statistics, University of St Andrews, St Andrews, Fife KY16 9SS, UK. {\it Email:} {\tt Nik.Ruskuc@st-andrews.ac.uk}}
\blfootnote{This work was supported by the following grants:
Future Fellowship FT190100632 of the Australian Research Council;
EP/V032003/1, EP/S020616/1 and EP/V003224/1 of the Engineering and Physical Sciences Research Council.
The first author thanks the Heilbronn Institute for partially funding his visit to St Andrews in 2025.
The second author thanks the Sydney Mathematical Research Institute, the University of Sydney
and Western Sydney University for partially funding his visit to Sydney in 2023.
}
}
\end{center}

\maketitle

\begin{abstract}
\noindent
A twisting of a monoid $S$ is a map $\Phi:S\times S\to\mathbb N$ satisfying the identity $\Phi(a,b) + \Phi(ab,c) = \Phi(a,bc) + \Phi(b,c)$.  Together with an additive commutative monoid $M$, and a fixed $q\in M$, this gives rise a so-called twisted product $M\times_\Phi^qS$, which has underlying set $M\times S$ and multiplication $(i,a)(j,b) = (i+j+\Phi(a,b)q,ab)$.  This construction has appeared in the special cases where $M$ is $\N$ or $\Z$ under addition, $S$ is a diagram monoid (e.g.~partition, Brauer or Temperley-Lieb), and  $\Phi$ counts floating components in concatenated diagrams.

In this paper we identify a special kind of `tight' twisting, and give a thorough structural description of the resulting twisted products.  This involves characterising Green's relations, (von Neumann) regular elements, idempotents, biordered sets, maximal subgroups, Sch\"utzenberger groups, and more.  We also consider a number of examples, including several apparently new ones, which take as their starting point certain  generalisations of  Sylvester's rank inequality from linear algebra.

\medskip

\noindent
\emph{Keywords}: Twistings, twisted products, diagram monoids, linear monoids, independence algebras.
\medskip

\noindent
MSC: 
20M10,  
20M20.  
\end{abstract}

\begin{center}
\emph{Dedicated to Prof.~Mikhail V.~Volkov on the occasion of his 70th birthday.}
\end{center}

\tableofcontents

\section{Introduction}\label{sect:intro}

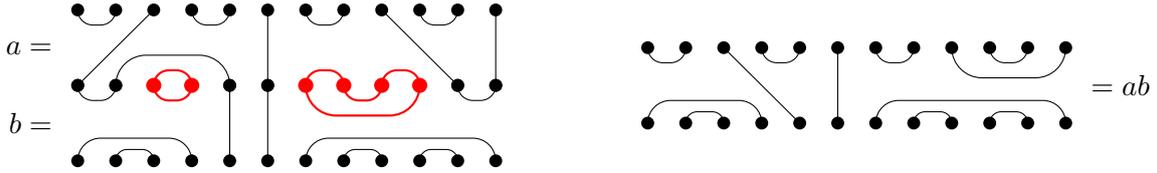
\begin{figure}[t]
\begin{center}
\begin{tikzpicture}[scale=.5]

\begin{scope}[shift={(0,0)}]	
\uvs{1,...,12}
\lvs{1,...,12}
\uarc12
\uarc45
\uarc78
\uarc{10}{11}
\darcx25{.8}
\darcxcol34{.4}{red,thick}
\darcxcol78{.4}{red,thick}
\darcxcol9{10}{.4}{red,thick}
\stline31
\stline66
\stline9{11}
\stline{12}{12}
\draw(0.6,1)node[left]{$a=$};
\end{scope}

\begin{scope}[shift={(0,-2)}]	
\uvs{1,...,12}
\lvs{1,...,12}
\uarc12
\uarcxcol34{.4}{red,thick}
\uarc{11}{12}
\uarcxcol7{10}{.8}{red,thick}
\uarcxcol89{.4}{red,thick}
\darcx14{.6}
\darcx7{12}{.6}
\darcx23{.3}
\darcx89{.3}
\darcx{10}{11}{.3}
\stline55
\stline66
\uvcol3{red}
\uvcol4{red}
\uvcol7{red}
\uvcol8{red}
\uvcol9{red}
\uvcol{10}{red}
\draw(0.6,1)node[left]{$b=$};
\end{scope}

\begin{scope}[shift={(15,-1)}]	
\uvs{1,...,12}
\lvs{1,...,12}
\uarc12
\uarc45
\uarc78
\uarcx9{12}{.8}
\uarc{10}{11}
\darcx14{.6}
\darcx7{12}{.6}
\darcx23{.3}
\darcx89{.3}
\darcx{10}{11}{.3}
\stline35
\stline66
\draw(12.4,1)node[right]{$=ab$};
\end{scope}

\end{tikzpicture}
\caption{Stacking Temperley--Lieb diagrams $a$ and $b$ (left).  Here there are two floating loops (shown in red), so~$\Phi(a,b)=2$.  Removing these components leads to the product $ab$ (right).}
\label{fig:TL}
\end{center}
\end{figure}

Twistings  arise in many parts of mathematics, under many different names, such as cocycles and multipliers \cite{Ebanks2015, DE1995, PS2002, Hosszu1971, Clark1967, ABS2024, KPS2015, Brown1994, ACCLMR2022, BFM2001}.  Notably,  diagram algebras are twisted semigroup algebras over the corresponding diagram monoids.  Let us give an example. If~$\Phi(a,b)$ denotes the number of `floating loops' when Temperley--Lieb diagrams $a$ and $b$ are stacked together to form the product $ab$ (see Figure \ref{fig:TL}), then the following identity holds:
\begin{equation}
\label{eq:twid}
\Phi(a,b) + \Phi(ab,c) = \Phi(a,bc) + \Phi(b,c).
\end{equation}
The set of all formal $\mathbb C$-linear combinations of such diagrams can then be turned into a twisted semigroup	 algebra, in which the product of basis elements is given by
\[
a \star b = \xi^{\Phi(a,b)} ab,
\]
for some fixed $\xi\in\mathbb C$; this is the Temperley--Lieb algebra \cite{TL1971}.
Analogous statements hold for other diagram algebras, such as partition, Brauer and Motzkin algebras
 \cite{Brauer1937, MM2014, Martin1994, Jones1994_2, BH2014}.  
This has proved to be a fruitful way to study diagram algebras, for example by creating a link between the cellular  structure of the algebra and  the ideals and Green's relations of the monoid \cite{Wilcox2007, JEgrpm, EG2017, DEG2017, HR2005, Kauffman1990}.

To mirror the fundamental idea behind  twisted algebras, one can define a so-called  twisted diagram monoid, with underlying set $\N\times S$, and product
\[
(i,a) ( j,b) = (i+j+\Phi(a,b),ab).
\]
Associativity follows from the twisting identity \eqref{eq:twid}.  Twisted diagram monoids have been studied frequently in the literature; see for example \cite{ER2022c, ER2022b, EGPAR2025, BDP2002,DE2018,KV2023,DE2017,KV2019,CHKLV2019,ACHLV2015,LF2006}.
It was advantageous in \cite{EGPAR2025,KV2019,KV2023} to embed these into larger twisted monoids with underlying sets $\Z\times S$ and the same multiplication.  It is in fact possible to change $\N$ into an arbitrary commutative monoid $M$, and then for a fixed parameter $q\in M$, define the product
\[
(i,a) ( j,b) = (i+j+\Phi(a,b)q,ab),
\]
resulting in what we will call a twisted product, $M\times_\Phi^qS$.
This can in turn be done for any monoid~$S$ with a twisting, i.e.~a function $\Phi:S\times S\to\N$ satisfying \eqref{eq:twid}, allowing one to capture many more important examples, including transformation monoids, linear monoids, and more general monoids of (partial) endomorphisms/automorphisms of independence algebras.

The purpose of the current paper is to initiate a systematic study of twisted products $M\times_\Phi^qS$.  
Specifically, we will investigate Green's relations, idempotents, subgroups, and regularity properties.
We obtain the sharpest results by identifying a key `tight' property of a twisting.  The above-mentioned `loop-counting' twisting is tight for the partition monoid $\P_n$, the planar partition monoid $\PP_n$, the Brauer monoid~$\B_n$ and the Temperley--Lieb monoid $\TL_n$, but loose for the partial Brauer monoid~$\PB_n$ and the Motzkin monoid $\M_n$.  The resulting difference in structure can be glimpsed by comparing Figure~\ref{fig:P2Phi}~(middle) to Figure \ref{fig:PB2Phi} (middle), which respectively show a tight twisted product $M\times_\Phi^q\P_2$ and a loose twisted product $M\times_\Phi^q\PB_2$.  Formal explanations of what these diagrams represent will be given below.

The article is organised as follows.    We begin in Section \ref{sect:prelim} with preliminary/background material, and then give definitions and basic properties of (tight) twistings and twisted products in Section \ref{sect:basic}.  We look at a number of examples of twistings in Section \ref{sect:eg}, starting with the canonical float-counting twisting of diagram monoids, and then a new family of twistings that we call `rigid'.  
The latter  are based on natural  generalisations of Sylvester's rank inequality from linear algebra, and apply to diagram monoids and to endomorphism monoids of certain independence algebras.  Curiously, these rigid twistings are tight for the Brauer monoid, but loose for the other diagram monoids listed above, including the partition and Temperley--Lieb monoids.
We then proceed to discuss the structure of a tight twisted product, specifically by characterising Green's relations (Section \ref{sect:GR}), idempotents and biordered sets (Section~\ref{sect:ET}), and (von Neumann) regularity (Section \ref{sect:reg}).  The results of these three sections  are fully general, in the sense that they apply to any tight twisted product.  We conclude in Section \ref{sect:IG} by considering a problem that 
crucially depends on the actual twisted product under consideration,
namely the determination of the idempotent-generated submonoid.  We completely describe these submonoids for rigid twisted products over groups, and for partition monoids with respect to the canonical twisting, and also compare these to existing results in the literature.

\section{Preliminaries}\label{sect:prelim}

We now give the preliminary definitions we need concerning semigroups (Section \ref{subsect:S}), transformation and diagram monoids (Section \ref{subsect:TD}) and independence algebras (Section \ref{subsect:A}).  For more details, see for example \cite{Howie1995,CPbook,EG2017,Gould1995}.
Throughout the paper we write $\N = \{0,1,2,\ldots\}$ and $\Z = \{0,\pm1,\pm2,\ldots\}$ for the sets of natural numbers and integers, respectively.


\subsection{Semigroups}\label{subsect:S}

Let $S$ be a semigroup, and denote by $S^1$ its monoid completion.  That is, $S^1 = S$ if $S$ is a monoid, or else $S = S\cup\{1\}$, where $1\not\in S$ acts as an adjoined identity element. Green's $\L$, $\R$ and $\J$ pre-orders and equivalences are defined, for $a,b\in S$, by
\begin{align*}
a \leqL b &\iff S^1a\sub S^1b , & a \L b &\iff S^1a= S^1b,\\
a \leqR b &\iff aS^1\sub bS^1 , & a \R b &\iff aS^1= bS^1 ,\\
a \leqJ b &\iff S^1aS^1\sub S^1bS^1 , & a \J b &\iff S^1aS^1= S^1bS^1.
\end{align*}
From these, we can also define the pre-order ${\leqH} = {\leqL}\cap{\leqR}$, and the equivalences ${\H} = {\L}\cap{\R}$ and ${\D} = {\L}\vee{\R}$.  The latter denotes the join of $\L$ and $\R$ in the lattice of all equivalence relations of $S$, i.e.~the transitive closure of ${\L}\cup{\R}$.  It is well known that ${\D} = {\L}\circ{\R} = {\R}\circ{\L}$ (where $\circ$ is the usual composition operation on binary relations), and that ${\D}={\J}$ when $S$ is finite.  If~$S$ is commutative, then all five of Green's equivalences are equal, i.e.~${\L}={\R}={\J}={\H}={\D}$, and similarly the four pre-orders are equal.

If $\K$ is any of $\L$, $\R$, $\J$, $\H$ or $\D$, we denote by $K_a$ the $\K$-class of an element $a\in S$.  If ${\K}\not={\D}$, then the set $S/{\K}$ of all such $\K$-classes is partially ordered by
\[
K_a \leq K_b \iff a\leqK b \qquad\text{for $a,b\in S$.}
\]

Aspects of the structure of a finite semigroup can be conveniently visualised by using a so-called \emph{egg-box diagram}.  Elements of each ${\D}={\J}$-class are drawn in a rectangular array, with rows containing $\R$-related elements, columns containing $\L$-related elements, and hence cells (intersections of rows and columns) containing $\H$-related elements.  The $\leq$-relationships between ${\D}={\J}$ classes are also indicated as a Hasse diagram.  Several examples can be seen in Figures \ref{fig:P2Phi}--\ref{fig:B34Phi}.

An element $a$ of a semigroup $S$ is (von Neumann) \emph{regular} if $a\in aSa$, i.e.~if $a=aba$ for some~$b\in S$.  
Note then that the element $c=bab$ satisfies $a=aca$ and $c=cac$; in this case,~$a$ and~$c$ are said to be (semigroup) \emph{inverses} of each other.  If $a$ is regular, then so too is every element of its $\D$-class $D_a$.  We write $\Reg(S)$ for the set of all regular elements of $S$, and we say that $S$ is regular if $S=\Reg(S)$.  The semigroup $S$ is said to be \emph{inverse} if every element $a$ has a unique inverse, denoted $a^{-1}$; it is well known that $a\mt a^{-1}$ is an \emph{involution} of~$S$, meaning that $(a^{-1})^{-1} = a$ and $(ab)^{-1} = b^{-1}a^{-1}$ for all $a,b\in S$.

An \emph{idempotent} of a semigroup $S$ is an element $e\in S$ satisfying $e=e^2$.  The set of all idempotents is denoted $E(S)$.  The $\H$-class of any idempotent is a group.  All group $\H$-classes contained in a common $\D$-class are isomorphic.  If $D$ is a regular $\D$-class, then every $\L$-class contained in $D$ contains at least one idempotent, and similarly for the $\R$-classes in $D$.  
In egg-box diagrams, cells containing idempotents (i.e.~the group $\H$-classes) are shaded grey; again see Figures \ref{fig:P2Phi}--\ref{fig:B34Phi}.
When $S$ is a monoid with identity $1$, the $\H$-class $H_1$ is called the \emph{group of units} of~$S$.
We say that $S$ is \emph{unipotent} if $E(S)=\{1\}$.

\subsection{Transformation and diagram monoids}\label{subsect:TD}

For a set $X$, we denote by
$\PT_X$ the \emph{partial transformation monoid},
$\T_X$ the \emph{full transformation monoid}, 
$\I_X$ the  \emph{symmetric inverse monoid}, and
$\S_X$ the  \emph{symmetric group}.
These consist, respectively, of all partial transformations $X\to X$, full transformations $X\to X$, partial bijections $X\to X$, and bijections $X\to X$, under composition in each case.  All four monoids are regular, and moreover $\I_X$ is inverse and $\S_X$ is a group.  When $X = \{1,\ldots,n\}$ for some integer $n\geq1$, we write $\PT_n = \PT_X$, and similarly for $\T_n$, $\I_n$ and $\S_n$.

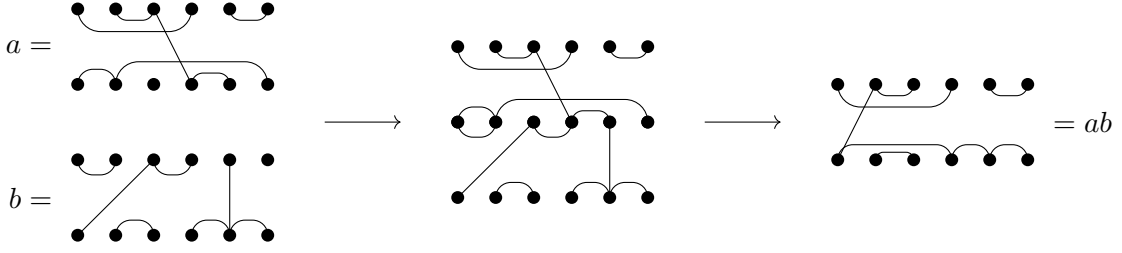
\begin{figure}[t]
\begin{center}
\begin{tikzpicture}[scale=.5]

\begin{scope}[shift={(0,0)}]	
\uvs{1,...,6}
\lvs{1,...,6}
\uarcx14{.6}
\uarcx23{.3}
\uarcx56{.3}
\darc12
\darcx26{.6}
\darcx45{.3}
\stline34
\draw(0.6,1)node[left]{$a=$};
\draw[->](7.5,-1)--(9.5,-1);
\end{scope}

\begin{scope}[shift={(0,-4)}]	
\uvs{1,...,6}
\lvs{1,...,6}
\uarc12
\uarc34
\darc45
\darc56
\darc23
\stline31
\stline55
\draw(0.6,1)node[left]{$b=$};
\end{scope}

\begin{scope}[shift={(10,-1)}]	
\uvs{1,...,6}
\lvs{1,...,6}
\uarcx14{.6}
\uarcx23{.3}
\uarcx56{.3}
\darc12
\darcx26{.6}
\darcx45{.3}
\stline34
\draw[->](7.5,0)--(9.5,0);
\end{scope}

\begin{scope}[shift={(10,-3)}]	
\uvs{1,...,6}
\lvs{1,...,6}
\uarc12
\uarc34
\darc45
\darc56
\stline31
\stline55
\darc23
\end{scope}

\begin{scope}[shift={(20,-2)}]	
\uvs{1,...,6}
\lvs{1,...,6}
\uarcx14{.6}
\uarcx23{.3}
\uarcx56{.3}
\darc14
\darc45
\darc56
\stline21
\darcx23{.2}
\draw(6.4,1)node[right]{$=ab$};
\end{scope}

\end{tikzpicture}
\caption{Multiplication of partitions $a,b\in\P_6$, with the product graph $\Pi(a,b)$ in the middle.}
\label{fig:P6}
\end{center}
\end{figure}

A further key class of examples are the \emph{diagram monoids}, which are all submonoids of the partition monoids, to which we now turn.  Fix a positive integer $n$, and write $\bn=\{1,\ldots,n\}$ and~${\bn'=\{1',\ldots,n'\}}$.  The  \emph{partition monoid} $\P_n$ consists of all set partitions of $\bn\cup\bn'$.  A partition $a\in\P_n$ is identified with any graph on vertex set $\bn\cup\bn'$ whose connected components are the blocks of $a$; such a graph is typically drawn with vertices $1<\cdots<n$ on an upper row and $1'<\cdots<n'$ on a lower row.  Some example partitions from $\P_6$ are shown in Figure \ref{fig:P6}:
\begin{equation}
\label{eq:partab}
\begin{split}
a = \big\{ \{1,4\} , \{2,3,4',5'\} , \{5,6\} , \{1',2',6'\} , \{3'\} \big\},\\
b = \big\{ \{1,2\} , \{3,4,1'\} , \{5,4',5',6'\} , \{6\} , \{2',3'\} \big\}.
\end{split}
\end{equation}

\newpage

The product of partitions $a,b\in\P_n$ is calculated as follows.  First, let $\bn''=\{1'',\ldots,n''\}$, and define three further graphs:
\bit
\item $a^\vee$ on vertex set $\bn\cup\bn''$, obtained by changing each lower vertex $x'$ of $a$ to $x''$,
\item $b^\wedge$ on vertex set $\bn''\cup\bn'$, obtained by changing each upper vertex $x$ of $b$ to $x''$,
\item $\Pi(a,b)$ on vertex set $\bn\cup\bn''\cup\bn'$, whose edge set is the union of those of $a^\vee$ and $b^\wedge$.
\eit
The graph $\Pi(a,b)$ is called the \emph{product graph} of $a$ and $b$; the vertices $1''<\cdots<n''$ are drawn in the middle row. 
The product $ab\in\P_n$ is then defined to be the partition for which $x,y\in\bn\cup\bn'$ belong to the same block of $ab$ if and only if $x$ and $y$ belong to the same connected component of $\Pi(a,b)$.   Calculation of the product $ab$ in $\P_6$, with $a$ and  $b$ as in \eqref{eq:partab}, is shown in Figure \ref{fig:P6}.

Note that the product graph $\Pi(a,b)$ may contain  \emph{floating components}, by which we mean components consisting entirely of double-dashed elements.  The number $\Phi(a,b)$ of such floating components will play a key role throughout.  For example, if $a,b\in \P_6$ are as in \eqref{eq:partab},  then~${\Phi(a,b) = 1}$, with the unique floating component of $\Pi(a,b)$ being $\{1'',2'',6''\}$; see Figure \ref{fig:P6}.

The partition monoid $\P_n$ contains many important submonoids, including:
\bit
\item the  \emph{Brauer monoid} $\B_n = \set{a\in\P_n}{\text{each block of $a$ has size $2$}}$,
\item the  \emph{partial Brauer monoid} $\PB_n = \set{a\in\P_n}{\text{each block of $a$ has size $\leq2$}}$,
\item the  \emph{planar partition monoid} $\PP_n = \set{a\in\P_n}{\text{$a$ is planar}}$,
where planarity here means that~$a$ can be drawn with no crossings, with all edges contained in the rectangle bounded by the vertices,
\item the  \emph{Temperley--Lieb monoid} $\TL_n = \B_n\cap\PP_n$, and
\item the  \emph{Motzkin monoid} $\M_n = \PB_n \cap \PP_n$.
\eit
An example of a planar partition is  $b\in\P_6$ pictured in Figure \ref{fig:P6}.

A block $A$ of a partition $a\in\P_n$ is called
an  \emph{upper} or  \emph{lower non-transversal} if $A\sub\bn$ or~$A\sub\bn'$, respectively.
Any other block is a \emph{transversal}.
We write
\[
a = \begin{partn}{6}A_1&\cdots&A_r&C_1&\cdots&C_s \\ \hhline{~|~|~|-|-|-} B_1&\cdots&B_r&D_1&\cdots&D_t\end{partn}
\]
to indicate that $a$ has transversals $A_i\cup B_i'$ (for $1\leq i\leq r$), upper non-transversals $C_i$ (for $1\leq i\leq s$) and lower non-transversals $D_i'$ (for $1\leq i\leq t$).  Any of $r$, $s$ or $t$ could be $0$, but not all three since we are assuming that $n\geq1$.

The partition monoid $\P_n$ comes equipped with a natural involution ${}^*$, defined by turning diagrams up-side down:
\[
a = \begin{partn}{6}A_1&\cdots&A_r&C_1&\cdots&C_s \\ \hhline{~|~|~|-|-|-} B_1&\cdots&B_r&D_1&\cdots&D_t\end{partn}  \mt
\begin{partn}{6} B_1&\cdots&B_r&D_1&\cdots&D_t\\ \hhline{~|~|~|-|-|-} A_1&\cdots&A_r&C_1&\cdots&C_s\end{partn}=a^*.
\]
This satisfies the identities $a^{**} = a = aa^*a$ and $(ab)^* = b^*a^*$, so that $\P_n$ is a  \emph{regular $*$-monoid}~\cite{NS1978}.  In particular, $\P_n$ is regular, but it is not inverse for $n\geq2$.  The submonoids~$\B_n$,~$\PB_n$,~$\PP_n$,~$\TL_n$ and~$\M_n$ are all closed under ${}^*$, so they are also regular $*$-monoids.

The  (\emph{co})\emph{domain}, (\emph{co})\emph{kernel} and  \emph{rank} of a partition $a\in\P_n$ are defined by
\begin{align*}
\dom(a) &= \set{x\in\bn}{x\text{ is contained in a transversal of }a},\\
\codom(a) &= \set{x\in\bn}{x'\text{ is contained in a transversal of }a},\\
\ker(a) &= \set{(x,y)\in\bn\times\bn}{\text{$x$ and $y$ belong to the same block of $a$}},\\
\coker(a) &= \set{(x,y)\in\bn\times\bn}{\text{$x'$ and $y'$ belong to the same block of $a$}},\\
\rank(a) &= \text{the number of transversals of $a$}.
\end{align*}
For example, $a\in\P_6$ in Figure \ref{fig:P6} has rank $1$, domain $\{2,3\}$, and kernel-classes $\{1,4\}$, $\{2,3\}$ and~$\{5,6\}$.

The above parameters can be used to characterise Green's relations on any of the above diagram monoids.  In particular, it will be useful to bear in mind the following:
\begin{itemize}
\item
$a\R b \iff [\dom(a)=\dom(b) \text{ and }\ker(a) = \ker(b)]$,
\item
$a\L b \iff [\codom(a)=\codom(b) \text{ and }\coker(a) = \coker(b)]$,
\item
$a\J b \iff \rank(a)=\rank(b)$,
\item
$a\leqJ b \iff \rank(a) \leq \rank(b)$.
\end{itemize}
For full details, see \cite{LF2006,Wilcox2007}.  In particular, if $S$ is any of $\P_n$, $\PB_n$, $\B_n$, $\PP_n$, $\M_n$ or $\TL_n$, then the ${\D}={\J}$-classes of $S$ are the sets
\[
D_r = D_r(S) = \set{a\in S}{\rank(a) = r} \qquad\text{for $0\leq r\leq n$,}
\]
with the additional constraint that $r\equiv n$ (mod $2$) if $S$ is either $\B_n$ or $\TL_n$.  The ordering on these is given by $D_r\leq D_s \iff r\leq s$.

Finally, we note that the transformation monoids $\T_n$ and $\I_n$ (and hence also $\S_n$)  can be regarded as submonoids of $\P_n$:
\[
\T_n = \set{a\in\P_n}{\dom(a)=\bn\text{ and }\ker(a) = \De_\bn}
\ANd
\I_n = \set{a\in\P_n}{\ker(a) = \De_\bn = \coker(a)}.
\]
Here and elsewhere we write $\De_X = \set{(x,x)}{x\in X}$ for the trivial/equality relation on a set $X$.  Note that $\I_n$ is closed under the involution ${}^*$ (which restricts to the inverse operation on $\I_n$), but $\T_n$ is not.

\subsection{Independence algebras and their (partial) endomorphism monoids}\label{subsect:A}

Another source of examples for us will be certain twistings on semigroups of transformations and matrices. They  can in fact be treated within a common framework of endomorphism monoids of independence algebras. These algebras were first introduced by Narkiewicz \cite{Nark1961} under the name of $v^*$-algebras. The study of their endomorphisms, and indeed  the term `independence algebras', were introduced by Gould \cite{Gould1995}. For more details on the history of this topic, and connections with other areas of mathematics, see \cite{Araujo2011}.

Let $A$ be an algebra (in the sense of universal algebra).  An  \emph{endomorphism} of~$A$ is a morphism $A\to A$, i.e.~a self-map of $A$ that respects the operations of $A$.  A  \emph{partial endomorphism} of~$A$ is a morphism $B\to A$ for some subalgebra $B\leq A$.  The sets $\End(A)$ and~$\PEnd(A)$ of all such (partial) endomorphisms are monoids under composition.  So too is ${\PAut(A) = \set{a\in\PEnd(A)}{a\text{ is injective}}}$, and the latter is inverse.  The automorphism group $\Aut(A) = \set{a\in\End(A)}{a\text{ is bijective}}$ is the group of units of all three of $\PEnd(A)$, $\End(A)$ and $\PAut(A)$.

For a subset $X$ of an algebra $A$ we write $\la X\ra$ for the subalgebra generated by $X$.  We say~$X$ is  \emph{independent} if $x\not\in\la X\sm \{x\}\ra$ for any $x\in X$.  A \emph{basis} of $A$ is an independent generating set.

Following \cite{Gould1995}, an  \emph{independence algebra} is an algebra $A$ that satisfies the  \emph{exchange property} and the  \emph{free basis property}:
\ben
\item[\textsf{(E)}] For all $X\sub A$ and $x,y\in A$, if $x\in\la X\cup \{y\}\ra \sm \la X\ra$, then $y\in\la X\cup \{x\}\ra$.
\item[\textsf{(F)}] If $X$ is a basis for $A$, then any map $X\to A$ can be extended (necessarily uniquely) to an endomorphism of $A$.
\een
It can be shown that any independent set of an independence algebra $A$ can be extended to a basis, and that all bases have the same size, called the  \emph{dimension} of $A$ and denoted $\dim(A)$.  It follows that \textsf{(F)} is equivalent to the ostensibly stronger condition:
\ben
\item[\textsf{(F)}$'$] If $X\sub A$ is independent, then any map $X\to A$ can be extended (necessarily uniquely) to a partial endomorphism $\la X\ra\to A$.
\een
The  \emph{rank} of a partial endomorphism $a\in\PEnd(A)$ is defined by $\rank(a) = \dim(\im(a))$, the dimension of the image of $a$.  (Note that any subalgebra of an independence algebra is itself an independence algebra, and therefore has a well-defined dimension.)

In our results we will need to assume that the underlying independence algebras are  \emph{strong}.  
This is a concept  introduced by Fountain and Lewin in their work \cite{FL1993}  on endomorphism monoids of independence algebras.
It formalises a well-known property of vector spaces, which turns out not to follow from the defining properties of an independence algebra. Fountain and Lewin give several equivalent formulations, of which we will work with the following
(see \cite[Lemma 1.6]{FL1993}):
\ben
\item[(S)] If $B,C\leq A$ are subalgebras, and if a basis $X$ for $B\cap C$ is extended to bases $X\sqcup Y$ and~$X\sqcup Z$ for $B$ and $C$, respectively, then $X\sqcup Y\sqcup Z$ is a basis for $B\vee C$.
\een
Here $\sqcup$ denotes disjoint union, and $B\vee C = \la B\cup C\ra$ is the join of $B$ and $C$ in the lattice of subalgebras of $A$.

Archetypal examples of strong independence algebras include vector spaces and (plain) sets.  (Note that the former have a separate unary operation for scalar multiplication by each element of the underlying field.)  If $A$ is simply a set, then
\[
\PEnd(A) = \PT_A \COMMA \End(A) = \T_A \COMMA \PAut(A) = \I_A \AND \Aut(A) = \S_A.
\]
If $A$ is a vector space over a field $F$, of finite dimension $n$, then $\End(A)$ is isomorphic to $M_n(F)$, the multiplicative monoid of all $n\times n$ matrices over $F$.  Structural similarities between the monoids $\T_n$ and $M_n(F)$ were one source of interest in  independence algebras in \cite{Gould1995}.  During the course of our investigations, we will uncover an intriguing point of difference between these two monoids, namely that a certain `rank-based twisting' is `tight' for~$M_n(F)$ but not for $\T_n$; see Proposition \ref{prop:MnF} and Remark \ref{rem:EndA}.

\section{Twistings: definitions and basic properties}\label{sect:basic}

Throughout this section we fix the following:
\begin{itemize}
\item
a monoid $S$, written multiplicatively, with identity $1$; and 
\item
a commutative monoid $M$, written additively, with identity $0$.
\end{itemize}

\begin{defn}\label{defn:Phi}
A \emph{twisting} of $S$ is a map $\Phi:S\times S\to \N$ satisfying the following:
\begin{enumerate}[label=\textup{\textsf{(T\arabic*)}},leftmargin=10mm]
\item \label{T1} $\Phi(a,b) + \Phi(ab,c) = \Phi(a,bc) + \Phi(b,c)$ for all $a,b,c\in S$.
\end{enumerate}
With a slight abuse of notation we write $\Phi(a,b,c)$ for the common value in \ref{T1}.  
\end{defn}

In the work on congruences of twisted partition monoids \cite{ER2022c,ER2022b} it transpired that the following property of multiplication in $\P_n$ was important in a number of technical arguments: for any partitions $a,b\in\P_n$ the product $ab$ can be obtained without introducing any floating components in the product graph by replacing $b$ by a suitable $b'\in\P_n$, and, dually, replacing $a$ by a suitable $a'\in\P_n$. Formalising to our general framework leads to the following:

\begin{defn}
A twisting $\Phi$ is \emph{tight} if it  satisfies the following:
\begin{enumerate}[label=\textup{\textsf{(T\arabic*)}},leftmargin=10mm]\addtocounter{enumi}{1}
\item \label{T2} 
\begin{enumerate}[label=\textup{\textsf{(\alph*)}},leftmargin=10mm]
\item \label{T2a} For all $a,b\in S$ there exist $a'\in S$ such that $ab = a'b$ and $\Phi(a',b) = 0$.
\item \label{T2b} For all $a,b\in S$ there exist $b'\in S$ such that $ab = ab'$ and $\Phi(a,b') = 0$.
\end{enumerate}
\end{enumerate}
Otherwise $\Phi$ is \emph{loose}.
\end{defn}

One can of course obtain a (tight) twisting by defining $\Phi(a,b) = 0$ for all $a,b\in S$.  We call this the \emph{trivial} twisting.  
Given a twisting $\Phi$ and a natural number $k$, we can define a new twisting~$\Phi'$ by $\Phi'(a,b) = k+\Phi(a,b)$; if $k>0$ then this is loose.
If $T$ is a subsemigroup of $S$, then any twisting $\Phi$ of $S$ restricts to a twisting of $T$.  If $\Phi$ is tight, then the restricted twisting need not be tight in general; similarly, a loose twisting could restrict to a tight one.
Several examples of twistings will be considered in Section \ref{sect:eg}.
In particular, in 	Remark \ref{rem:PEndA} we will see an example where \ref{T2}\ref{T2b} holds, but \ref{T2}\ref{T2a} does not.

Here is the key definition, using a twisting to combine $M$ and $S$ into a single semigroup:

\begin{defn}\label{defn:T}
Let $\Phi:S\times S\to\N$ be a twisting, and fix some $q\in M$.  The \emph{twisted product} $M\times_\Phi^qS$ is defined to be the semigroup with underlying set $M\times S$, and operation
\[
(i,a)(j,b) = (i+j+\Phi(a,b)q,ab) \qquad\text{for $i,j\in M$ and $a,b\in S$.}
\]
The product $M\times_\Phi^qS$ is said to be \emph{tight} if $\Phi$ is tight, or \emph{loose} otherwise.
\end{defn}

Associativity of the above operation follows quickly from \ref{T1}, and in fact we have
\[
(i,a)(j,b)(k,c) = (i+j+k+\Phi(a,b,c)q,abc) \qquad\text{for $i,j,k\in M$ and $a,b,c\in S$.}
\]
Note that if $\Phi$ is the trivial twisting, or if $q=0$, then $M\times_\Phi^qS$ is the classical direct product~$M\times S$, with operation $(i,a)(j,b) = (i+j,ab)$.

The following lemma gathers some basic consequences of \ref{T1} and \ref{T2}.

\begin{lemma}\label{lem:T}
If $\Phi:S\times S\to\N$ is a tight twisting, then for all $a,b,c\in S$, the following hold:
\begin{enumerate}[label=\textup{\textsf{(T\arabic*)}},leftmargin=10mm]\addtocounter{enumi}{2}
\item \label{T3} $a\leqL b \implies \Phi(a,c) \geq \Phi(b,c)$, and $a\L b \implies \Phi(a,c) = \Phi(b,c)$,
\item \label{T4} $a\leqR b \implies \Phi(c,a) \geq \Phi(c,b)$, and $a\R b \implies \Phi(c,a) = \Phi(c,b)$, 
\item \label{T5} $\Phi(1,a) = \Phi(a,1) = 0$,
\item \label{T6} there exist $a',c'\in S$ such that $abc=a'bc'$ and $\Phi(a',b,c')=0$.
\end{enumerate}
\end{lemma}

\pf
\firstpfitem{\ref{T3}}  For the first implication, suppose $a\leqL b$.  Then~$a=sb$ for some~$s\in S$, and by~\ref{T2} we can assume that $\Phi(s,b)=0$.  Combining this with \ref{T1}, we have
\[
\Phi(a,c) = 0 + \Phi(a,c) = \Phi(s,b)+\Phi(sb,c) = \Phi(s,bc)+\Phi(b,c) \geq \Phi(b,c).
\]
The second implication follows from the first.

\pfitem{\ref{T4}}  This is dual to \ref{T3}.

\pfitem{\ref{T5}}  By \ref{T2}, there exists $a'\in M$ such that $a\cdot1 = a'\cdot 1$ and $\Phi(a',1) = 0$.  But $a\cdot1 = a'\cdot 1$ simply says that $a=a'$, so indeed $\Phi(a,1) = \Phi(a',1) = 0$.  We obtain $\Phi(1,a)=0$ by symmetry.

\pfitem{\ref{T6}}  By \ref{T2} there exist $a',c'\in S$ such that $ab=a'b$ and $a'b\cdot c = a'b\cdot c'$, with $\Phi(a',b)=0$ and $\Phi(a'b,c')=0$.  We then have
\[
a'bc' = a'bc = abc \AND \Phi(a',b,c') = \Phi(a',b)+\Phi(a'b,c') = 0+0 = 0.  \qedhere
\]
\epf

We now list some fundamental properties of a tight twisted product $T = M\times_\Phi^qS$, which we will use without explicit reference.  The third is easily checked, and the rest are simple consequences of \ref{T5}.
\begin{enumerate}[label=\textup{\textsf{(P\arabic*)}},leftmargin=10mm]
\item $T$ is a monoid with identity $(0,1)$.
\item The set $M\times\{1\}$ is a submonoid of $T$, and is isomorphic to $M$.
\item If $M$ has an absorbing element $\infty$ (so $i+\infty=\infty$ for all $i\in M$), then $\{\infty\}\times S$ is a submonoid of $T$, and is isomorphic to $S$.  (This holds more generally if $\infty$ is replaced by any element $w\in M$ for which $2w = w = w+q$.)
\item $T$ is generated by the union of $M\times\{1\}$ and $\{0\}\times S$.  Specifically, we have
\[
(i,a) = (i,1)(0,a) = (0,a)(i,1) \qquad\text{for any $i\in M$ and $a\in S$.}
\]
\end{enumerate}
While the subset $\{0\}\times S$ of $T$ is in one-one correspondence with $S$, it is typically not a submonoid, as $(0,a)(0,b)=(\Phi(a,b)q,ab)$ for $a,b\in S$.

The next result will be used when we consider idempotents in Section \ref{sect:ET}.

\begin{lemma}\label{lem:ae}
If $\Phi$ is a twisting of a monoid $S$, and if $a\in S$ and $e\in E(S)$, then
\[
a\leqL e \implies \Phi(a,e)=\Phi(e,e)
\AND
a\leqR e \implies \Phi(e,a)=\Phi(e,e).
\]
\end{lemma}

\pf
For the first implication (the second is dual), note that $a\leqL e$ implies $a=ae$.  Combining this with \ref{T1} yields
\[
\Phi(a,e) + \Phi(a,e) = \Phi(a,e) + \Phi(ae,e) = \Phi(a,ee) + \Phi(e,e) = \Phi(a,e) + \Phi(e,e),
\]
and the assertion follows.
\epf

Some of our motivating examples involve semigroups $S$ with an involution ${}^*$.  We say that a twisting $\Phi$ of such a semigroup is \emph{$*$-symmetric} if
\begin{equation}\label{eq:SS}
\Phi(a,b) = \Phi(b^*,a^*) \qquad\text{for all $a,b\in S$.}
\end{equation}

\begin{lemma}\label{lem:T2}
If $\Phi$ is a $*$-symmetric twisting of a monoid $S$ with involution, then \ref{T2}\ref{T2a} and \ref{T2}\ref{T2b} are equivalent.  Consequently, $\Phi$ is tight if and only if either of these holds.
\end{lemma}

\pf
It suffices by symmetry to show that \ref{T2}\ref{T2a} implies \ref{T2}\ref{T2b}, so suppose \ref{T2}\ref{T2a} holds.  To verify~\ref{T2}\ref{T2b}, let~$a,b\in S$.  By assumption there exists $c\in S$ such that $b^*a^* = ca^*$ and $\Phi(c,a^*) = 0$.  We then take $b' = c^*$, and we note that
\[
ab' = ac^* = (ca^*)^* = (b^*a^*)^* = ab \AND \Phi(a,b') = \Phi(a,c^*) = \Phi(c,a^*) = 0.  \qedhere
\]
\epf

\section{Motivating examples}\label{sect:eg}

The purpose of this section is to describe some natural examples of twistings.  These include the canonical twistings of diagram monoids (Section \ref{subsect:PhiPn}), a new twisting for monoids of partial endomorphisms of strong independence algebras (Section \ref{subsect:PEndA}), including full transformation monoids and matrix monoids, and analogous new twistings for diagram monoids (Section \ref{subsect:newPhiPn}).  

\subsection{Diagram monoids and their canonical twistings}\label{subsect:PhiPn}

As just noted, our first motivating example concerns diagram monoids, and we begin with the partition monoid~$\P_n$.
For partitions $a,b\in\P_n$, we let $\Phi(a,b)$ be the number of floating components in the product graph $\Pi(a,b)$.  It is a non-trivial fact that this $\Phi$ is a tight twisting.  Properties~\ref{T1} and~\ref{T2} are \cite[Lemma 4.1]{FL2011} and \cite[Lemma 2.9(i)]{ER2022b}, respectively.  We call this~$\Phi$ the \emph{canonical twisting} of $\P_n$.  It is clear from the definition that $\Phi$ is $*$-symmetric with respect to the standard involution ${}^*$ of $\P_n$.  

Taking $M$ to be either $(\N,+)$ or $(\Z,+)$, and $q=1$ in either case, leads to two twisted products:
$\N\times_\Phi^1\P_n$ and $\Z\times_\Phi^1\P_n$.
The first of these monoids was the subject of the papers \cite{ER2022c,ER2022b}.  
The second has not been explicitly studied, to the best of our knowledge, but will play a pivotal role in the forthcoming paper~\cite{EGPAR2025}.

Of course $\Phi$ restricts to a twisting of any submonoid of $\P_n$.  Curiously, this restriction is not always tight, as we will see in Proposition \ref{prop:canPhi}.  For the proof we require a technical lemma about planar partitions:

\begin{lemma}\label{lem:PPn}
For any $a\in\PP_n$ with $\rank(a)\geq1$, there exists $c\in\PP_n$ such that $ac\R a$, $\codom(ac)=\bn$ and $\Phi(a,c)=0$.
\end{lemma}

\pf
Let $r=\rank(a)$.
If $r=1$ we can simply take $c = \binom\bn\bn$, the partition with a single block.  Now suppose $r\geq2$, and write $a = \begin{partn}{6}A_1&\cdots&A_r&C_1&\cdots&C_s \\ \hhline{~|~|~|-|-|-} B_1&\cdots&B_r&D_1&\cdots&D_t\end{partn}$.  By \cite[Lemma 7.1]{EMRT2018} we can assume that ${B_1<\cdots<B_r}$ (where $X<Y$ means that $x<y$ for all $x\in X$ and $y\in Y$).  Write $m_i = \max(B_i)$ for each~$i$, and define the sets
\[
M_1 = \{1,\ldots,m_1\} \COMMa M_i = \{m_{i-1}+1,\ldots,m_i\} \text{ for $1<i<r$} \ANd M_r = \{m_{r-1}+1,\ldots,n\}.
\]
Then the conditions are easily checked for $c = \begin{partn}{3}M_1&\cdots&M_r\\M_1&\cdots&M_r\end{partn}$, noting that $ac =  \begin{partn}{6}A_1&\cdots&A_r&C_1&\cdots&C_s \\ \hhline{~|~|~|-|-|-} M_1&\cdots&M_r&\multicolumn{3}{c}{}\end{partn}$.  
\epf

\begin{prop}\label{prop:canPhi}
The canonical twisting is tight for $\P_n$, $\B_n$, $\PP_n$ and $\TL_n$, but loose for~$\PB_n$ and~$\M_n$ if $n\geq2$.
\end{prop}

\pf
Property \ref{T2} for $\B_n$ and $\TL_n$ is established in \cite{EF2025}.  

For $\PP_n$, it suffices by Lemma \ref{lem:T2} to establish \ref{T2}\ref{T2b}.  To do so, let ${a,b\in\PP_n}$.  Suppose first that $\rank(a)=0$, and let the lower non-transversals of $ab$ be $A_1',\ldots,A_k'$, where we assume that $1\in A_1$.  Then $ab = ab'$ for $b' = \begin{partn}{4}\bn&\multicolumn{3}{c}{} \\ \hhline{~|-|-|-} A_1&A_2&\cdots&A_k\end{partn}\in\PP_n$, and we have $\Phi(a,b')=0$ because $\dom(b')=\bn$.
Now suppose $\rank(a)\geq1$, and let $c\in\PP_n$ be as in Lemma \ref{lem:PPn}.  Now, ${ab \leqR a \R ac}$, so we have $ab = acd$ for some $d\in\PP_n$.  Thus, it remains to check that $\Phi(a,b')=0$ for $b' = cd$.  Now, 
\[
\Phi(a,c,d) = \Phi(a,c) + \Phi(ac,d) = 0,
\]
as $\Phi(a,c)=0$, and $\Phi(ac,d)=0$ since $\codom(ac)=\bn$.  But then it follows that
\[
\Phi(a,b') = \Phi(a,cd) \leq \Phi(a,cd)+\Phi(c,d) = \Phi(a,c,d) = 0,
\]
so certainly $\Phi(a,b')=0$.

For the failure of \ref{T2}\ref{T2b} in $\PB_n$ and $\M_n$, take
\begin{equation}\label{eq:ab}
a=\custpartn{1,2,3,6}{1,2,3,6}{\uarc12\stline33\stline66\udotted36\ldotted36}
\AND
b=\custpartn{1,2,3,6}{1,2,3,6}{\uarc12\darc12\stline33\stline66\udotted36\ldotted36},
\end{equation}
both from $\M_n$.  Then $ab=b$, and the only other $b'\in\PB_n$ satisfying $ab'=b$ is $b'=\custpartn{1,2,3,6}{1,2,3,6}{\darc12\stline33\stline66\udotted36\ldotted36}$, but we have $\Phi(a,b)=1$ and $\Phi(a,b')=2$.
\epf

The products $\N\times_\Phi^1\B_n$ and $\N\times_\Phi^1\TL_n$ have appeared frequently in the literature \cite{BDP2002,DE2018,KV2023,DE2017,KV2019,CHKLV2019,ACHLV2015,LF2006}, with the latter typically called the \emph{Kauffman monoid}.  The larger monoids $\Z\times_\Phi^1\B_n$ and $\Z\times_\Phi^1\TL_n$ have also featured in \cite{KV2019,KV2023}.  Twisted products involving $\PB_n$ and $\M_n$ have more intricate structures, due to their twistings being loose.  For example, the middle diagrams in Figures \ref{fig:P2Phi} and \ref{fig:PB2Phi} show twisted products involving $\P_2$ and $\PB_2$, respectively, with respect to the canonical twistings; the former is tight, and the latter is loose.

\subsection{Rigid twistings}\label{subsect:rigid}

The remaining twistings we consider here arise from a general construction:

\begin{lemma}\label{lem:Phirm}
Let $S$ be a monoid, $r:S\to\Z$ an arbitrary function, and $m\in\Z$ an arbitrary integer.  Then the function
\[
\Phi_{r,m}:S\times S\to\Z \GIVENBY \Phi_{r,m}(a,b) = m - r(a) - r(b) + r(ab) \qquad\text{for $a,b\in S$}
\]
satisfies \ref{T1}.  Consequently, $\Phi_{r,m}$ is a twisting of $S$ if and only if the following inequality holds:
\begin{equation}\label{eq:Sylvester}
r(a) + r(b) \leq r(ab) + m \qquad\text{for all $a,b\in S$.}
\end{equation}
If additionally $S$ has an involution ${}^*$ for which the identity $r(a^*)=r(a)$ holds, then $\Phi_{r,m}$ is $*$-symmetric.
\end{lemma}

\pf
It is easy to check that both sides of \ref{T1} evaluate to $2m-r(a)-r(b)-r(c)+r(abc)$.
It is also clear that the inequality \eqref{eq:Sylvester} is equivalent to having $\Phi_{r,m}(a,b)\geq0$ for all $a,b\in S$.
Finally, if the identity $r(a^*)=r(a)$ holds, then for any $a,b\in S$ we have
\[
\Phi_{r,m}(b^*,a^*) = m - r(b^*) - r(a^*) + r(b^*a^*) = m - r(b) - r(a) + r(ab) = \Phi_{r,m}(a,b).  \qedhere
\]
\epf

A twisting of the above form $\Phi_{r,m}$ will be called \emph{rigid}.  Note that
\[
\Phi_{r,m}(a,1) = \Phi_{r,m}(1,a) = m - r(1) \qquad\text{for any $a\in S$.}
\]
In particular, if $r(1)=m$, then $M\times_{\Phi_{r,m}}^qS$ is a monoid with identity $(0,1)$, whether $\Phi_{r,m}$ is tight or not.  Note also that $r(1)=m$ is a necessary condition for a rigid twisting $\Phi_{r,m}$ to be tight.

Rigid twistings can be thought of as generalisations of the \emph{coboundaries} from the paper \cite{DE1995}, which studied analogues of twistings that map from a cancellative commutative monoid to a group (rather than from an arbitrary monoid to $\N$).  The main result of \cite{DE1995} was that every twisting of the type they study is a coboundary.
This is not the case for twistings in general.  In particular, we will see later that the canonical twistings on diagram monoids are not rigid; see Remark \ref{rem:rigidPn}.

\subsection{(Partial) endomorphism monoids of strong independence algebras}\label{subsect:PEndA}

A special case of \eqref{eq:Sylvester} is the well-known \emph{Sylvester rank inequality} in linear algebra, which says that for $n\times n$ matrices~$a$ and $b$ over a field $F$ we have
\[
\rank(a) + \rank(b) \leq \rank(ab) + n.
\]
Because of Lemma \ref{lem:Phirm}, this leads to a rigid twisting on the multiplicative monoid $M_n(F)$ of all such matrices, given by $\Phi(a,b) = n - \rank(a) - \rank(b) + \rank(ab)$.  Among other things, we will see that this twisting is tight; see Proposition \ref{prop:MnF}.

In fact, Sylvester's rank inequality holds in much greater generality, namely for partial endomorphisms of finite-dimensional strong independence algebras:

\begin{prop}\label{prop:PEndA}
If $A$ is a strong independence algebra of finite dimension $n$, then 
\[
\rank(a) + \rank(b) \leq \rank(ab) + n \qquad\text{for all $a,b\in\PEnd(A)$.}
\]
\end{prop}

\pf
Write $k = \rank(a)$ and $l = \rank(ab)$.  We must show that
\begin{equation}\label{eq:l+n-k}
\rank(b) \leq l+n-k.
\end{equation}
Also define the subalgebras
\[
B = \im(a) \AND C = \dom(b) .
\]
Choose a basis $\{w_1,\ldots,w_l\}$ for $\im(ab)$, and choose $u_1,\ldots,u_l\in A$ such that $u_iab = w_i$ for each~$i$.  Also set $v_i = u_ia$ for each~$i$.  Since $\{w_1,\ldots,w_l\}$ is independent, so too are $\{u_1,\ldots,u_l\}$ and $\{v_1,\ldots,v_l\}$, and we note that $\{u_1,\ldots,u_l\}\sub \dom(a)$ and $\{v_1,\ldots,v_l\}\sub \im(a)\cap\dom(b) = B\cap C$.  Next we extend $\{v_1,\ldots,v_l\}$ to bases:
\bit
\item $\{v_1,\ldots,v_m\} = \{v_1,\ldots,v_l\}\cup\{v_{l+1},\ldots,v_m\}$ of $B\cap C$,
\item $\{v_1,\ldots,v_k\} = \{v_1,\ldots,v_m\}\cup\{v_{m+1},\ldots,v_k\}$ of $B$, and
\item $\{v_1,\ldots,v_m\}\cup\{x_1,\ldots,x_h\}$ of $C = \dom(b)$.
\eit
Now, $\im(b)$ is generated by the image of this last basis, meaning that
\[
\im(b) = \la v_1b,\ldots,v_mb,x_1b,\ldots,x_hb\ra = \la v_1b,\ldots,v_mb\ra \vee \la x_1b,\ldots,x_hb\ra.
\]
It follows that
\[
\rank(b) = \dim(\im(b)) \leq \dim\la v_1b,\ldots,v_mb\ra +\dim \la x_1b,\ldots,x_hb\ra \leq \dim\la v_1b,\ldots,v_mb\ra + h.
\]
Since $v_i\in B\cap C = \im(a)\cap\dom(b)$ for each $1\leq i\leq m$, we have
\[
\la v_1b,\ldots,v_mb\ra \sub \im(ab) = \la w_1,\ldots,w_l\ra = \la v_1b,\ldots,v_lb\ra \sub \la v_1b,\ldots,v_mb\ra.
\]
It follows that $\la v_1b,\ldots,v_mb\ra = \im(ab)$, and so $\dim\la v_1b,\ldots,v_mb\ra = \rank(ab) = l$.  Thus, continuing from above, we have
\[
\rank(b) \leq \dim\la v_1b,\ldots,v_mb\ra + h = l +h.
\]
Thus, we can complete the proof of \eqref{eq:l+n-k}, and hence of the proposition, by showing that $h\leq n-k$.  But this follows from the strong property in $A$.  Indeed, looking at the above bases for $B$, $C$ and~$B\cap C$, it follows that
\[
\{v_1,\ldots,v_m\} \cup \{v_{m+1},\ldots,v_k\} \cup \{x_1,\ldots,x_h\} = \{v_1,\ldots,v_k\} \cup \{x_1,\ldots,x_h\}
\]
is a basis of $B\vee C$.  In particular, this set is independent, and hence $k+h\leq n$, i.e.~$h\leq n-k$, as required.
\epf

The inequality just established leads to a rigid twisting of $\PEnd(A)$.  Although this need not be tight in general, it does satisfy one half of axiom \ref{T2}.

\begin{thm}\label{thm:PEndA}
Suppose $A$ is a strong independence algebra of finite dimension $n$.
\ben
\item \label{PEndA1} The partial endomorphism monoid $\PEnd(A)$ has a twisting given by
\[
\Phi(a,b) = n - \rank(a) - \rank(b) + \rank(ab).  
\]
\item \label{PEndA2} This twisting satisfies \ref{T2}\ref{T2b}.  That is, for any $a,b\in\PEnd(A)$ we have $ab=ab'$ for some $b'\in\PEnd(A)$ with $\Phi(a,b') = 0$.  
\item \label{PEndA3} In the previous part, if $b$ belongs to $\End(A)$ or to $\PAut(A)$, then such an element $b'$ exists in $\End(A)$ or $\PAut(A)$, respectively.
\een
\end{thm}

\pf
\firstpfitem{\ref{PEndA1}}  This follows from Lemma \ref{lem:Phirm} and Proposition \ref{prop:PEndA}.

\pfitem{\ref{PEndA2}}  Fix $a,b\in\PEnd(A)$, and let $k,l\in\N$ and $B,C\leq A$ be as in the proof of Proposition \ref{prop:PEndA}.  Also fix the bases
\[
\text{$\{w_1,\ldots,w_l\}$ of $\im(ab)$\COMMA
$\{v_1,\ldots,v_m\}$ of $B\cap C$ \AND
$\{v_1,\ldots,v_k\}$ of $B$,}
\]
as in the same proof, where $l\leq m\leq k$, and where $v_ib=w_i$ for each $1\leq i\leq l$.  Now extend $\{v_1,\ldots,v_k\}$ to a basis $\{v_1,\ldots,v_n\}$ of $A$, and extend $\{w_1,\ldots,w_l\}$ to an independent set $\{w_1,\ldots,w_l\} \cup \{z_{k+1},\ldots,z_n\}$.  (The latter is possible because $l\leq k$.)  Let
\[
D = \la v_1,\ldots,v_m,v_{k+1},\ldots,v_n\ra
\AND
E = \la w_1,\ldots,w_l,z_{k+1},\ldots,z_n\ra,
\]
and define $b':D\to E$ in $\PEnd(A)$ by
\[
v_ib' = v_ib \text{ for $1\leq i\leq m$} \AND v_jb' = z_j \text{ for $k+1\leq j\leq n$.}
\]
Note that $b'$ agrees with $b$ on the basis $\{v_1,\ldots,v_k\} = \{v_1,\ldots,v_m\} \cup \{v_{m+1},\ldots,v_k\}$ of $B = \im(a)$, in the sense that $v_ib$ and $v_ib'$ are equal for $1\leq i\leq m$, but are both undefined for $m+1\leq i\leq k$.  It follows from this that $ab=ab'$.  By construction, we have
\begin{align*}
\rank(b') = \dim(\im(b')) &= \dim \la v_1b,\ldots,v_mb,z_{k+1},\ldots,z_n\ra \\
&= \dim \la w_1,\ldots,w_l,z_{k+1},\ldots,z_n\ra = l+n-k,
\end{align*}
where we used the fact that 
\[
\la w_1,\ldots,w_l\ra = \im(ab) = (B\cap C)b = \la v_1,\ldots,v_m\ra b = \la v_1b,\ldots,v_mb\ra.
\]
It follows that indeed
\[
\Phi(a,b') = n - k - (l+n-k) + l = 0.
\]

\pfitem{\ref{PEndA3}}  Suppose first that $b\in\PAut(A)$.  We must then have $l=m$, as $b$ maps $B\cap C = \la v_1,\ldots,v_m\ra$ bijectively onto $\la w_1,\ldots,w_l\ra$.  It follows from this that
\[
\dim(D) = m+n-k = l+n-k = \dim(E),
\]
so that $b':D\to E$ is indeed an isomorphism, and hence belongs to $\PAut(A)$.

On the other hand, if $b\in\End(A)$, then $C = \dom(b) = A$, and so $B\cap C = B$, which forces $m=k$.  We then have
\[
\dim(D) = m+n-k = k+n-k = n,
\]
so that $\dom(b') = D = A$, and $b'\in\End(A)$.
\epf

\begin{rem}
The assumption that the independence algebra $A$ is \emph{strong} was used in the proof of Proposition \ref{prop:PEndA}, which itself feeds into the proof of part \ref{PEndA1} of Theorem \ref{thm:PEndA}.  However, this assumption on $A$ was not used in the proofs of parts \ref{PEndA2} or \ref{PEndA3} of the theorem.
\end{rem}

\begin{rem}\label{rem:PEndA}
To see that \ref{T2}\ref{T2a} does not hold (in general) for the twisting $\Phi$ of $\PEnd(A)$ from Theorem \ref{thm:PEndA}, consider the case that $A=\{1,2,3\}$ with no operations, so that $\PEnd(A) = \PT_3$ is the partial transformation semigroup of degree~$3$.  Then with
\[
a = \trans{1&2&3\\1&1&3} \AND b = \trans{1&2&3\\1&2&2},
\]
both from $\PT_3$, we have $ab = \trans{1&2&3\\1&1&2}$.  The only other element $a'$ of $\PT_3$ satisfying $ab=a'b$ is $a'=\trans{1&2&3\\1&1&2} (=ab)$, and we have $\Phi(a,b) = \Phi(a',b) = 1$.
\end{rem}

The twisting of $\PEnd(A)$ from Theorem \ref{thm:PEndA} restricts to a twisting on any submonoid of $\PEnd(A)$.  The next result concerns two of the most important such submonoids, namely $\End(A)$ and $\PAut(A)$.  Since the latter is inverse, the inversion map $a\mt a^{-1}$ is an involution.

\begin{prop}\label{prop:EndA}
If $A$ is a strong independence algebra of finite dimension $n$, then $\End(A)$ and $\PAut(A)$ have twistings given by
\[
\Phi(a,b) = n - \rank(a) - \rank(b) + \rank(ab) .  
\]
This is tight and ${}^{-1}$-symmetric for $\PAut(A)$.
\end{prop}

\pf
It remains only to prove the assertion concerning $\PAut(A)$.  Symmetry follows from Lemma \ref{lem:Phirm} and the identity $\rank(a)=\rank(a^{-1})$.  Tightness then follows from Lemma \ref{lem:T2} and Theorem \ref{thm:PEndA}\ref{PEndA3}.
\epf

\begin{rem}\label{rem:EndA}
For $A=\{1,2,3\}$ with no operations, the failure of \ref{T2}\ref{T2a} for $\End(A)=\T_3$ can be deduced from Remark \ref{rem:PEndA}, as the elements $a,b\in\PT_3=\PEnd(A)$ used there actually belong to~$\T_3$.
\end{rem}

As special cases, consider again the algebra $A = \{1,\ldots,n\}$ with no operations.  Here we have
\[
\PEnd(A) = \PT_n \COMMA \End(A) = \T_n \AND \PAut(A) = \I_n.
\]
It follows from Proposition \ref{prop:EndA} that the twisting on $\I_n$ is tight, but we observed in Remarks~\ref{rem:PEndA} and~\ref{rem:EndA} that those on $\PT_n$ and $\T_n$ are loose.  The next result gives an important special case in which the twisting on $\End(A)$ \emph{is} tight.
It is well known that the linear monoid $M_n(F)$, consisting of all $n\times n$ matrices over an arbitrary field $F$, under multiplication, is isomorphic to $\End(V)$ for any $n$-dimensional vector space $V$ over $F$, and that $V$ is a strong independence algebra (of dimension $n$).  It is also well known that the transpose map $a\mt a^T$ is an involution on $M_n(F)$.

\begin{prop}\label{prop:MnF}
For any $n\geq1$ and any field $F$, the linear monoid $M_n(F)$ has a tight,~${}^T$-symmetric twisting, given by
\[
\Phi(a,b) = n - \rank(a) - \rank(b) + \rank(ab).
\]
\end{prop}

\pf
By Proposition \ref{prop:EndA}, and the above-mentioned isomorphism $M_n(F)\cong\End(V)$, this $\Phi$ is a twisting.
Since $\rank(a^T)=\rank(a)$ for all $a\in M_n(F)$, we again obtain ${}^T$-symmetry from Lemma~\ref{lem:Phirm}.  Thus, another appeal to Lemma \ref{lem:T2} and Theorem \ref{thm:PEndA} shows that $\Phi$ is tight.
\epf

\subsection{New twistings for diagram monoids}\label{subsect:newPhiPn}

We now prove a Sylvester-style rank inequality for partition monoids, which will allow us to introduce a new family of (rigid) twistings for diagram monoids; see Theorem \ref{thm:altPhiPn}.  Curiously, these are tight for the Brauer monoids $\B_n$, but not for any of $\P_n$, $\PP_n$, $\PB_n$, $\M_n$ or $\TL_n$ (apart from trivially small $n$), as we show in Proposition \ref{prop:altPhiBn}.  Before we begin, we need a technical result concerning joins of equivalence relations.

Let $\Eq(X)$ be the $\vee$-semilattice of equivalence relations on a finite set $X$.  
For $\ve\in\Eq(X)$ we write $\Vert\ve\Vert = |X/\ve|$ for the number of $\ve$-classes.  
For distinct $x,y\in X$, let $\ve_{xy}$ be the equivalence whose only non-trivial class is $\{x,y\}$.  Note then that for any $\ve\in\Eq(X)$, the join $\ve\vee\ve_{xy}$ is either equal to $\ve$ or else is obtained by merging precisely two $\ve$-classes.  Thus, we have
\begin{equation}\label{eq:veveeve}
\Vert\ve\vee\ve_{xy}\Vert \geq \Vert\ve\Vert-1.
\end{equation}

\begin{lemma}\label{lem:SylvesterEqX}
If $\ve,\eta\in\Eq(X)$ for a finite set $X$, then $\Vert\ve\Vert+\Vert\eta\Vert \leq \Vert\ve\vee\eta\Vert + |X|$.
\end{lemma}

\pf
Write $k = \Vert\ve\Vert$, $l = \Vert\eta\Vert$ and $n=|X|$.  We show by descending induction on $l$ that $\Vert\ve\vee\eta\Vert \geq k+l-n$.  If $l=n$, then $\eta=\De_X$, and the claim is obvious since $\ve\vee\eta=\ve$.  So now suppose $l<n$.  We can then write $\eta = \eta' \vee \ve_{xy}$ for some $\eta'\in\Eq(X)$ with $\Vert\eta'\Vert = l+1$, and for distinct $x,y\in X$.  Using \eqref{eq:veveeve} and induction, we then have
\[
\Vert\ve\vee\eta\Vert = \Vert\ve\vee\eta'\vee\ve_{xy}\Vert \geq \Vert\ve\vee\eta'\Vert - 1 \geq (k+(l+1)-n) - 1 = k+l-n.  \qedhere
\]
\epf

\begin{rem}
Lemma \ref{lem:SylvesterEqX} says that when $|X|=n$ is finite, the semilattice $\Eq(X)$ satisfies the Sylvester-type inequality \eqref{eq:Sylvester} with $m=n$ and $r(\ve)=\Vert\ve\Vert$.  It then follows from Lemma \ref{lem:Phirm} that~$\Eq(X)$ has a (rigid) twisting given by
\[
\Phi(\ve,\eta) = n - \Vert\ve\Vert - \Vert\eta\Vert + \Vert\ve\vee\eta\Vert.
\]
This is easily seen to be tight.  Indeed, let $\ve,\eta\in\Eq(X)$, and let $k = \Vert\ve\Vert$ and $l=\Vert\ve\vee\eta\Vert$.  Also let $\bn/(\ve\vee\eta) = \{A_1,\ldots,A_l\}$.  Now, each $A_i$ is a union of $\ve$-classes, say $A_i = A_{i1}\cup\cdots\cup A_{ik_i}$.
For each $1\leq i\leq l$ and $1\leq j\leq k_i$, choose some $a_{ij}\in A_{ij}$, and let $\eta'\in\Eq(X)$ be such that $\set{a_{ij}}{1\leq j\leq k_i}$ is an $\eta'$-class for each $1\leq i\leq l$, with every other $\eta'$-class a singleton.  Then
\[
\ve\vee\eta' = \ve\vee\eta \AND \Vert\eta'\Vert = l + n-(k_1+\cdots+k_l) = l+n-k,
\]
so that $\Phi(\ve,\eta') = n - k - (l+n-k) + l = 0$.
\end{rem}

We now return to the partition monoid $\P_n$:

\begin{thm}\label{thm:altPhiPn}
For any $a,b\in\P_n$ we have 
\[\rank(a) + \rank(b) \leq \rank(ab) + n.
\]
Consequently,~$\P_n$ has a $*$-symmetric twisting given by
\[
\Phi(a,b) = n - \rank(a) - \rank(b) + \rank(ab) .
\]
\end{thm}

\pf
Because of Lemma \ref{lem:Phirm} and the identity $\rank(a^*)=\rank(a)$, it suffices to prove the claimed inequality.
So fix $a,b\in\P_n$, and write $k = \rank(a)$, $l = \rank(b)$ and $m = \rank(ab)$; we must show that $m \geq k+l-n$.

During the proof we will be interested in three kinds of connected components of the product graph $\Pi(a,b)$.  Specifically, we call a component $C$:
\bit
\item a \emph{transversal} if $C\cap\bn\not=\es$ and $C\cap\bn'\not=\es$ (in which case necessarily also $C\cap\bn''\not=\es$), 
\item an \emph{upper semi-transversal} if $C\cap\bn\not=\es$, $C\cap\bn''\not=\es$ and $C\cap\bn'=\es$,
\item a \emph{lower semi-transversal} if $C\cap\bn=\es$, $C\cap\bn''\not=\es$ and $C\cap\bn'\not=\es$.
\eit
Note that for any such $C$, the intersection $C\cap\bn''$ has the form $A''$ for some $\coker(a)\vee\ker(b)$-class~$A$.  Note also that ${m = \rank(ab)}$ is equal to the number of transversals of $\Pi(a,b)$.  Denote the unions of all transversals, all upper semi-transversals and all lower semi-transversals of $\Pi(a,b)$, respectively, by:
\[
X\cup Y''\cup Z' \COMMA T\cup U'' \AND V''\cup W', \WHERE T,U,V,W,X,Y,Z\sub\bn.
\]
By construction, the sets $Y$, $U$ and $V$ are pairwise disjoint, so it follows that
\[
|U|+|V|+|Y| \leq n.
\]
Now consider the restrictions $\ve = \coker(a)\rest_Y$ and $\eta = \ker(b)\rest_Y$.  For any transversal $C$ of $\Pi(a,b)$, we have $C\cap\bn'' = A''$ for some $\ve\vee\eta$-class $A$.  It follows that
\[
m = \rank(ab) = \Vert\ve\vee\eta\Vert.
\]
Also let:
\bit
\item $k_1$ (resp.~$l_1$) be the number of transversals of $a$ (resp.~$b$) contained in a transversal of $\Pi(a,b)$, 
\item $k_2$ (resp.~$l_2$) be the number of transversals of $a$ (resp.~$b$) contained in a semi-transversal. 
\eit
We then have
\[
k = k_1+k_2 \COMMA l = l_1+l_2 \COMMA k_1\leq\Vert\ve\Vert \COMMA l_1\leq\Vert\eta\Vert \COMMA k_2\leq|U| \AND l_2\leq|V|.
\]
Combining the above information with Lemma \ref{lem:SylvesterEqX}, it follows that indeed
\begin{align*}
m = \Vert\ve\vee\eta\Vert \geq \Vert\ve\Vert + \Vert\eta\Vert - |Y| \geq k_1+l_1 - |Y| &= (k-k_2) + (l-l_2) - |Y| \\
&\geq k+l -|U|-|V|-|Y| \geq k+l-n.  \qedhere
\end{align*}
\epf

We call $\Phi$ from Theorem \ref{thm:altPhiPn} the \emph{rank-based twisting} of $\P_n$.  As ever, this restricts to a twisting on any submonoid of $\P_n$.

\begin{prop}\label{prop:altPhiBn}
The rank-based twisting is tight for $\B_n$, but loose for $\P_n$, $\PP_n$, $\PB_n$ and~$\M_n$ for $n\geq2$, and for $\TL_n$ for $n\geq3$.
\end{prop}

\pf
For the failure of \ref{T2} in $\P_n$, $\PP_n$, $\PB_n$ and $\M_n$, consider $a,b\in\M_n$ as in \eqref{eq:ab}, and note that $ab = b$.  For any $b'\in\P_n$ with $ab'=ab$, we have
\[
\Phi(a,b') = n - (n-2) - \rank(b') + (n-2) = n - \rank(b').
\]
For this to equal $0$ we would need $\rank(b') = n$, meaning that $b'\in\S_n$ is a permutation.  But for any permutation $b'\in\S_n$ we have $\coker(ab') = \De_\bn \not= \coker(b)$, so that $ab'\not=b=ab$.

For $\TL_n$, consider
\[
a = \custpartn{1,2,3,4,7}{1,2,3,4,7}{\uarc12\stline31\darc23\stline44\stline77\udotted47\ldotted47}
\AND
b = \custpartn{1,2,3,4,7}{1,2,3,4,7}{\uarc12\stline33\darc12\stline44\stline77\udotted47\ldotted47}.
\]
As in the previous paragraph, any $b'\in\TL_n$ satisfying $ab'=ab(=b)$ and $\Phi(a,b')=0$ must have $\rank(b')=n$.  But there is a unique Temperley--Lieb diagram of rank~$n$, namely the identity $1$, and we do not have $a\cdot1=b$.

By Lemma \ref{lem:T2}, it remains to verify \ref{T2}\ref{T2b} in $\B_n$.  To do so, fix~${a,b\in\B_n}$, let $k = \rank(a)$ and $l = \rank(ab)$, and write
\[
a = \begin{partn}{9}
x_1&\cdots&x_l&x_{l+1}&\cdots&x_k&x_{k+1},x_{k+2}&\cdots&x_{n-1},x_n \\ \hhline{~|~|~|~|~|~|-|-|-} y_1&\cdots&y_l&y_{l+1}&\cdots&y_k&y_{k+1},y_{k+2}&\cdots&y_{n-1},y_n
\end{partn}.
\]
Re-ordering if necessary, we can assume that $ab$ has the form
\[
ab = \begin{partn}{9}
x_1&\cdots&x_l&x_{l+1},x_{l+2}&\cdots&x_{k-1},x_k&x_{k+1},x_{k+2}&\cdots&x_{n-1},x_n \\ \hhline{~|~|~|-|-|-|-|-|-} z_1&\cdots&z_l&z_{l+1},z_{l+2}&\cdots&z_{k-1},z_k&z_{k+1},z_{k+2}&\cdots&z_{n-1},z_n
\end{partn}.
\]
We then have $ab = ab'$ for
\[
b' = 
\begin{partn}{9}
y_1&\cdots&y_l&y_{k+1}&\cdots&y_n&y_{l+1},y_{l+2}&\cdots&y_{k-1},y_k \\ \hhline{~|~|~|~|~|~|-|-|-} z_1&\cdots&z_l&z_{k+1}&\cdots&z_n&z_{l+1},z_{l+2}&\cdots&z_{k-1},z_k
\end{partn},
\]
and $\Phi(a,b') = n - k - (l+n-k) + l = 0$.
\epf

\begin{rem}\label{rem:rigidPn}
The rank-based twisting coincides with the canonical (float-counting) twisting on the submonoid $\I_n$ of $\P_n$.  Indeed, consider $a,b\in\I_n$.  A floating component in the product graph $\Pi(a,b)$ is simply a point $x''$, and there is one for each $x\in\bn\sm(\codom(a)\cup\dom(b))$.  Thus, writing $\Phi$ for the canonical twisting, we have
\begin{align*}
\Phi(a,b) &= n - |\codom(a)\cup\dom(b)| \\
&= n - |\codom(a)| - |\dom(b)| + |\codom(a)\cap\dom(b)| \\
&= n - \rank(a) - \rank(b) + \rank(ab).
\end{align*}

The rank-based and canonical twistings are distinct for the monoids $\P_n$, $\PP_n$, $\PB_n$, $\B_n$, $\M_n$ and $\TL_n$, however.  In fact, we claim that for any of these monoids, the canonical twisting~$\Phi$ is not rigid.  That is, it does not have the form $\Phi_{r,m}$ (cf.~Lemma \ref{lem:Phirm}) for \emph{any} function $r$ and integer~$m$ (apart from trivially small $n$).  Indeed, consider any submonoid $\TL_n\leq S\leq\P_n$ for $n\geq3$, and suppose to the contrary that $\Phi = \Phi_{r,m}$ for some $r:S\to\Z$ and $m\in\Z$.  Note then that for any idempotent $e\in E(S)$ we have
\[
\Phi(e,e) = \Phi_{r,m}(e,e) = m - r(e) - r(e) + r(ee) = m - r(e), \qquad\text{so that}\qquad r(e) = m - \Phi(e,e).
\]
Now consider the idempotents
\[
e = \custpartn{1,2,3,4,7}{1,2,3,4,7}{\uarc12\stline33\darc12\stline44\stline77\udotted47\ldotted47}
\COMMA
f = \custpartn{1,2,3,4,7}{1,2,3,4,7}{\uarc23\stline11\darc23\stline44\stline77\udotted47\ldotted47}
\AND
ef = \custpartn{1,2,3,4,7}{1,2,3,4,7}{\uarc12\stline31\darc23\stline44\stline77\udotted47\ldotted47},
\]
all from $\TL_n(\sub S)$.  Since $\Phi(e,e) = \Phi(f,f)=1$ and $\Phi(ef,ef) = 0$, we then have
\[
0 = \Phi(e,f) = \Phi_{r,m}(e,f) = m - r(e) - r(f) + r(ef) = m - (m-1)-(m-1)+m = 2,
\]
a contradiction.
\end{rem}

\begin{rem}
Some products arising from the rank-based twistings are shown in the right-hand columns of Figures~\ref{fig:P2Phi}--\ref{fig:B34Phi} for $\P_2$, $\PB_2$, $\B_3$ and $\B_4$.  These are loose for $\P_2$ and $\PB_2$, but tight for~$\B_3$ and~$\B_4$.  More information about these diagrams will be presented in Remark \ref{rem:eggbox}.
\end{rem}

\section{Green's relations}\label{sect:GR}

The rest of the paper is devoted to describing the structure and properties of a tight twisted product~${M\times_\Phi^qS}$.  We begin here with a description of Green's relations.
Throughout this section, we fix a monoid $S$, a tight twisting $\Phi:S\times S\to\N$, and an (additive) commutative monoid $M$ with a fixed element $q\in M$.  To simplify notation, we will denote the arising twisted product $M\times_\Phi^qS$ by $T$.  

The next result shows how Green's relations on $T = M\times_\Phi^qS$ are built from the corresponding relations on $M$ and~$S$.  To avoid confusion, we write $\xi^M$, $\xi^S$ and $\xi^T$ for the $\xi$-relation on~$M$,~$S$ and~$T$, respectively, where here $\xi$ is any of Green's relations (pre-orders or equivalences).  
Note that  ${\L^M} = {\R^M} = {\J^M} = {\D^M} = {\H^M}$ and ${\leqL^M} = {\leqR^M} = {\leqJ^M} = {\leqH^M}$, because of commutativity.

\begin{thm}\label{thm:GR}
Let $T=M\times_\Phi^q S$ be a tight twisted product.
If $\xi$ is any of Green's relations (pre-orders or equivalences), then 
\[
(i,a) \mr\xi^T (j,b) \ \ \iff\ \  i\mr\xi^M j \ \ \text{and}\ \  a\mr\xi^Sb \qquad\text{for $i,j\in M$ and $a,b\in S$.}
\]
\end{thm}

\pf
We just prove the result when $\xi$ is $\leqJ$.  The proofs for $\leqL$ and $\leqR$ are analogous, and then the remaining ones follow by combining these.

For the forward implication, suppose $(i,a) \leqJ^T(j,b)$, so that
\[
(i,a) = (k,s)(j,b)(l,t) = (j+k+l+\Phi(s,b,t)q,sbt) \qquad\text{for some $k,l\in M$ and $s,t\in S$.}
\]
Then $i = j+k+l+\Phi(s,b,t)q \leqJ^M j$ and $a=sbt\leqJ^S b$.

Conversely, suppose $i\leqJ^M j$ and $a\leqJ^S b$, so that $i = j+l$ for some $l\in M$ (recall that $M$ is commutative) and $a=sbt$ for some $s,t\in S$.  By \ref{T6} we can assume that $\Phi(s,b,t) = 0$.  We then have
\[
(0,s)(j,b)(l,t) = (j+l+\Phi(s,b,t)q,sbt) = (i+0q,a) = (i,a),
\]
and so ${(i,a) \leqJ^T(j,b)}$.
\epf

We state some immediate consequences:

\begin{cor}
\label{co:GR}
Let $T=M\times_\Phi^q S$ be a tight twisted product, and $\K$ any of Green's equivalences.
\ben
\item \label{it:Kia}
For any $i\in M$ and $a\in S$ we have  $K_{(i,a)}^T = K_i^M \times K_a^S = H_i^M \times K_a^S$.
\item\label{it:TKTleq}
For ${\K}\neq {\D}$  we have the poset isomorphism
\[
(T/{\K^T},\leq) \cong (M/{\K^M},\leq)\times(S/{\K^S},\leq) = (M/{\H^M},\leq)\times(S/{\K^S},\leq).
\]
\een
In particular, if $M$ is a group, then $K_{(i,a)}^T = M\times K_a^S$, and $(T/{\K^T},\leq) \cong (S/{\K^S},\leq)$ for ${\K}\not={\D}$. \epfres
\end{cor}

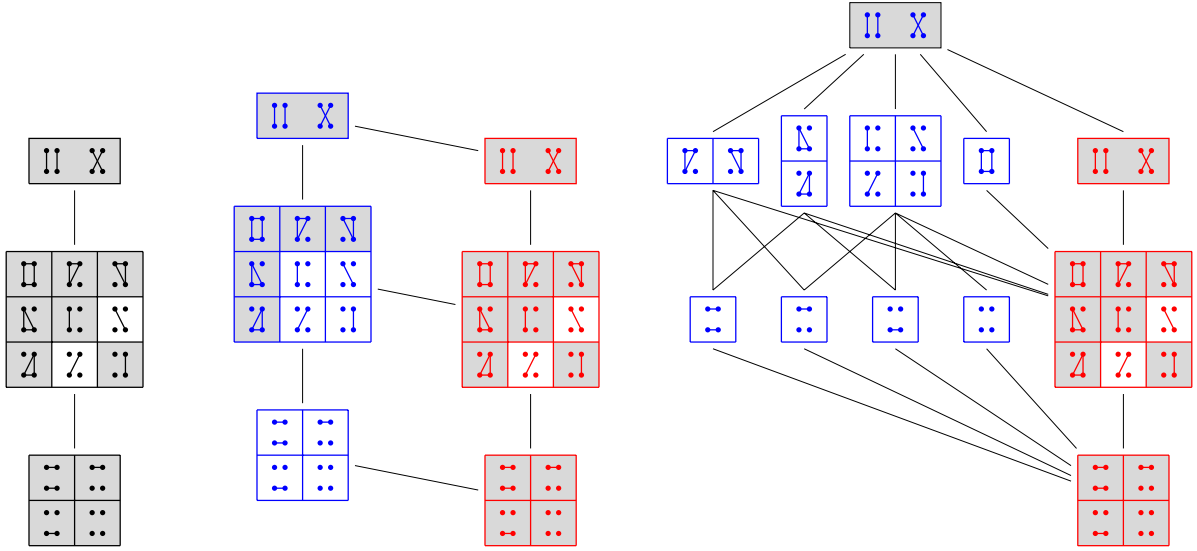
\begin{figure}[t]
\begin{center}
\scalebox{.6}{
\begin{tikzpicture}[scale=1]
%
\begin{scope}[shift={(0,-1)}]
\node (D2) at (0,7.5) {
\begin{tikzpicture}[scale=1]
\nc\xx{1}
\foreach \x/\y in {0/0,1/0} {\fill[thick,lightgray!60] (\x,\y)--(\x,\y+1)--(\x+1,\y+1)--(\x+1,\y)--(\x,\y);}
\draw[thick](0,0)--(2,0)--(2,1)--(0,1)--(0,0);
\node () at (0.5,0.5) {\Ptwo{1/1,2/2}{}{}{}};
\node () at (1.5,0.5) {\Ptwo{1/2,2/1}{}{}{}};
\end{tikzpicture}
};
\node (D1) at (0,4) {
\begin{tikzpicture}[scale=1]
\nc\xx{3}
\foreach \x/\y in {0/2,1/2,2/2,0/1,1/1,0/0,2/0} {\fill[lightgray!60] (\x,\y)--(\x,\y+1)--(\x+1,\y+1)--(\x+1,\y)--(\x,\y);}
\foreach \x in {0,...,\xx} {\draw [thick](\x,0)--(\x,\xx) (0,\x)--(\xx,\x); }
\node () at (0.5,2.5) {\Ptwo{1/1,2/2}{1/2}{1/2}};
\node () at (1.5,2.5) {\Ptwo{1/1,2/1}{1/2}{}};
\node () at (2.5,2.5) {\Ptwo{1/2,2/2}{1/2}{}};
\node () at (0.5,1.5) {\Ptwo{1/1,1/2}{}{1/2}};
\node () at (1.5,1.5) {\Ptwo{1/1}{}{}};
\node () at (2.5,1.5) {\Ptwo{1/2}{}{}};
\node () at (0.5,0.5) {\Ptwo{2/1,2/2}{}{1/2}{}};
\node () at (1.5,0.5) {\Ptwo{2/1}{}{}{}};
\node () at (2.5,0.5) {\Ptwo{2/2}{}{}};
\end{tikzpicture}
};
\node (D0) at (0,0) {
\begin{tikzpicture}[scale=1]
\nc\xx{2}
\foreach \x/\y in {0/0,0/1,1/0,1/1} {\fill[lightgray!60] (\x,\y)--(\x,\y+1)--(\x+1,\y+1)--(\x+1,\y)--(\x,\y);}
\foreach \x in {0,...,\xx} {\draw[thick] (\x,0)--(\x,\xx) (0,\x)--(\xx,\x); }
\node () at (0.5,1.5) {\Ptwo{}{1/2}{1/2}};
\node () at (1.5,1.5) {\Ptwo{}{1/2}{}};
\node () at (0.5,0.5) {\Ptwo{}{}{1/2}{}};
\node () at (1.5,0.5) {\Ptwo{}{}{}};
\end{tikzpicture}
};
\draw
(D2)--(D1)
(D1)--(D0)
;
%
%
\end{scope}
%
%
\begin{scope}[shift={(-1,0)}]
%
%
\begin{scope}[shift={(6,0)}]
\node (D2L) at (0,7.5) {
\begin{tikzpicture}[scale=1]
\nc\xx{1}
\foreach \x/\y in {0/0,1/0} {\fill[thick,blue, lightgray!60] (\x,\y)--(\x,\y+1)--(\x+1,\y+1)--(\x+1,\y)--(\x,\y);}
\draw[thick,blue] (0,0)--(2,0)--(2,1)--(0,1)--(0,0);
\node () at (0.5,0.5) {\BluePtwo{1/1,2/2}{}{}{}};
\node () at (1.5,0.5) {\BluePtwo{1/2,2/1}{}{}{}};
\end{tikzpicture}
};
\node (D1L) at (0,4) {
\begin{tikzpicture}[scale=1]
\nc\xx{3}
\foreach \x/\y in {0/2,1/2,2/2,0/1,0/0} {\fill[thick,blue, lightgray!60] (\x,\y)--(\x,\y+1)--(\x+1,\y+1)--(\x+1,\y)--(\x,\y);}
\foreach \x in {0,...,\xx} {\draw[thick,blue] (\x,0)--(\x,\xx) (0,\x)--(\xx,\x); }
\node () at (0.5,2.5) {\BluePtwo{1/1,2/2}{1/2}{1/2}};
\node () at (1.5,2.5) {\BluePtwo{1/1,2/1}{1/2}{}};
\node () at (2.5,2.5) {\BluePtwo{1/2,2/2}{1/2}{}};
\node () at (0.5,1.5) {\BluePtwo{1/1,1/2}{}{1/2}};
\node () at (1.5,1.5) {\BluePtwo{1/1}{}{}};
\node () at (2.5,1.5) {\BluePtwo{1/2}{}{}};
\node () at (0.5,0.5) {\BluePtwo{2/1,2/2}{}{1/2}{}};
\node () at (1.5,0.5) {\BluePtwo{2/1}{}{}{}};
\node () at (2.5,0.5) {\BluePtwo{2/2}{}{}};
\end{tikzpicture}
};
\node (D0L) at (0,0) {
\begin{tikzpicture}[scale=1]
\nc\xx{2}
\foreach \x/\y in {} {\fill[thick,blue, lightgray!60] (\x,\y)--(\x,\y+1)--(\x+1,\y+1)--(\x+1,\y)--(\x,\y);}
\foreach \x in {0,...,\xx} {\draw[thick,blue] (\x,0)--(\x,\xx) (0,\x)--(\xx,\x); }
\node () at (0.5,1.5) {\BluePtwo{}{1/2}{1/2}};
\node () at (1.5,1.5) {\BluePtwo{}{1/2}{}};
\node () at (0.5,0.5) {\BluePtwo{}{}{1/2}{}};
\node () at (1.5,0.5) {\BluePtwo{}{}{}};
\end{tikzpicture}
};
\node (D2R) at (5,7.5-1) {
\begin{tikzpicture}[scale=1]
\nc\xx{1}
\foreach \x/\y in {0/0,1/0} {\fill[thick,red, lightgray!60] (\x,\y)--(\x,\y+1)--(\x+1,\y+1)--(\x+1,\y)--(\x,\y);}
\draw[thick,red](0,0)--(2,0)--(2,1)--(0,1)--(0,0);
\node () at (0.5,0.5) {\RedPtwo{1/1,2/2}{}{}{}};
\node () at (1.5,0.5) {\RedPtwo{1/2,2/1}{}{}{}};
\end{tikzpicture}
};
\node (D1R) at (5,4-1) {
\begin{tikzpicture}[scale=1]
\nc\xx{3}
\foreach \x/\y in {0/2,1/2,2/2,0/1,1/1,0/0,2/0} {\fill[thick,red, lightgray!60] (\x,\y)--(\x,\y+1)--(\x+1,\y+1)--(\x+1,\y)--(\x,\y);}
\foreach \x in {0,...,\xx} {\draw[thick,red] (\x,0)--(\x,\xx) (0,\x)--(\xx,\x); }
\node () at (0.5,2.5) {\RedPtwo{1/1,2/2}{1/2}{1/2}};
\node () at (1.5,2.5) {\RedPtwo{1/1,2/1}{1/2}{}};
\node () at (2.5,2.5) {\RedPtwo{1/2,2/2}{1/2}{}};
\node () at (0.5,1.5) {\RedPtwo{1/1,1/2}{}{1/2}};
\node () at (1.5,1.5) {\RedPtwo{1/1}{}{}};
\node () at (2.5,1.5) {\RedPtwo{1/2}{}{}};
\node () at (0.5,0.5) {\RedPtwo{2/1,2/2}{}{1/2}{}};
\node () at (1.5,0.5) {\RedPtwo{2/1}{}{}{}};
\node () at (2.5,0.5) {\RedPtwo{2/2}{}{}};
\end{tikzpicture}
};
\node (D0R) at (5,0-1) {
\begin{tikzpicture}[scale=1]
\nc\xx{2}
\foreach \x/\y in {0/0,0/1,1/0,1/1} {\fill[thick,red, lightgray!60] (\x,\y)--(\x,\y+1)--(\x+1,\y+1)--(\x+1,\y)--(\x,\y);}
\foreach \x in {0,...,\xx} {\draw[thick,red] (\x,0)--(\x,\xx) (0,\x)--(\xx,\x); }
\node () at (0.5,1.5) {\RedPtwo{}{1/2}{1/2}};
\node () at (1.5,1.5) {\RedPtwo{}{1/2}{}};
\node () at (0.5,0.5) {\RedPtwo{}{}{1/2}{}};
\node () at (1.5,0.5) {\RedPtwo{}{}{}};
\end{tikzpicture}
};
\draw
(D2L)--(D1L)
(D1L)--(D0L)
(D2R)--(D1R)
(D1R)--(D0R)
(D2L)--(D2R)
(D1L)--(D1R)
(D0L)--(D0R)
;
\end{scope}
%
%
%
\end{scope}
%
%
%
\begin{scope}[shift={(23,-1)}]
%
%
%
%
%
%
%
\nc\xDtwoR0 \nc\yDtwoR{7.5}
\nc\xDoneR0 \nc\yDoneR4
\nc\xDzeroR0 \nc\yDzeroR0
\nc\xDtwoL{-5} \nc\yDtwoL{10.5}
\nc\xa{-9} \nc\ya4
\nc\xb{-7} \nc\yb4
\nc\xc{-5} \nc\yc4
\nc\xd{-3} \nc\yd4
\nc\xA{-9} \nc\yA{7.5}
\nc\xB{-7} \nc\yB{7.5}
\nc\xC{-5} \nc\yC{7.5}
\nc\xD{-3} \nc\yD{7.5}
\node (D2L) at (\xDtwoL,\yDtwoL) {
\begin{tikzpicture}[scale=1]
\nc\xx{1}
\foreach \x/\y in {0/0,1/0} {\fill[thick,blue, lightgray!60] (\x,\y)--(\x,\y+1)--(\x+1,\y+1)--(\x+1,\y)--(\x,\y);}
\draw(0,0)--(2,0)--(2,1)--(0,1)--(0,0);
\node () at (0.5,0.5) {\BluePtwo{1/1,2/2}{}{}{}};
\node () at (1.5,0.5) {\BluePtwo{1/2,2/1}{}{}{}};
\end{tikzpicture}
};
\node (A) at (\xA,\yA) {
\begin{tikzpicture}[scale=1]
\nc\xx{2}
\nc\yy{1}
\foreach \x/\y in {} {\fill[thick,blue, lightgray!60] (\x,\y)--(\x,\y+1)--(\x+1,\y+1)--(\x+1,\y)--(\x,\y);}
\foreach \x in {0,...,\xx} {\draw[thick,blue] (\x,0)--(\x,\yy);}
\foreach \y in {0,...,\yy} {\draw[thick,blue] (0,\y)--(\xx,\y); }
\node () at (0.5,0.5) {\BluePtwo{1/1,2/1}{1/2}{}};
\node () at (1.5,0.5) {\BluePtwo{1/2,2/2}{1/2}{}};
\end{tikzpicture}
};
\node (B) at (\xB,\yB) {
\begin{tikzpicture}[scale=1]
\nc\xx{1}
\nc\yy{2}
\foreach \x/\y in {} {\fill[thick,blue, lightgray!60] (\x,\y)--(\x,\y+1)--(\x+1,\y+1)--(\x+1,\y)--(\x,\y);}
\foreach \x in {0,...,\xx} {\draw[thick,blue] (\x,0)--(\x,\yy);}
\foreach \y in {0,...,\yy} {\draw[thick,blue] (0,\y)--(\xx,\y); }
\node () at (0.5,1.5) {\BluePtwo{1/1,1/2}{}{1/2}};
\node () at (0.5,0.5) {\BluePtwo{2/1,2/2}{}{1/2}{}};
\end{tikzpicture}
};
\node (C) at (\xC,\yC) {
\begin{tikzpicture}[scale=1]
\nc\xx{2}
\nc\yy{2}
\foreach \x/\y in {} {\fill[thick,blue, lightgray!60] (\x,\y)--(\x,\y+1)--(\x+1,\y+1)--(\x+1,\y)--(\x,\y);}
\foreach \x in {0,...,\xx} {\draw[thick,blue] (\x,0)--(\x,\yy);}
\foreach \y in {0,...,\yy} {\draw[thick,blue] (0,\y)--(\xx,\y); }
\node () at (0.5,1.5) {\BluePtwo{1/1}{}{}};
\node () at (1.5,1.5) {\BluePtwo{1/2}{}{}};
\node () at (0.5,0.5) {\BluePtwo{2/1}{}{}{}};
\node () at (1.5,0.5) {\BluePtwo{2/2}{}{}};
\end{tikzpicture}
};
\node (D) at (\xD,\yD) {
\begin{tikzpicture}[scale=1]
\nc\xx{1}
\nc\yy{1}
\foreach \x/\y in {} {\fill[thick,blue, lightgray!60] (\x,\y)--(\x,\y+1)--(\x+1,\y+1)--(\x+1,\y)--(\x,\y);}
\foreach \x in {0,...,\xx} {\draw[thick,blue] (\x,0)--(\x,\yy);}
\foreach \y in {0,...,\yy} {\draw[thick,blue] (0,\y)--(\xx,\y); }
\node () at (0.5,0.5) {\BluePtwo{1/1,2/2}{1/2}{1/2}};
\end{tikzpicture}
};
\node (a) at (\xa,\ya) {
\begin{tikzpicture}[scale=1]
\nc\xx{1}
\foreach \x/\y in {} {\fill[thick,blue, lightgray!60] (\x,\y)--(\x,\y+1)--(\x+1,\y+1)--(\x+1,\y)--(\x,\y);}
\foreach \x in {0,...,\xx} {\draw[thick,blue] (\x,0)--(\x,\xx) (0,\x)--(\xx,\x); }
\node () at (0.5,0.5) {\BluePtwo{}{1/2}{1/2}{}};
\end{tikzpicture}
};
\node (b) at (\xb,\yb) {
\begin{tikzpicture}[scale=1]
\nc\xx{1}
\foreach \x/\y in {} {\fill[thick,blue, lightgray!60] (\x,\y)--(\x,\y+1)--(\x+1,\y+1)--(\x+1,\y)--(\x,\y);}
\foreach \x in {0,...,\xx} {\draw[thick,blue] (\x,0)--(\x,\xx) (0,\x)--(\xx,\x); }
\node () at (0.5,0.5) {\BluePtwo{}{1/2}{}{}};
\end{tikzpicture}
};
\node (c) at (\xc,\yc) {
\begin{tikzpicture}[scale=1]
\nc\xx{1}
\foreach \x/\y in {} {\fill[thick,blue, lightgray!60] (\x,\y)--(\x,\y+1)--(\x+1,\y+1)--(\x+1,\y)--(\x,\y);}
\foreach \x in {0,...,\xx} {\draw[thick,blue] (\x,0)--(\x,\xx) (0,\x)--(\xx,\x); }
\node () at (0.5,0.5) {\BluePtwo{}{}{1/2}{}};
\end{tikzpicture}
};
\node (d) at (\xd,\yd) {
\begin{tikzpicture}[scale=1]
\nc\xx{1}
\foreach \x/\y in {} {\fill[thick,blue, lightgray!60] (\x,\y)--(\x,\y+1)--(\x+1,\y+1)--(\x+1,\y)--(\x,\y);}
\foreach \x in {0,...,\xx} {\draw[thick,blue] (\x,0)--(\x,\xx) (0,\x)--(\xx,\x); }
\node () at (0.5,0.5) {\BluePtwo{}{}{}{}};
\end{tikzpicture}
};
\node (D2R) at (\xDtwoR,\yDtwoR) {
\begin{tikzpicture}[scale=1]
\nc\xx{1}
\foreach \x/\y in {0/0,1/0} {\fill[thick,red, lightgray!60] (\x,\y)--(\x,\y+1)--(\x+1,\y+1)--(\x+1,\y)--(\x,\y);}
\draw[thick,red](0,0)--(2,0)--(2,1)--(0,1)--(0,0);
\node () at (0.5,0.5) {\RedPtwo{1/1,2/2}{}{}{}};
\node () at (1.5,0.5) {\RedPtwo{1/2,2/1}{}{}{}};
\end{tikzpicture}
};
\node (D1R) at (\xDoneR,\yDoneR) {
\begin{tikzpicture}[scale=1]
\nc\xx{3}
\foreach \x/\y in {0/2,1/2,2/2,0/1,1/1,0/0,2/0} {\fill[thick,red, lightgray!60] (\x,\y)--(\x,\y+1)--(\x+1,\y+1)--(\x+1,\y)--(\x,\y);}
\foreach \x in {0,...,\xx} {\draw[thick,red] (\x,0)--(\x,\xx) (0,\x)--(\xx,\x); }
\node () at (0.5,2.5) {\RedPtwo{1/1,2/2}{1/2}{1/2}};
\node () at (1.5,2.5) {\RedPtwo{1/1,2/1}{1/2}{}};
\node () at (2.5,2.5) {\RedPtwo{1/2,2/2}{1/2}{}};
\node () at (0.5,1.5) {\RedPtwo{1/1,1/2}{}{1/2}};
\node () at (1.5,1.5) {\RedPtwo{1/1}{}{}};
\node () at (2.5,1.5) {\RedPtwo{1/2}{}{}};
\node () at (0.5,0.5) {\RedPtwo{2/1,2/2}{}{1/2}{}};
\node () at (1.5,0.5) {\RedPtwo{2/1}{}{}{}};
\node () at (2.5,0.5) {\RedPtwo{2/2}{}{}};
\end{tikzpicture}
};
\node (D0R) at (\xDzeroR,\yDzeroR) {
\begin{tikzpicture}[scale=1]
\nc\xx{2}
\foreach \x/\y in {0/0,0/1,1/0,1/1} {\fill[thick,red, lightgray!60] (\x,\y)--(\x,\y+1)--(\x+1,\y+1)--(\x+1,\y)--(\x,\y);}
\foreach \x in {0,...,\xx} {\draw[thick,red] (\x,0)--(\x,\xx) (0,\x)--(\xx,\x); }
\node () at (0.5,1.5) {\RedPtwo{}{1/2}{1/2}};
\node () at (1.5,1.5) {\RedPtwo{}{1/2}{}};
\node () at (0.5,0.5) {\RedPtwo{}{}{1/2}{}};
\node () at (1.5,0.5) {\RedPtwo{}{}{}};
\end{tikzpicture}
};
\draw
(D2R)--(D1R)
(D1R)--(D0R)
(a.south)--(D0R)
(b.south)--(D0R)
(c.south)--(D0R)
(d.south)--(D0R)
(D2L)--(A.north)
(D2L)--(B.north)
(D2L)--(C.north)
(D2L)--(D.north)
(D2L)--(D2R.north)
(A.south)--(D1R)
(B.south)--(D1R)
(C.south)--(D1R)
(D.south)--(D1R)
(A.south)--(a.north)
(A.south)--(b.north)
(B.south)--(a.north)
(B.south)--(c.north)
(C.south)--(b.north)
(C.south)--(c.north)
(C.south)--(d.north)
;
%
%
%
%
%
\end{scope}
%
%
%
%
%
%
\end{tikzpicture}
}
\caption{Egg-box diagrams of the partition monoid $\P _2$, and two of its twisted products.  See Remark \ref{rem:eggbox} for more details.}
\label{fig:P2Phi}
\end{center}
\end{figure}

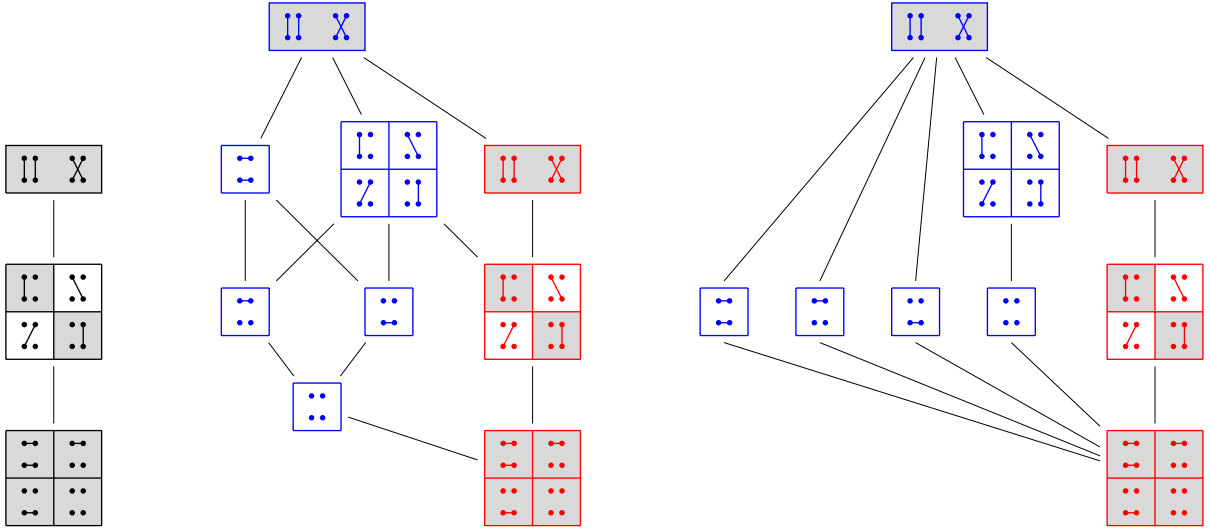
\begin{figure}[t]
\begin{center}
\scalebox{0.63}{
\begin{tikzpicture}[scale=1]
\nc\xDtwoR0 \nc\yDtwoR{6.5}
\nc\xDoneR0 \nc\yDoneR{3.5}
\nc\xDzeroR0 \nc\yDzeroR0
\nc\xDtwoL{-4.5} \nc\yDtwoL{9.5}
\nc\xa{-6} \nc\ya{6.5}
\nc\xb{-6} \nc\yb{3.5}
\nc\xc{-3} \nc\yc{3.5}
\nc\xd{-4.5} \nc\yd{1.5}
\nc\xC{-3} \nc\yC{6.5}
\node (D2L) at (\xDtwoL,\yDtwoL) {
\begin{tikzpicture}[scale=1]
\nc\xx{1}
\foreach \x/\y in {0/0,1/0} {\fill[lightgray!60] (\x,\y)--(\x,\y+1)--(\x+1,\y+1)--(\x+1,\y)--(\x,\y);}
\draw[thick,blue](0,0)--(2,0)--(2,1)--(0,1)--(0,0);
\node () at (0.5,0.5) {\BluePtwo{1/1,2/2}{}{}{}};
\node () at (1.5,0.5) {\BluePtwo{1/2,2/1}{}{}{}};
\end{tikzpicture}
};
\node (C) at (\xC,\yC) {
\begin{tikzpicture}[scale=1]
\nc\xx{2}
\nc\yy{2}
\foreach \x/\y in {} {\fill[lightgray!60] (\x,\y)--(\x,\y+1)--(\x+1,\y+1)--(\x+1,\y)--(\x,\y);}
\foreach \x in {0,...,\xx} {\draw[thick,blue] (\x,0)--(\x,\yy);}
\foreach \y in {0,...,\yy} {\draw [thick,blue](0,\y)--(\xx,\y); }
\node () at (0.5,1.5) {\BluePtwo{1/1}{}{}};
\node () at (1.5,1.5) {\BluePtwo{1/2}{}{}};
\node () at (0.5,0.5) {\BluePtwo{2/1}{}{}{}};
\node () at (1.5,0.5) {\BluePtwo{2/2}{}{}};
\end{tikzpicture}
};
\node (a) at (\xa,\ya) {
\begin{tikzpicture}[scale=1]
\nc\xx{1}
\foreach \x/\y in {} {\fill[lightgray!60] (\x,\y)--(\x,\y+1)--(\x+1,\y+1)--(\x+1,\y)--(\x,\y);}
\foreach \x in {0,...,\xx} {\draw[thick,blue] (\x,0)--(\x,\xx) (0,\x)--(\xx,\x); }
\node () at (0.5,0.5) {\BluePtwo{}{1/2}{1/2}{}};
\end{tikzpicture}
};
\node (b) at (\xb,\yb) {
\begin{tikzpicture}[scale=1]
\nc\xx{1}
\foreach \x/\y in {} {\fill[lightgray!60] (\x,\y)--(\x,\y+1)--(\x+1,\y+1)--(\x+1,\y)--(\x,\y);}
\foreach \x in {0,...,\xx} {\draw[thick,blue] (\x,0)--(\x,\xx) (0,\x)--(\xx,\x); }
\node () at (0.5,0.5) {\BluePtwo{}{1/2}{}{}};
\end{tikzpicture}
};
\node (c) at (\xc,\yc) {
\begin{tikzpicture}[scale=1]
\nc\xx{1}
\foreach \x/\y in {} {\fill[lightgray!60] (\x,\y)--(\x,\y+1)--(\x+1,\y+1)--(\x+1,\y)--(\x,\y);}
\foreach \x in {0,...,\xx} {\draw[thick,blue] (\x,0)--(\x,\xx) (0,\x)--(\xx,\x); }
\node () at (0.5,0.5) {\BluePtwo{}{}{1/2}{}};
\end{tikzpicture}
};
\node (d) at (\xd,\yd) {
\begin{tikzpicture}[scale=1]
\nc\xx{1}
\foreach \x/\y in {} {\fill[lightgray!60] (\x,\y)--(\x,\y+1)--(\x+1,\y+1)--(\x+1,\y)--(\x,\y);}
\foreach \x in {0,...,\xx} {\draw[thick,blue] (\x,0)--(\x,\xx) (0,\x)--(\xx,\x); }
\node () at (0.5,0.5) {\BluePtwo{}{}{}{}};
\end{tikzpicture}
};
\node (D2R) at (\xDtwoR,\yDtwoR) {
\begin{tikzpicture}[scale=1]
\nc\xx{1}
\foreach \x/\y in {0/0,1/0} {\fill[lightgray!60] (\x,\y)--(\x,\y+1)--(\x+1,\y+1)--(\x+1,\y)--(\x,\y);}
\draw[thick,red](0,0)--(2,0)--(2,1)--(0,1)--(0,0);
\node () at (0.5,0.5) {\RedPtwo{1/1,2/2}{}{}{}};
\node () at (1.5,0.5) {\RedPtwo{1/2,2/1}{}{}{}};
\end{tikzpicture}
};
\node (D1R) at (\xDoneR,\yDoneR) {
\begin{tikzpicture}[scale=1]
\nc\xx{2}
\foreach \x/\y in {0/1,1/0} {\fill[lightgray!60] (\x,\y)--(\x,\y+1)--(\x+1,\y+1)--(\x+1,\y)--(\x,\y);}
\foreach \x in {0,...,\xx} {\draw[thick,red] (\x,0)--(\x,\xx) (0,\x)--(\xx,\x); }
\node () at (0.5,1.5) {\RedPtwo{1/1}{}{}};
\node () at (1.5,1.5) {\RedPtwo{1/2}{}{}};
\node () at (0.5,0.5) {\RedPtwo{2/1}{}{}{}};
\node () at (1.5,0.5) {\RedPtwo{2/2}{}{}};
\end{tikzpicture}
};
\node (D0R) at (\xDzeroR,\yDzeroR) {
\begin{tikzpicture}[scale=1]
\nc\xx{2}
\foreach \x/\y in {0/0,0/1,1/0,1/1} {\fill[lightgray!60] (\x,\y)--(\x,\y+1)--(\x+1,\y+1)--(\x+1,\y)--(\x,\y);}
\foreach \x in {0,...,\xx} {\draw[thick,red] (\x,0)--(\x,\xx) (0,\x)--(\xx,\x); }
\node () at (0.5,1.5) {\RedPtwo{}{1/2}{1/2}};
\node () at (1.5,1.5) {\RedPtwo{}{1/2}{}};
\node () at (0.5,0.5) {\RedPtwo{}{}{1/2}{}};
\node () at (1.5,0.5) {\RedPtwo{}{}{}};
\end{tikzpicture}
};
\draw
(D2L)--(D2R)--(D1R)--(D0R)
(a)--(b)--(d)
(a)--(c)--(d)
(d)--(D0R)
(D2L)--(a)
(D2L)--(C)
(C)--(b)
(C)--(c)
(C)--(D1R)
;
%
%
\begin{scope}[shift={(13,0)}]
%
%
\rnc\xDtwoR0 \rnc\yDtwoR{6.5}
\rnc\xDoneR0 \rnc\yDoneR{3.5}
\rnc\xDzeroR0 \rnc\yDzeroR0
\rnc\xDtwoL{-4.5} \rnc\yDtwoL{9.5}
\rnc\xa{-9} \rnc\ya{3.5}
\rnc\xb{-7} \rnc\yb{3.5}
\rnc\xc{-5} \rnc\yc{3.5}
\rnc\xd{-3} \rnc\yd{3.5}
\rnc\xC{-3} \rnc\yC{6.5}
\node (D2L) at (\xDtwoL,\yDtwoL) {
\begin{tikzpicture}[scale=1]
\nc\xx{1}
\foreach \x/\y in {0/0,1/0} {\fill[lightgray!60] (\x,\y)--(\x,\y+1)--(\x+1,\y+1)--(\x+1,\y)--(\x,\y);}
\draw[thick,blue](0,0)--(2,0)--(2,1)--(0,1)--(0,0);
\node () at (0.5,0.5) {\BluePtwo{1/1,2/2}{}{}{}};
\node () at (1.5,0.5) {\BluePtwo{1/2,2/1}{}{}{}};
\end{tikzpicture}
};
\node (C) at (\xC,\yC) {
\begin{tikzpicture}[scale=1]
\nc\xx{2}
\nc\yy{2}
\foreach \x/\y in {} {\fill[lightgray!60] (\x,\y)--(\x,\y+1)--(\x+1,\y+1)--(\x+1,\y)--(\x,\y);}
\foreach \x in {0,...,\xx} {\draw[thick,blue] (\x,0)--(\x,\yy);}
\foreach \y in {0,...,\yy} {\draw[thick,blue] (0,\y)--(\xx,\y); }
\node () at (0.5,1.5) {\BluePtwo{1/1}{}{}};
\node () at (1.5,1.5) {\BluePtwo{1/2}{}{}};
\node () at (0.5,0.5) {\BluePtwo{2/1}{}{}{}};
\node () at (1.5,0.5) {\BluePtwo{2/2}{}{}};
\end{tikzpicture}
};
\node (a) at (\xa,\ya) {
\begin{tikzpicture}[scale=1]
\nc\xx{1}
\foreach \x/\y in {} {\fill[lightgray!60] (\x,\y)--(\x,\y+1)--(\x+1,\y+1)--(\x+1,\y)--(\x,\y);}
\foreach \x in {0,...,\xx} {\draw[thick,blue] (\x,0)--(\x,\xx) (0,\x)--(\xx,\x); }
\node () at (0.5,0.5) {\BluePtwo{}{1/2}{1/2}{}};
\end{tikzpicture}
};
\node (b) at (\xb,\yb) {
\begin{tikzpicture}[scale=1]
\nc\xx{1}
\foreach \x/\y in {} {\fill[lightgray!60] (\x,\y)--(\x,\y+1)--(\x+1,\y+1)--(\x+1,\y)--(\x,\y);}
\foreach \x in {0,...,\xx} {\draw[thick,blue] (\x,0)--(\x,\xx) (0,\x)--(\xx,\x); }
\node () at (0.5,0.5) {\BluePtwo{}{1/2}{}{}};
\end{tikzpicture}
};
\node (c) at (\xc,\yc) {
\begin{tikzpicture}[scale=1]
\nc\xx{1}
\foreach \x/\y in {} {\fill[lightgray!60] (\x,\y)--(\x,\y+1)--(\x+1,\y+1)--(\x+1,\y)--(\x,\y);}
\foreach \x in {0,...,\xx} {\draw[thick,blue] (\x,0)--(\x,\xx) (0,\x)--(\xx,\x); }
\node () at (0.5,0.5) {\BluePtwo{}{}{1/2}{}};
\end{tikzpicture}
};
\node (d) at (\xd,\yd) {
\begin{tikzpicture}[scale=1]
\nc\xx{1}
\foreach \x/\y in {} {\fill[lightgray!60] (\x,\y)--(\x,\y+1)--(\x+1,\y+1)--(\x+1,\y)--(\x,\y);}
\foreach \x in {0,...,\xx} {\draw[thick,blue] (\x,0)--(\x,\xx) (0,\x)--(\xx,\x); }
\node () at (0.5,0.5) {\BluePtwo{}{}{}{}};
\end{tikzpicture}
};
\node (D2R) at (\xDtwoR,\yDtwoR) {
\begin{tikzpicture}[scale=1]
\nc\xx{1}
\foreach \x/\y in {0/0,1/0} {\fill[lightgray!60] (\x,\y)--(\x,\y+1)--(\x+1,\y+1)--(\x+1,\y)--(\x,\y);}
\draw[thick,red](0,0)--(2,0)--(2,1)--(0,1)--(0,0);
\node () at (0.5,0.5) {\RedPtwo{1/1,2/2}{}{}{}};
\node () at (1.5,0.5) {\RedPtwo{1/2,2/1}{}{}{}};
\end{tikzpicture}
};
\node (D1R) at (\xDoneR,\yDoneR) {
\begin{tikzpicture}[scale=1]
\nc\xx{2}
\foreach \x/\y in {0/1,1/0} {\fill[lightgray!60] (\x,\y)--(\x,\y+1)--(\x+1,\y+1)--(\x+1,\y)--(\x,\y);}
\foreach \x in {0,...,\xx} {\draw[thick,red] (\x,0)--(\x,\xx) (0,\x)--(\xx,\x); }
\node () at (0.5,1.5) {\RedPtwo{1/1}{}{}};
\node () at (1.5,1.5) {\RedPtwo{1/2}{}{}};
\node () at (0.5,0.5) {\RedPtwo{2/1}{}{}{}};
\node () at (1.5,0.5) {\RedPtwo{2/2}{}{}};
\end{tikzpicture}
};
\node (D0R) at (\xDzeroR,\yDzeroR) {
\begin{tikzpicture}[scale=1]
\nc\xx{2}
\foreach \x/\y in {0/0,0/1,1/0,1/1} {\fill[lightgray!60] (\x,\y)--(\x,\y+1)--(\x+1,\y+1)--(\x+1,\y)--(\x,\y);}
\foreach \x in {0,...,\xx} {\draw[thick,red] (\x,0)--(\x,\xx) (0,\x)--(\xx,\x); }
\node () at (0.5,1.5) {\RedPtwo{}{1/2}{1/2}};
\node () at (1.5,1.5) {\RedPtwo{}{1/2}{}};
\node () at (0.5,0.5) {\RedPtwo{}{}{1/2}{}};
\node () at (1.5,0.5) {\RedPtwo{}{}{}};
\end{tikzpicture}
};
\draw
(D2L)--(D2R)--(D1R)--(D0R)
(a.south)--(D0R)
(b.south)--(D0R)
(c.south)--(D0R)
(d.south)--(D0R)
(D2L)--(a.north)
(D2L)--(b.north)
(D2L)--(c.north)
(D2L)--(C)--(d)
;
\end{scope}
%
%
%
%
%
\begin{scope}[shift={(-10,0)}]
%
%
\rnc\xDtwoR0 \rnc\yDtwoR{6.5}
\rnc\xDoneR0 \rnc\yDoneR{3.5}
\rnc\xDzeroR0 \rnc\yDzeroR0
\rnc\xDtwoL{-4.5} \rnc\yDtwoL{9.5}
\rnc\xa{-9} \rnc\ya{3.5}
\rnc\xb{-7} \rnc\yb{3.5}
\rnc\xc{-5} \rnc\yc{3.5}
\rnc\xd{-3} \rnc\yd{3.5}
\rnc\xC{-3} \rnc\yC{6.5}
\node (D2R) at (\xDtwoR,\yDtwoR) {
\begin{tikzpicture}[scale=1]
\nc\xx{1}
\foreach \x/\y in {0/0,1/0} {\fill[lightgray!60] (\x,\y)--(\x,\y+1)--(\x+1,\y+1)--(\x+1,\y)--(\x,\y);}
\draw[thick](0,0)--(2,0)--(2,1)--(0,1)--(0,0);
\node () at (0.5,0.5) {\Ptwo{1/1,2/2}{}{}{}};
\node () at (1.5,0.5) {\Ptwo{1/2,2/1}{}{}{}};
\end{tikzpicture}
};
\node (D1R) at (\xDoneR,\yDoneR) {
\begin{tikzpicture}[scale=1]
\nc\xx{2}
\foreach \x/\y in {0/1,1/0} {\fill[lightgray!60] (\x,\y)--(\x,\y+1)--(\x+1,\y+1)--(\x+1,\y)--(\x,\y);}
\foreach \x in {0,...,\xx} {\draw[thick] (\x,0)--(\x,\xx) (0,\x)--(\xx,\x); }
\node () at (0.5,1.5) {\Ptwo{1/1}{}{}};
\node () at (1.5,1.5) {\Ptwo{1/2}{}{}};
\node () at (0.5,0.5) {\Ptwo{2/1}{}{}{}};
\node () at (1.5,0.5) {\Ptwo{2/2}{}{}};
\end{tikzpicture}
};
\node (D0R) at (\xDzeroR,\yDzeroR) {
\begin{tikzpicture}[scale=1]
\nc\xx{2}
\foreach \x/\y in {0/0,0/1,1/0,1/1} {\fill[lightgray!60] (\x,\y)--(\x,\y+1)--(\x+1,\y+1)--(\x+1,\y)--(\x,\y);}
\foreach \x in {0,...,\xx} {\draw[thick] (\x,0)--(\x,\xx) (0,\x)--(\xx,\x); }
\node () at (0.5,1.5) {\Ptwo{}{1/2}{1/2}};
\node () at (1.5,1.5) {\Ptwo{}{1/2}{}};
\node () at (0.5,0.5) {\Ptwo{}{}{1/2}{}};
\node () at (1.5,0.5) {\Ptwo{}{}{}};
\end{tikzpicture}
};
\draw
(D2R)--(D1R)--(D0R)
;
\end{scope}
%
%
%
%
%
%
%
\end{tikzpicture}
}
\caption{Egg-box diagrams of the partial Brauer monoid $\PB _2$, and two of its twisted products.  See Remark \ref{rem:eggbox} for more details.}
\label{fig:PB2Phi}
\end{center}
\end{figure}

\begin{figure}[t]
\begin{center}
\scalebox{.5}{
\begin{tikzpicture}[scale=1]
\begin{scope}[shift={(0,-1.5)}]
\node (D1) at (0,2) {
\begin{tikzpicture}[scale=.5]
\nc\xx{1}
\foreach \x/\y in {0/0} {\fill[lightgray!60] (\x,\y)--(\x,\y+1)--(\x+1,\y+1)--(\x+1,\y)--(\x,\y);}
\foreach \x in {0,...,\xx} {\draw[thick] (\x,0)--(\x,\xx) (0,\x)--(\xx,\x); }
\end{tikzpicture}
};
\node (D0) at (0,0) {
\begin{tikzpicture}[scale=.5]
\nc\xx{3}
\foreach \x/\y in {0/0,0/1,0/2,1/0,1/1,1/2,2/0,2/1,2/2} {\fill[lightgray!60] (\x,\y)--(\x,\y+1)--(\x+1,\y+1)--(\x+1,\y)--(\x,\y);}
\foreach \x in {0,...,\xx} {\draw[thick] (\x,0)--(\x,\xx) (0,\x)--(\xx,\x); }
\end{tikzpicture}
};
\draw
(D1)--(D0)
;
\end{scope}
\begin{scope}[shift={(6,0)}]
\node (D1L) at (0,2) {
\begin{tikzpicture}[scale=.5]
\nc\xx{1}
\foreach \x/\y in {0/0} {\fill[lightgray!60] (\x,\y)--(\x,\y+1)--(\x+1,\y+1)--(\x+1,\y)--(\x,\y);}
\foreach \x in {0,...,\xx} {\draw[thick,blue] (\x,0)--(\x,\xx) (0,\x)--(\xx,\x); }
\end{tikzpicture}
};
\node (D0L) at (0,0) {
\begin{tikzpicture}[scale=.5]
\nc\xx{3}
\foreach \x/\y in {1/2,2/2,0/1,2/1,0/0,1/0} {\fill[lightgray!60] (\x,\y)--(\x,\y+1)--(\x+1,\y+1)--(\x+1,\y)--(\x,\y);}
\foreach \x in {0,...,\xx} {\draw[thick,blue] (\x,0)--(\x,\xx) (0,\x)--(\xx,\x); }
\end{tikzpicture}
};
\end{scope}
\begin{scope}[shift={(10.5,-1.5)}]
\node (D1R) at (0,2) {
\begin{tikzpicture}[scale=.5]
\nc\xx{1}
\foreach \x/\y in {0/0} {\fill[lightgray!60] (\x,\y)--(\x,\y+1)--(\x+1,\y+1)--(\x+1,\y)--(\x,\y);}
\foreach \x in {0,...,\xx} {\draw[thick,red] (\x,0)--(\x,\xx) (0,\x)--(\xx,\x); }
\end{tikzpicture}
};
\node (D0R) at (0,0) {
\begin{tikzpicture}[scale=.5]
\nc\xx{3}
\foreach \x/\y in {0/0,0/1,0/2,1/0,1/1,1/2,2/0,2/1,2/2} {\fill[lightgray!60] (\x,\y)--(\x,\y+1)--(\x+1,\y+1)--(\x+1,\y)--(\x,\y);}
\foreach \x in {0,...,\xx} {\draw[thick,red] (\x,0)--(\x,\xx) (0,\x)--(\xx,\x); }
\end{tikzpicture}
};
\draw
(D1L)--(D0L)
(D1R)--(D0R)
(D1L)--(D1R)
(D0L)--(D0R)
;
\end{scope}
\begin{scope}[shift={(16.5,0)}]
\node (D1L) at (0,2) {
\begin{tikzpicture}[scale=.5]
\nc\xx{1}
\foreach \x/\y in {0/0} {\fill[lightgray!60] (\x,\y)--(\x,\y+1)--(\x+1,\y+1)--(\x+1,\y)--(\x,\y);}
\foreach \x in {0,...,\xx} {\draw[thick,blue] (\x,0)--(\x,\xx) (0,\x)--(\xx,\x); }
\end{tikzpicture}
};
\node (D0L) at (0,0) {
\begin{tikzpicture}[scale=.5]
\nc\xx{3}
\foreach \x in {0,...,\xx} {\draw[thick,blue] (\x,0)--(\x,\xx) (0,\x)--(\xx,\x); }
\end{tikzpicture}
};
\end{scope}
\begin{scope}[shift={(21,-1.5)}]
\node (D1R) at (0,2) {
\begin{tikzpicture}[scale=.5]
\nc\xx{1}
\foreach \x/\y in {0/0} {\fill[lightgray!60] (\x,\y)--(\x,\y+1)--(\x+1,\y+1)--(\x+1,\y)--(\x,\y);}
\foreach \x in {0,...,\xx} {\draw[thick,red] (\x,0)--(\x,\xx) (0,\x)--(\xx,\x); }
\end{tikzpicture}
};
\node (D0R) at (0,0) {
\begin{tikzpicture}[scale=.5]
\nc\xx{3}
\foreach \x/\y in {0/0,0/1,0/2,1/0,1/1,1/2,2/0,2/1,2/2} {\fill[lightgray!60] (\x,\y)--(\x,\y+1)--(\x+1,\y+1)--(\x+1,\y)--(\x,\y);}
\foreach \x in {0,...,\xx} {\draw[thick,red] (\x,0)--(\x,\xx) (0,\x)--(\xx,\x); }
\end{tikzpicture}
};
\draw
(D1L)--(D0L)
(D1R)--(D0R)
(D1L)--(D1R)
(D0L)--(D0R)
;
\end{scope}
%
%
%
%
\begin{scope}[shift={(0,-11)}]
%
\begin{scope}[shift={(0,-1.5)}]
\node (D2) at (0,6.5) {
\begin{tikzpicture}[scale=.5]
\nc\xx{1}
\foreach \x/\y in {0/0} {\fill[lightgray!60] (\x,\y)--(\x,\y+1)--(\x+1,\y+1)--(\x+1,\y)--(\x,\y);}
\foreach \x in {0,...,\xx} {\draw [thick](\x,0)--(\x,\xx) (0,\x)--(\xx,\x); }
\end{tikzpicture}
};
\node (D1) at (0,3.5) {
\begin{tikzpicture}[scale=.5]
\nc\xx{6}
\foreach \x in {0,1,2,3,4} \foreach \y in {5} {\fill[lightgray!60] (\x,\y)--(\x,\y+1)--(\x+1,\y+1)--(\x+1,\y)--(\x,\y);}
\foreach \x in {0,1,2,3,5} \foreach \y in {4} {\fill[lightgray!60] (\x,\y)--(\x,\y+1)--(\x+1,\y+1)--(\x+1,\y)--(\x,\y);}
\foreach \x in {0,1,2,4,5} \foreach \y in {3} {\fill[lightgray!60] (\x,\y)--(\x,\y+1)--(\x+1,\y+1)--(\x+1,\y)--(\x,\y);}
\foreach \x in {0,1,3,4,5} \foreach \y in {2} {\fill[lightgray!60] (\x,\y)--(\x,\y+1)--(\x+1,\y+1)--(\x+1,\y)--(\x,\y);}
\foreach \x in {0,2,3,4,5} \foreach \y in {1} {\fill[lightgray!60] (\x,\y)--(\x,\y+1)--(\x+1,\y+1)--(\x+1,\y)--(\x,\y);}
\foreach \x in {1,2,3,4,5} \foreach \y in {0} {\fill[lightgray!60] (\x,\y)--(\x,\y+1)--(\x+1,\y+1)--(\x+1,\y)--(\x,\y);}
\foreach \x in {0,...,\xx} {\draw[thick] (\x,0)--(\x,\xx) (0,\x)--(\xx,\x); }
\end{tikzpicture}
};
\node (D0) at (0,0) {
\begin{tikzpicture}[scale=.5]
\nc\xx{3}
\foreach \x/\y in {0/0,0/1,0/2,1/0,1/1,1/2,2/0,2/1,2/2} {\fill[lightgray!60] (\x,\y)--(\x,\y+1)--(\x+1,\y+1)--(\x+1,\y)--(\x,\y);}
\foreach \x in {0,...,\xx} {\draw[thick] (\x,0)--(\x,\xx) (0,\x)--(\xx,\x); }
\end{tikzpicture}
};
\draw
(D2)--(D1)
(D1)--(D0)
;
\end{scope}
\begin{scope}[shift={(6,0)}]
\node (D2L) at (0,6.5) {
\begin{tikzpicture}[scale=.5]
\nc\xx{1}
\foreach \x/\y in {0/0} {\fill[lightgray!60] (\x,\y)--(\x,\y+1)--(\x+1,\y+1)--(\x+1,\y)--(\x,\y);}
\foreach \x in {0,...,\xx} {\draw[thick,blue] (\x,0)--(\x,\xx) (0,\x)--(\xx,\x); }
\end{tikzpicture}
};
\node (D1L) at (0,3.5) {
\begin{tikzpicture}[scale=.5]
\nc\xx{6}
\foreach \x in {1,2,3,4} \foreach \y in {5} {\fill[lightgray!60] (\x,\y)--(\x,\y+1)--(\x+1,\y+1)--(\x+1,\y)--(\x,\y);}
\foreach \x in {0,2,3,5} \foreach \y in {4} {\fill[lightgray!60] (\x,\y)--(\x,\y+1)--(\x+1,\y+1)--(\x+1,\y)--(\x,\y);}
\foreach \x in {0,1,4,5} \foreach \y in {3} {\fill[lightgray!60] (\x,\y)--(\x,\y+1)--(\x+1,\y+1)--(\x+1,\y)--(\x,\y);}
\foreach \x in {0,1,4,5} \foreach \y in {2} {\fill[lightgray!60] (\x,\y)--(\x,\y+1)--(\x+1,\y+1)--(\x+1,\y)--(\x,\y);}
\foreach \x in {0,2,3,5} \foreach \y in {1} {\fill[lightgray!60] (\x,\y)--(\x,\y+1)--(\x+1,\y+1)--(\x+1,\y)--(\x,\y);}
\foreach \x in {1,2,3,4} \foreach \y in {0} {\fill[lightgray!60] (\x,\y)--(\x,\y+1)--(\x+1,\y+1)--(\x+1,\y)--(\x,\y);}
\foreach \x in {0,...,\xx} {\draw[thick,blue] (\x,0)--(\x,\xx) (0,\x)--(\xx,\x); }
\end{tikzpicture}
};
\node (D0L) at (0,0) {
\begin{tikzpicture}[scale=.5]
\nc\xx{3}
\foreach \x in {0,...,\xx} {\draw[thick,blue] (\x,0)--(\x,\xx) (0,\x)--(\xx,\x); }
\end{tikzpicture}
};
\end{scope}
\begin{scope}[shift={(10.5,-1.5)}]
\node (D2R) at (0,6.5) {
\begin{tikzpicture}[scale=.5]
\nc\xx{1}
\foreach \x/\y in {0/0} {\fill[lightgray!60] (\x,\y)--(\x,\y+1)--(\x+1,\y+1)--(\x+1,\y)--(\x,\y);}
\foreach \x in {0,...,\xx} {\draw[thick,red] (\x,0)--(\x,\xx) (0,\x)--(\xx,\x); }
\end{tikzpicture}
};
\node (D1R) at (0,3.5) {
\begin{tikzpicture}[scale=.5]
\nc\xx{6}
\foreach \x in {0,1,2,3,4} \foreach \y in {5} {\fill[lightgray!60] (\x,\y)--(\x,\y+1)--(\x+1,\y+1)--(\x+1,\y)--(\x,\y);}
\foreach \x in {0,1,2,3,5} \foreach \y in {4} {\fill[lightgray!60] (\x,\y)--(\x,\y+1)--(\x+1,\y+1)--(\x+1,\y)--(\x,\y);}
\foreach \x in {0,1,2,4,5} \foreach \y in {3} {\fill[lightgray!60] (\x,\y)--(\x,\y+1)--(\x+1,\y+1)--(\x+1,\y)--(\x,\y);}
\foreach \x in {0,1,3,4,5} \foreach \y in {2} {\fill[lightgray!60] (\x,\y)--(\x,\y+1)--(\x+1,\y+1)--(\x+1,\y)--(\x,\y);}
\foreach \x in {0,2,3,4,5} \foreach \y in {1} {\fill[lightgray!60] (\x,\y)--(\x,\y+1)--(\x+1,\y+1)--(\x+1,\y)--(\x,\y);}
\foreach \x in {1,2,3,4,5} \foreach \y in {0} {\fill[lightgray!60] (\x,\y)--(\x,\y+1)--(\x+1,\y+1)--(\x+1,\y)--(\x,\y);}
\foreach \x in {0,...,\xx} {\draw[thick,red] (\x,0)--(\x,\xx) (0,\x)--(\xx,\x); }
\end{tikzpicture}
};
\node (D0R) at (0,0) {
\begin{tikzpicture}[scale=.5]
\nc\xx{3}
\foreach \x/\y in {0/0,0/1,0/2,1/0,1/1,1/2,2/0,2/1,2/2} {\fill[lightgray!60] (\x,\y)--(\x,\y+1)--(\x+1,\y+1)--(\x+1,\y)--(\x,\y);}
\foreach \x in {0,...,\xx} {\draw[thick,red] (\x,0)--(\x,\xx) (0,\x)--(\xx,\x); }
\end{tikzpicture}
};
\draw
(D2L)--(D1L)--(D0L)
(D2R)--(D1R)--(D0R)
(D2L)--(D2R)
(D1L)--(D1R)
(D0L)--(D0R)
;
\end{scope}
\begin{scope}[shift={(7+9.5,0)}]
\node (D2L) at (0,6.5) {
\begin{tikzpicture}[scale=.5]
\nc\xx{1}
\foreach \x/\y in {0/0} {\fill[lightgray!60] (\x,\y)--(\x,\y+1)--(\x+1,\y+1)--(\x+1,\y)--(\x,\y);}
\foreach \x in {0,...,\xx} {\draw[thick,blue] (\x,0)--(\x,\xx) (0,\x)--(\xx,\x); }
\end{tikzpicture}
};
\node (D1L) at (0,3.5) {
\begin{tikzpicture}[scale=.5]
\nc\xx{6}
\foreach \x in {0,...,\xx} {\draw[thick,blue] (\x,0)--(\x,\xx) (0,\x)--(\xx,\x); }
\end{tikzpicture}
};
\node (D0L) at (0,0) {
\begin{tikzpicture}[scale=.5]
\nc\xx{3}
\foreach \x in {0,...,\xx} {\draw[thick,blue] (\x,0)--(\x,\xx) (0,\x)--(\xx,\x); }
\end{tikzpicture}
};
\end{scope}
\begin{scope}[shift={(11.5+9.5,-1.5)}]
\node (D2R) at (0,6.5) {
\begin{tikzpicture}[scale=.5]
\nc\xx{1}
\foreach \x/\y in {0/0} {\fill[lightgray!60] (\x,\y)--(\x,\y+1)--(\x+1,\y+1)--(\x+1,\y)--(\x,\y);}
\foreach \x in {0,...,\xx} {\draw[thick,red] (\x,0)--(\x,\xx) (0,\x)--(\xx,\x); }
\end{tikzpicture}
};
\node (D1R) at (0,3.5) {
\begin{tikzpicture}[scale=.5]
\nc\xx{6}
\foreach \x in {0,1,2,3,4} \foreach \y in {5} {\fill[lightgray!60] (\x,\y)--(\x,\y+1)--(\x+1,\y+1)--(\x+1,\y)--(\x,\y);}
\foreach \x in {0,1,2,3,5} \foreach \y in {4} {\fill[lightgray!60] (\x,\y)--(\x,\y+1)--(\x+1,\y+1)--(\x+1,\y)--(\x,\y);}
\foreach \x in {0,1,2,4,5} \foreach \y in {3} {\fill[lightgray!60] (\x,\y)--(\x,\y+1)--(\x+1,\y+1)--(\x+1,\y)--(\x,\y);}
\foreach \x in {0,1,3,4,5} \foreach \y in {2} {\fill[lightgray!60] (\x,\y)--(\x,\y+1)--(\x+1,\y+1)--(\x+1,\y)--(\x,\y);}
\foreach \x in {0,2,3,4,5} \foreach \y in {1} {\fill[lightgray!60] (\x,\y)--(\x,\y+1)--(\x+1,\y+1)--(\x+1,\y)--(\x,\y);}
\foreach \x in {1,2,3,4,5} \foreach \y in {0} {\fill[lightgray!60] (\x,\y)--(\x,\y+1)--(\x+1,\y+1)--(\x+1,\y)--(\x,\y);}
\foreach \x in {0,...,\xx} {\draw[thick,red] (\x,0)--(\x,\xx) (0,\x)--(\xx,\x); }
\end{tikzpicture}
};
\node (D0R) at (0,0) {
\begin{tikzpicture}[scale=.5]
\nc\xx{3}
\foreach \x/\y in {0/0,0/1,0/2,1/0,1/1,1/2,2/0,2/1,2/2} {\fill[lightgray!60] (\x,\y)--(\x,\y+1)--(\x+1,\y+1)--(\x+1,\y)--(\x,\y);}
\foreach \x in {0,...,\xx} {\draw[thick,red] (\x,0)--(\x,\xx) (0,\x)--(\xx,\x); }
\end{tikzpicture}
};
\draw
(D2L)--(D1L)--(D0L)
(D2R)--(D1R)--(D0R)
(D2L)--(D2R)
(D1L)--(D1R)
(D0L)--(D0R)
;
\end{scope}
%
%
%
%
\end{scope}
\end{tikzpicture}
}
\caption{Egg-box diagrams of the Brauer monoids~$\B_3$ (top row) and $\B_4$ (bottom row), and two each of their twisted products.  See Remark \ref{rem:eggbox} for more details.}
\label{fig:B34Phi}
\end{center}
\end{figure}
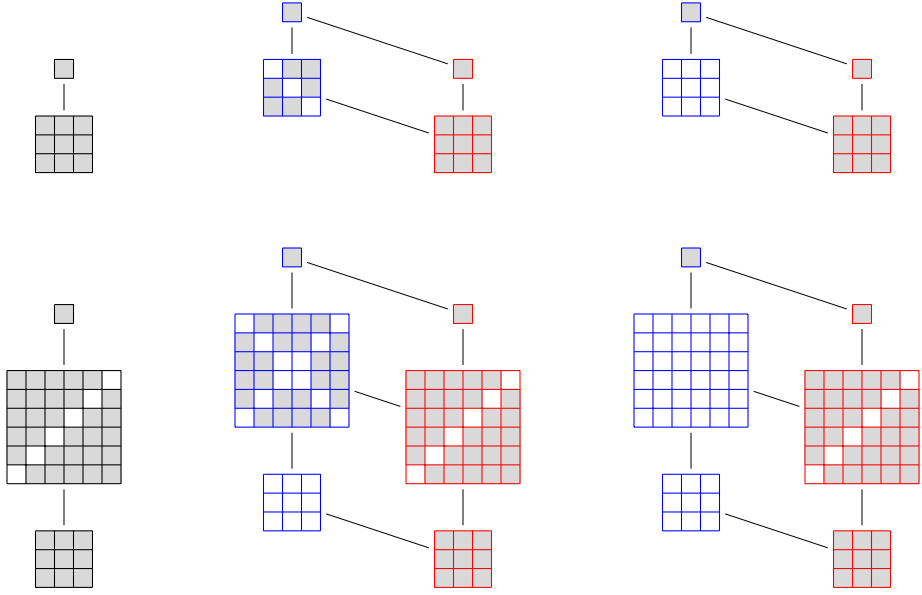

\begin{rem}\label{rem:eggbox}
We now discuss some examples of tight twisted products, to illustrate Theorem~\ref{thm:GR} and Corollary \ref{co:GR}. For comparison we also  consider some loose products, where these results no longer apply.  All examples will involve the two-element additive monoid $M = \{0,\infty\}$, and we will take $q=\infty$. In each case $S$ is a diagram monoid of small degree. The tightness or looseness of the considered twistings always follows from  Proposition~\ref{prop:canPhi} or \ref{prop:altPhiBn}.
The Green structure for the loose examples was computed using the Semigroups package for GAP \cite{GAP,Semigroups}.

We start with $S = \P_2$, the partition monoid of degree $2$, whose egg-box diagram is shown in Figure \ref{fig:P2Phi} (left).  The tight twisted product $T = M\times_\Phi^\infty S$, arising from the canonical twisting~$\Phi$, is shown in Figure~\ref{fig:P2Phi} (middle).  Here blue and red diagrams of $a\in\P_2$ represent $(0,a)$ and~$(\infty,a)$, respectively, and group $\H$-classes are shaded grey; the same colouring conventions apply to the other diagrams discussed below.
The posets $(S/{\J^S},\leq) \cong {\bf3}$ and ${(M/{\J^M},\leq) \cong {\bf2}}$ are chains of sizes~$3$ and~$2$, respectively, while $(T/{\J^T},\leq) \cong {\bf3}\times{\bf2}$ is the direct product, as can be readily seen in the diagram.

By contrast, Figure \ref{fig:P2Phi} (right) pictures the egg-box diagram of the twisted product $M\times_\Phi^\infty \P_2$, where this time $\Phi$ is the rank-based twisting.  Since this twisting is loose, Theorem \ref{thm:GR} does not apply here, and one can immediately see from the diagram that the $\J$-class poset is not isomorphic to the direct product ${\bf3}\times{\bf2}$.  Indeed, while the sets $\{0\}\times D_1$ and $\{0\}\times D_0$ are whole $\D$-classes in the product arising from the canonical twisting, they `fall apart' into multiple $\D$-classes for the rank-based twisting.

Figure \ref{fig:PB2Phi} follows the same pattern for the partial Brauer monoid $\PB_2$, showing egg-box diagrams for $\PB_2$ itself (left), and the products $M\times_\Phi^\infty\PB_2$ arising from the canonical twisting (middle) and the rank-based twisting (right).  Both of these products are loose.  Although they have the same ${\D}={\J}$-classes as each other, the orderings on these are different.

Figure \ref{fig:B34Phi} treats  the Brauer monoids $\B_3$ (top row) and $\B_4$ (bottom row).  For reasons of space, we have not pictured individual elements here; $\B_4$ has size $105$.  The canonical and rank-based twistings are both tight, and in each case the poset of $\J$-classes decomposes as a direct product.
\end{rem}

As an application of the above results, we can also deduce a useful result concerning the property of stability.
A semigroup $S$ is said to be \emph{stable} if 
\begin{equation}\label{eq:stable}
a \J ab \iff a\R ab \AND a\J ba \iff a\L ba \qquad\text{for all $a,b\in S$.}
\end{equation}
Stable semigroups have many important and useful properties, as explained in \cite{EH2020}.  It is known that all finite semigroups are stable, and that ${\D}={\J}$ in any stable semigroup; see for example \cite[Section A.2]{RS2009}.  Any commutative semigroup is stable because for such a semigroup we have ${\J}={\R}={\L}$.

\begin{prop}\label{prop:stability}
A tight twisted product $T = M\times_\Phi^qS$ is stable if and only if $S$ is stable.
\end{prop}

\pf
Throughout the proof we use Theorem \ref{thm:GR} extensively, without explicit mention.  
Suppose first that $S$ is stable.  To demonstrate stability of $T$, it suffices by symmetry to establish the first equivalence in \eqref{eq:stable}, and since ${\R}\sub{\J}$, only the forward implication is required.  So suppose
\[
(i,a) \J^T (i,a)(j,b) = (i+j+\Phi(a,b)q,ab) \qquad\text{for some $i,j\in M$ and $a,b\in S$.}
\]
We then have $i\J^Mi+j+\Phi(a,b)q$ and $a\J^Sab$.  But $M$ and $S$ are both stable (by assumption for $S$, and by commutativity for $M$), so it follows that $i\R^Mi+j+\Phi(a,b)q$ and $a\R^Sab$.  But this then gives
\[
(i,a) \R^T (i+j+\Phi(a,b)q,ab) = (i,a)(j,b).
\]

Conversely, suppose $T$ is stable.  Again aiming to prove the forward implication in the first equivalence in \eqref{eq:stable}, suppose $a\J^Sab$ for some $a,b\in S$.  By \ref{T2}\ref{T2b}, we can assume that~${\Phi(a,b)=0}$.  We then have
\[
(0,a) \J^T (0,ab) = (0,a)(0,b),
\]
so it follows from stability of $T$ that $(0,a) \R^T (0,a)(0,b) = (0,ab)$, and hence that $a\R^Sab$.
\epf

In \cite[Corollary 3]{KV2023} it is proved that the twisted Brauer monoid $\Z\times_\Phi^1\B_n$ is stable when $\Phi$ is the canonical twisting.  This also follows from Proposition \ref{prop:stability} because~$\B_n$ is finite and hence stable.

\section{Idempotents, Sch\"utzenberger groups and biordered sets}\label{sect:ET}

We continue to fix a tight twisted product $T = M\times_\Phi^qS$, where $M$, $S$, $\Phi$ and $q$ have their usual meanings.
In this section we describe the set $E(T)$ of idempotents of $T$ (Section \ref{subsect:E}), and its structure as a biordered set (Section \ref{subsect:ET}).  
In between, we show how the Sch\"utzenberger groups of $T$ are related to those of $M$ and $S$ (Section \ref{subsect:Ga}).

\subsection{The idempotents}\label{subsect:E}

Before we look at the idempotents of the tight twisted product $T=M\times_\Phi^qS$ directly, it is convenient to describe the group $\H^T$-classes.  We begin with a technical lemma concerning the additive monoid $(M,+)$, whose simple proof is omitted; for similar results see \cite[Section~I.4]{Grillet2001}.

\begin{lemma}\label{lem:Hk}
If $H$ is a group $\H^M$-class of $M$, and if $k\in M$, then the following are equivalent:
\ben
\item \label{Hki} $H = H+k$,
\item \label{Hkii} $i\leqJ^M k$ for some $i\in H$,
\item \label{Hkiii} $i\leqJ^M k$ for all $i\in H$.  \epfres
\een
\end{lemma}

If $H$ is an $\H^S$-class (of $S$), then it quickly follows from \ref{T3} and \ref{T4} that $\Phi(a,b) = \Phi(c,d)$ for all $a,b,c,d\in H$.  We write $\Phi(H)$ for this common value (so $\Phi(H) = \Phi(a,b)$ for any $a,b\in H$).

In the following proof, we will make use of the fact that an $\H$-class $H$ of a semigroup is a group if and only if $xy\in H$ for some $x,y\in H$; see for example \cite[Theorem 2.2.5]{Howie1995}.

\begin{thm}\label{thm:H}
Let $T = M\times_\Phi^qS$ be a tight twisted product.  If $H$ is an $\H^M$-class, and $H'$ an $\H^S$-class, then the following are equivalent:
\ben
\item \label{Hi} the $\H^T$-class $H\times H'$ is a group,
\item \label{Hii} $H$ and $H'$ are groups, and $H=H+\Phi(H')q$,
\item \label{Hiii} $H$ and $H'$ are groups, and $i\leqJ^M\Phi(H')q$ for some (equivalently all) $i\in H$,
\item \label{Hiv} $H$ and $H'$ are groups, and either $\Phi(H')=0$ or else $i\leqJ^Mq$ for some (equivalently all)~${i\in H}$.
\een
When the above conditions hold, the group $\H^T$-class $H\times H'$ is isomorphic to the direct product $(H,+)\times(H',\cdot)$.
\end{thm}

\pf
\firstpfitem{\ref{Hi}$\implies$\ref{Hiii}}  Suppose $H\times H'$ is a group, say with identity $(i,a)$.  Since
\[
(i,a) = (i,a)^2 = (2i+\Phi(a,a)q,a^2),
\]
it follows that $a=a^2$ and $i = 2i+\Phi(a,a)q = 2i+\Phi(H')q$.  The former tells us that $H_a^S = H'$ is a group.  The latter tells us that $i \leqJ^M\Phi(H')q$, and also that $i\R^M2i$.  Since ${\R^M}={\H^M}$, this gives $i\H^M2i$, which implies that $H_i^M=H$ is a group (as mentioned before the proof).

\pfitem{\ref{Hiii}$\implies$\ref{Hii}}  This follows from Lemma \ref{lem:Hk}.

\pfitem{\ref{Hii}$\implies$\ref{Hi}}  
Suppose $H$ and $H'$ are groups, with $H=H+\Phi(H')q$.  Let the identities of these groups be $i\in H$ and $e\in H'$, and let $p\in H$ be such that $i = p + \Phi(H')q = p + \Phi(e,e)q$.  Since $i$ is the identity of $H$, we have
\[
(p,e)(p,e) = (2p+\Phi(e,e)q,e^2) = (p+i,e) = (p,e),
\]
and hence $(p,e)\in H\times H'$ is an idempotent.

\pfitem{\ref{Hiii}$\implies $\ref{Hiv}}  This is clear.

\pfitem{\ref{Hiv}$\implies $\ref{Hiii}}
We only need to check that $i\leqJ^M \Phi(H')q$. If $\Phi(H')=0$ this is obvious.
When $\Phi(H')\neq 0$ we have $\Phi(H')i\H^M i$ because $H$ is a group, and hence
$i\leqJ^M \Phi(H')i\leqJ^M\Phi(H')q$, as requied.

\aftercases For the final assertion, suppose $H\times H'$ is a group.  Let~$i$ and $e$ be the identities of the groups~$H$ and $H'$, respectively, and write $m = \Phi(H')$, so that
\[
i \leqJ^M mq
\COMMA
H = H+mq
\AND
\Phi(a,b) = m \text{ for all $a,b\in H'$.}
\]
Also let $t = i + mq$; since $H=H+mq$ we have $t\in H$.  For any $j,k\in H$ and $a,b\in H'$ we have
\[
(j,a)(k,b) = (j+k+\Phi(a,b)q,ab) = (j+k+i+mq,ab) = (j+t+k,ab).
\]
It follows that the $\H^T$-class $H\times H'$ is isomorphic to the direct product $(H,\oplus)\times(H',\cdot)$, where~$\oplus$ is the \emph{variant} operation \cite{Hickey1983} given by $j\oplus k = j+t+k$.  Since the map $j\mt j+t$ is easily seen to be an isomorphism $(H,\oplus) \to (H,+)$, the proof is complete.
\epf

We now define the set 
\begin{equation}\label{eq:Lam}
\begin{split}
 \Om &= \bigset{(i,e)\in E(M)\times E(S)}{H_i^M\times H_e^S \text{ is a group}} \\
 &= \bigset{(i,e)\in E(M)\times E(S)}{i \leqJ^M\Phi(e,e)q},
\end{split}
\end{equation}
where the second equality follows from Theorem \ref{thm:H}.  For $(i,e)\in\Om$, we write $\ve(i,e)$ for the unique idempotent in $H_i^M\times H_e^S$, which is of course the identity element of this group $\H^T$-class.  The next result follows from Theorem \ref{thm:H} and its proof.

\begin{prop}\label{prop:ET}
We have $E(T) = E(M\times_\Phi^qS) = \set{\ve(i,e)}{(i,e)\in\Om}$, where $\Om$ is as in \eqref{eq:Lam}.  For any $(i,e)\in\Om$, we have
\[
\ve(i,e) = (p,e), \WHERE p = p(i,e) \in H_i^M\text{ is such that $i = p + \Phi(e,e)q$.}  \epfreseq
\]
\end{prop}

\begin{rem}\label{rem:ETunit}
If $q$ is a unit, then so too is $\Phi(e,e)q$ for any $e\in E(S)$.  In this case it follows that $\Om = E(M)\times E(S)$, and that
\[
\ve(i,e) = (i-\Phi(e,e)q,e) \qquad\text{for all $i\in E(M)$ and $e\in E(S)$.}
\]
In particular, if $M$ is a group, so that $E(M)=\{0\}$, the idempotents of $T$ all have the form $\ve(0,e) = (-\Phi(e,e)q,e)$ for $e\in E(S)$.

On the other hand, if $e\in E(S)$ is such that $\Phi(e,e)=0$, then $(i,e)\in\Om$ for all $i\in E(M)$.  In this case we have $p(i,e)=i$, and so $\ve(i,e)=(i,e)$.

If $M$ is unipotent, i.e.~$E(M)=\{0\}$, but $q$ is not a unit, then $0\not\leqJ^Mq$, and so
\[
\Om = \set{(0,e)}{e\in E(S),\ \Phi(e,e) = 0} = E(T).
\]
\end{rem}

\begin{rem}
Figure \ref{fig:P2Phi} pictures the partition monoid $\P_2$ (left) as well as the tight twisted product $T = M\times_\Phi^\infty\P_2$ (middle), where $M=\{0,\infty\}$, and where $\Phi$ is the canonical twisting.  Since~$\infty$ is an absorbing element of $M$, we see that $(\infty,e)$ is an idempotent of $T$ for every idempotent~$e$ of~$\P_2$.  On the other hand, $(0,e)$ is an idempotent of $T$ for 
$e = \custpartn{1,2}{1,2}{\stline11\stline22} $,~%
$\custpartn{1,2}{1,2}{\stline11\stline22\uline12\dline12} $,~%
$\custpartn{1,2}{1,2}{\stline11\stline21\uline12} $,~%
$\custpartn{1,2}{1,2}{\stline12\stline22\uline12} $,~%
$\custpartn{1,2}{1,2}{\stline11\stline12\dline12} $ and~%
$\custpartn{1,2}{1,2}{\stline21\stline22\dline12}$,
but not for
$e = \custpartn{1,2}{1,2}{\stline11} $, 
$\custpartn{1,2}{1,2}{\stline22} $, 
$\custpartn{1,2}{1,2}{\uline12\dline12} $, 
$\custpartn{1,2}{1,2}{\uline12} $, 
$\custpartn{1,2}{1,2}{\dline12} $ or
$\custpartn{1,2}{1,2}{} $.
This is because $0 \not\leqJ^M \infty = q$, and $\Phi(e,e)=0$ for all the idempotents in the first list, but $\Phi(e,e)>0$ for those in the second.

Similar comments apply to both (tight) twisted products $T=M\times_\Phi^\infty\B_4$ pictured in Figure~\ref{fig:B34Phi} (and also to $M\times_\Phi^\infty\B_3$).  For the canonical twisting we have, among many examples, ${(0,e)\in E(T)}$ and ${(0,f)\not\in E(T)}$ for the idempotents $e = \custpartn{1,2,3,4}{1,2,3,4}{\uarc12\darc23\stline31\stline44}$ and $f = \custpartn{1,2,3,4}{1,2,3,4}{\uarc12\darc12\stline33\stline44}$ from $\B_4$, as $\Phi(e,e)=0$ and $\Phi(f,f)=1$.  On the other hand, for the rank-based twisting we have $\Phi(e,e) = 4-\rank(e) >0$ for any idempotent $e\in E(\B_4)$ of rank less than $4$.
\end{rem}

\subsection{Sch\"utzenberger groups}\label{subsect:Ga}

An arbitrary $\H$-class $H$ of a semigroup can be given a group structure, even if $H$ does not contain an idempotent; 
these are known as Sch\"{u}tzenberger groups. 
In this section we want to describe the Sch\"{u}tzenberger groups of a tight twisted product.

Let us first briefly review the basic definitions and facts; for details and proofs see \cite[Section~2.4]{CPbook} and \cite[pp.~166--167]{Arbib1968}.
Let $S$ be an arbitrary semigroup.
Define the set
\[
P = P(H) = \set{u\in S^1}{Hu\cap H\not=\es} = \set{u\in S^1}{Hu=H}.
\]
Each $u\in P$ determines a bijection $\rho_u:H\to H:a\mt au$, and the set of all these forms the Sch\"utzenberger group of $H$:
\[
\Ga = \Ga(H) = \set{\rho_u}{u\in P(H)}.
\]
This is a simply transitive group of permutations of $H$, and in the case that $H$ is a subgroup of~$S$ we have $\Ga(H)\cong H$.
We can also give $H$ itself a group operation $\star$, as follows.  For this, fix some element $h\in H$.  For $a,b\in H$ we define
\begin{equation}\label{eq:star}
a \star b = a\be \qquad\text{where $\be\in\Ga$ is such that $b = h\be$.}
\end{equation}
Note that $\be$ here is uniquely determined, due to $\Ga$ being simply transitive.  The group $(H,\star)$ is isomorphic to $\Ga(H)$, so we refer to both groups as `the Sch\"utzenberger group' of $H$.
All $\H$-classes in a common $\D$-class have isomorphic Sch\"utzenberger groups.

Returning now to tight twisted products,
we saw in Theorem \ref{thm:H} that the group $\H$-classes in $T = M\times_\Phi^qS$ are direct products of group $\H$-classes of $M$ and $S$.  This holds more generally for Sch\"utzenbeger groups:

\begin{thm}\label{thm:Ga}
For any $\H^T$-class $H\times H'$ of a tight twisted product $T = M\times_\Phi^qS$, we have
\[
\Ga(H\times H') \cong \Ga(H) \times \Ga(H').
\]
\end{thm}

\pf
Fix elements $i\in H$ and $a\in H'$, so that
\[
H = H_i^M \COMMA H' = H_a^S \AND H\times H' = H_{(i,a)}^T.
\]
Let
\begin{align*}
P = P(H) &= \set{j\in M}{H+j= H} &\text{and} && P' = P(H') &= \set{u\in S}{H'u=H'}\\
&= \set{j\in M}{i+j\in H} & && &= \set{u\in S}{au\in H'},
\end{align*}
and also define
\[
P'' = \set{u\in P'}{\Phi(a,u) = 0} = \set{u\in P'}{\Phi(b,u) = 0 \text{ for all $b\in H'$}},
\]
noting that the second equality follows from \ref{T3}.  We first claim that
\begin{equation}\label{eq:GaH'}
\Ga(H') = \set{\rho_u}{u\in P''}.
\end{equation}
To prove this, we need to show that for any $u\in P'$ we have $\rho_u = \rho_v$ for some $v\in P''$.  By~\ref{T2} we have $au = av$ for some $v\in S$ with $\Phi(a,v) = 0$.
Since also $av = au \in H'$, we have $v\in P'$, so indeed $v\in P''$.  Now for any $b\in H'$ we have $b = sa$ for some $s\in S$, and so $b\rho_v = sav = sau = b\rho_u$, which gives $\rho_v=\rho_u$.  This completes the proof of \eqref{eq:GaH'}.

Returning now to the main proof, we will show that
\[
(H\times H',\star) \cong(H,\star) \times (H',\star),
\]
where here we write $\star$ for the operations in all three Sch\"utzenberger groups, as in~\eqref{eq:star}, with respect to the fixed elements $i\in H$, $a\in H'$, and $(i,a)\in H\times H'$.  Now consider two elements $(j,b)$ and $(k,c)$ of $H\times H'$.  We must show that
\[
(j,b) \star (k,c) = (j\star k,b\star c).
\]
By definition, we have
\[
j\star k = j\ka \ANd b\star c = b\ga, \quad\text{where $\ka\in\Ga(H)$ and $\ga\in\Ga(H')$ are such that $k=i\ka$ and $c=a\ga$.}
\]
Now, $\ka = \rho_h$ for some $h\in P$, and by \eqref{eq:GaH'} we have $\ga = \rho_u$ for some $u\in P''$.  Note that
\[
(i,a)(h,u) = (i+h+\Phi(a,u)q,au) = (i\rho_h,a\rho_u) = (i\ka,a\ga) = (k,c).
\]
Since $(i,a)$ and $(k,c)$ both belong to $H\times H'$, this shows that $(h,u) \in P(H\times H')$, and also that $(i,a)\rho_{(h,u)} = (k,c)$.  It follows that
\begin{align*}
(j,b) \star (k,c) = (j,b)\rho_{(h,u)} = (j,b)(h,u) = (j+h+\Phi(b,u)q,bu) = (j\rho_h,b\rho_u) &= (j\ka,b\ga) \\
&= (j\star k,b\star c),
\end{align*}
as required.
\epf

\subsection{The biordered set}\label{subsect:ET}

The set of idempotents of a semigroup need not be a subsemigroup, but it has a partial operation giving it the structure of a so-called biordered set, as defined by Nambooripad~\cite{Nambooripad1979}.  
Having characterised the idempotents of a tight twisted product, here we wish to describe the associated biordered structure.  

First we give the relevant definitions.
Let $S$ be an arbitrary semigroup.
Define two pre-orders, $\larr$ and $\rarr$ on $ E(S)$ by
\[
e\larr f \iff e=ef  \AND e\rarr f \iff e=fe.
\]
A pair $(e,f)$ of idempotents is called \emph{basic} if (at least) one of $e\larr f$, $e\rarr f$, $f\larr e$ or $f\rarr e$ holds.  In this case, $ef$ and $fe$ are both idempotents (at least one of which is equal to $e$ or $f$).  
The \emph{biordered set} of $S$ is the set $E(S)$ with the partial product which is the restriction of the product of $S$ to the set of all basic pairs.

For two biordered sets $E(S)$ and $E(T)$, their direct product $E(S)\times E(T)$ is the partial algebra in which $(e_1,f_1)(e_2,f_2)$ is defined if and only if $e_1e_2$ is defined in $E(S)$ and $f_1f_2$ is defined in $E(T)$, in which case $(e_1,f_1)(e_2,f_2)=(e_1e_2,f_1f_2)$.

A \emph{morphism} of biordered sets is a map $\phi:E(S)\to E(T)$, for semigroups $S$ and $T$, such that whenever $ef$ is defined in $E(S)$, so too is $e\phi\cdot f\phi$ in $E(T)$, and $e\phi\cdot f\phi = (ef)\phi$.  An \emph{isomorphism} of biordered sets is a bijective morphism whose inverse is also a morphism.

Now, returning to our earlier notation, let $T = M\times_\Phi^qS$ be a tight twisted product.
We wish to describe the biordered set $E(T)$ in terms of $E(S)$ and $E(M)$.
This amounts to describing all the products that exist in $E(T)$, which involves characterising the $\larr$ and $\rarr$ pre-orders.  
This characterisation is in terms of the corresponding pre-orders in $E(S)$ and $E(M)$.  But note that since $M$ is commutative, the biordered set $E(M)$ is a semilattice, and so ${\larr}={\rarr}$ in $E(M)$; this coincides with the usual partial order $\leq$, given by ${i\leq j \iff i = i+j(=j+i)}$, which is in turn the restriction of $\leqJ^M$ to $E(M)$.

\newpage

\begin{prop}\label{prop:BS}
Let $(i,e),(j,f)\in\Om$.  Then
\ben
\item \label{BS1} $\ve(i,e) \larr \ve(j,f) \iff i\leq j$ and $e\larr f$, in which case $\ve(j,f)\ve(i,e) = \ve(i,fe)$,
\item \label{BS2} $\ve(i,e) \rarr \ve(j,f) \iff i\leq j$ and $e\rarr f$, in which case $\ve(i,e)\ve(j,f) = \ve(i,ef)$.
\een
Consequently, the product $\ve(i,e)\ve(j,f)$ is defined in the biordered set $E(T)$ if and only if $(i,e)(j,f)$ is defined in the direct product $E(M)\times E(S)$, in which case $\ve(i,e)\ve(j,f) = \ve(i+j,ef)$.
\end{prop}

\pf
The final assertion follows by combining parts \ref{BS1} and \ref{BS2}.  Thus, by symmetry, it suffices to prove \ref{BS1}.  For the rest of the proof we write $s=p(i,e)\in H_i^M$ and ${t=p(j,f)\in H_j^M}$, so that
\[
\ve(i,e) = (s,e) \COMMA \ve(j,f) = (t,f) \COMMA i = s+\Phi(e,e)q \AND j = t+\Phi(f,f)q.
\]
Suppose first that $\ve(i,e) \larr \ve(j,f)$.  Then
\[
(s,e) = \ve(i,e) = \ve(i,e)\ve(j,f) = (s,e)(t,f) = (s+t+\Phi(e,f)q,ef).
\]
It follows that $e=ef$ (i.e.~$e\larr f$) and $s=s+t+\Phi(e,f)q \leqJ^M t$.  Since $s\J^Mi$ and $t\J^Mj$, we deduce that $i \leqJ^M j$, so indeed $i\leq j$.

Conversely, suppose $i\leq j$ and $e\larr f$, meaning that ${i = i+j}$ and $e=ef$.  We must show that
\begin{align}
\label{eq:iejf1} \ve(i,e)\ve(j,f) &= \ve(i,e), \qquad\text{i.e.}\qquad \ve(i,e) \larr \ve(j,f),\\
\label{eq:iejf2} \text{and}\qquad \ve(j,f)\ve(i,e) &= \ve(i,fe).
\end{align}
For \eqref{eq:iejf1}, we first note that
\begin{equation}\label{eq:iejf11}
\ve(i,e)\ve(j,f) = (s,e)(t,f) = (s+t+\Phi(e,f)q,ef) = (s+t+\Phi(e,f)q,e).
\end{equation}
From $e=ef$ we have $e\leqL f$, so it follows from Lemma \ref{lem:ae} that ${\Phi(e,f)=\Phi(f,f)}$.  Since $i$ is the identity element of the group $H_i^M=H_s^M$, we also have $s+j = (s+i)+j = s+(i+j) = s+i = s$.  Combining the previous two conclusions with \eqref{eq:iejf11} and the definition of $t$ yields
\[
\ve(i,e)\ve(j,f) = (s+t+\Phi(e,f)q,e) = (s+t+\Phi(f,f)q,e) = (s+j,e) = (s,e) = \ve(i,e).
\]
We can now deduce \eqref{eq:iejf2} fairly quickly.  Since $\ve(i,e) \larr \ve(j,f)$, it follows that $\ve(j,f)\ve(i,e)$ is an idempotent in the $\L^T$-class of $\ve(i,e)$.  This product therefore has the form $\ve(i,g)$ for some idempotent $g\in E(S)$.  This idempotent must in fact be $g=fe$, since the second coordinate in $\ve(j,f)\ve(i,e) = (t,f)(s,e)$ is simply $fe$.
\epf

As a consequence, we have the following:

\begin{thm}\label{thm:BS1}
For a tight twisted product $T = M\times_\Phi^qS$, the map
\[
E(T) \to E(M) \times E(S): \ve(i,e) \mt (i,e) \qquad\text{for $(i,e)\in\Om$}
\]
is a well-defined embedding of biordered sets.  It is an isomorphism if and only if one of the following holds:
\ben
\item \label{BS1a} $q$ is a unit of~$M$, or 
\item \label{BS1b} $\Phi(e,e)=0$ for all idempotents $e\in E(S)$.
\een
\end{thm}

\pf
Since $\ve(i,e)$ is the unique idempotent in the $\H^T$-class of $(i,e)\in\Om$, the map is well defined.  It follows from Proposition \ref{prop:BS} that it is a morphism, and it is clearly injective.

For the second assertion, first note that the map is surjective (and hence bijective) if and only if ${\Om = E(M)\times E(S)}$.  By~Proposition \ref{prop:ET}, this is equivalent to having
\[
i \leqJ^M \Phi(e,e)q \qquad\text{for all $i\in E(M)$ and $e\in E(S)$.}
\]
Since $i\leqJ^M0$ for all $i\in E(M)$, this is equivalent to having $0\leqJ^M\Phi(e,e)q$ for all $e\in E(S)$, which is in turn equivalent to $\Phi(e,e)q$ being a unit for all such $e$.  This is equivalent to one of~\ref{BS1a} or~\ref{BS1b} holding.  Finally we note that when the map is surjective, its inverse $(i,e)\mt\ve(i,e)$ is a morphism as well, as follows from the final assertion of Proposition \ref{prop:BS}.
\epf

\begin{cor}\label{cor:BS2}
If $M$ is unipotent and $q$ is a unit (e.g.~if $M$ is a group), then the biordered sets $E(T) = E(M\times_\Phi^qS)$ and $E(S)$ are isomorphic.
\end{cor}

\pf
By Theorem \ref{thm:BS1} we have $E(T) \cong E(M)\times E(S) = \{0\}\times E(S) \cong E(S)$.
\epf

\begin{rem}\label{rem:0e}
As in Remark \ref{rem:ETunit}, note that when $M$ is a group, each idempotent of $T$ has the form $\ve(0,e) = (-\Phi(e,e)q,e)$ for $e\in E(S)$.  Denoting this idempotent by $\ol e = \ve(0,e)$, the map $E(S)\to E(T):e\mt\ol e$ determines an isomorphism of biordered sets.  Thus, $\ol e\ol f = \ol{ef}$ for any basic pair $(e,f)$ in $E(S)$.  This does not extend to non-basic products in general, however.  That is, if~$ef$ happens to be an idempotent of $S$, then we need not have $\ol e\ol f = \ol{ef}$ in $T$; in fact, $\ol e\ol f$ might not even be an idempotent of $T$.

For a concrete instance of this, consider $T = \Z\times_\Phi^1\P_n$ for $n\geq3$, where $\Phi$ is the canonical twisting on the partition monoid $\P_n$, and where the group $M$ is $(\Z,+)$ and we take $q=1$.  Also consider the idempotents
\[
e = \custpartn{1,2,3,4,7}{1,2,3,4,7}{\uarc12\darc12\stline33\stline44\stline77\udotted47\ldotted47}
\COMMA
f = \custpartn{1,2,3,4,7}{1,2,3,4,7}{\uarc23\darc23\stline11\stline44\stline77\udotted47\ldotted47}
\AND
ef = \custpartn{1,2,3,4,7}{1,2,3,4,7}{\uarc12\darc23\stline31\stline44\stline77\udotted47\ldotted47},
\]
all from $\P_n$.  We then have $\ol e=(-1,e)$, $\ol f=(-1,f)$ and $\ol{ef} = (0,ef)$, yet $\ol e\ol f = (-2,ef)$.
\end{rem}

\section{Regularity}\label{sect:reg}

We again fix a tight twisted product $T = M\times_\Phi^qS$, and consider the regular elements of~$T$, as well as regularity of $T$ itself.  
Analogously to our treatment of idempotents in Section \ref{sect:ET}, where we first described the group $\H^T$-classes, here we begin by characterising the regular $\D^T$-classes before focussing on the regular elements themselves.

For a regular $\D^S$-class $D$ (of $S$), let
\[
\Phi(D) = \min\set{\Phi(e,e)}{e\in E(D)}.
\]
Note that this is somewhat different in nature to the $\Phi(H)$ parameters defined in Section \ref{sect:ET}; specifically, it is possible for the set ${\set{\Phi(e,e)}{e\in E(D)}}$ to have size greater than $1$.

\begin{thm}\label{thm:D}
Let $T = M\times_\Phi^qS$ be a tight twisted product.  For an ${\H^M}(={\D^M})$-class $H$ and a $\D^S$-class $D$, the following are equivalent:
\ben
\item \label{D1} the $\D^T$-class $H\times D$ is regular,
\item \label{D2} $H$ is a group, $D$ is regular, and $H = H+\Phi(D)q$,
\item \label{D3} $H$ is a group, $D$ is regular, and $i\leqJ^M\Phi(D)q$ for some (equivalently all) $i\in H$,
\item \label{D4} $H$ is a group, $D$ is regular, and either $\Phi(D)=0$ or else $i\leqJ^Mq$ for some (equivalently all)~${i\in H}$.
\een
Every group $\H^T$-class contained in a regular $\D^T$-class $H\times D$ is isomorphic to the direct product $(H,+)\times(G,\cdot)$ for some (equivalently, any) group $\H^S$-class $G\sub D$.    
\end{thm}

\pf
By considering a product of the form $(i,a)(j,b)(i,a)$, it is clear that for $H\times D$ to be regular, it is necessary for $H$ and $D$ to be regular, and in particular for~$H$ to be a group.  So suppose this is the case, and fix some $i\in H$.  
Keeping in mind that $H$ is a group, we see that
\begin{align*}
H\times D\text{ is regular} &\iff H\times D\text{ contains an idempotent} \\
&\iff H\times H_e^S \text{ is a group for some $e\in E(D)$}\\
&\iff  i\leqJ^M\Phi(e,e)q \text{ for some $e\in E(D)$} && \text{by Theorem \ref{thm:H}}\\
&\iff  i\leqJ^M\Phi(D)q &&\text{as $\Phi(D)\leq\Phi(e,e)$}\\
&\iff H=H+\Phi(D)q && \text{by Lemma \ref{lem:Hk}.}
\end{align*}
This shows that \ref{D1}$\iff$\ref{D2}$\iff$\ref{D3}.  The equivalence \ref{D3}$\iff$\ref{D4} is exactly as in the proof of Theorem~\ref{thm:H}, which also gives the final assertion.
\epf

\begin{rem}
Consider the case that a $\D^T$-class $H\times D$ is regular and no element of~$H$ is $\leqJ^M$-below $q$.  It then follows from Theorem \ref{thm:D} that $\Phi(D) = 0$, and so $\Phi(e,e) = 0$ for some $e\in E(D)$.  Since $H\times D$ is regular, it follows from general semigroup structure theory that every $\R^T$- and $\L^T$-class contained in $H\times D$ contains an idempotent.  Such an idempotent has the form $(j,f)$ for some $j\in H$ and $f\in E(D)$.  Expanding the product $(j,f)(j,f)$, we see that $j = 2j + \Phi(f,f)q$; since $j\not\leqJ^M q$, this forces $\Phi(f,f) = 0$.  So we are led to the following statement (and its dual):
\bit
\item If $\Phi(e,e)=0$ for some idempotent $e\in E(S)$, then every $\L^S$-class contained in $D_e$ contains an idempotent $f$ such that $\Phi(f,f) = 0$.
\eit
This statement refers only to the monoid $S$ and its (tight) twisting $\Phi$, and says nothing about the monoid $M$ (or $T$).  One can in fact prove this statement directly without reference to the structure of $M\times_\Phi^qS$, as follows.  

Fix some $\L^S$-class $L\sub D_e$.  Since $e\in D_e$ we can fix some $a\in L$ with $a\R^Se$.  We then have $e = ab$ for some $b\in S$, and by \ref{T2} we can assume that $\Phi(a,b) = 0$.  We then set $f = ba$.  Note first that $a=ea$ (as $a\R^Se$) and so $af = aba = ea = a$, so that $f\in L_a^S = L$; we also have $f^2=baf=ba=f$.  We additionally claim that $\Phi(f,f) = 0$.  Indeed, we have
\begin{align*}
\Phi(f,f) = \Phi(ba,ba) &= \Phi(a,ba) &&\text{by \ref{T3}, as $a\L^S f=ba$}\\
&\leq \Phi(a,ba)+\Phi(b,a) \\
&= \Phi(a,b) + \Phi(ab,a) &&\text{by \ref{T1}}\\
&= 0 + \Phi(ab,ab) &&\text{by \ref{T4}, as $a\R^S e=ab$}\\
&= \Phi(e,e) = 0.
\end{align*}
Since also $\Phi(f,f) \geq 0$, it follows that indeed $\Phi(f,f) = 0$.
\end{rem}

We now give a direct characterisation of the regular elements of $T$:

\begin{cor}\label{cor:RegT}
For a tight twisted product $T = M\times_\Phi^qS$, we have
\[
\Reg(T) = \set{(i,a)\in\Reg(M)\times\Reg(S)}{i\leqJ^M\Phi(D_a^S)q}.
\]
\end{cor}

\pf
Observe that $(i,a)\in T$ is regular if and only if its $\D^T$-class $D_{(i,a)}^T = H_i^M\times D_a^S$ is regular.  We now apply Theorem \ref{thm:D}.
\epf

We can also say precisely when the entire twisted product $T$ is regular:

\begin{thm}\label{thm:reg}
For a tight twisted product $T = M\times_\Phi^qS$, the following are equivalent:
\ben
\item $T$ is regular,
\item $M$ and $S$ are both regular, and either $q$ is a unit or else $\Phi(D)=0$ for every $\D^S$-class $D$.
\een
\end{thm}

\pf
It follows from Corollary \ref{cor:RegT} that $T$ is regular if and only if $M$ and $S$ are regular and $i\leqJ^M\Phi(D)q$ for every $i\in M$ and every $\D^S$-class $D$.  
As in the proof of Theorem \ref{thm:BS1}, this is equivalent to $\Phi(D)q$ being a unit for every $\D^S$-class $D$, and hence to $q$ being a unit or having $\Phi(D)=0$ for all $D$.
\epf

\begin{rem}
One can examine Figures \ref{fig:P2Phi} and \ref{fig:B34Phi}, and locate the regular $\D$-classes in the various tight twisted products pictured there.  For example, consider $M \times_\Phi^\infty \P_2$, where $M=\{0,\infty\}$ and~$\Phi$ is the canonical twisting, as shown in Figure \ref{fig:P2Phi} (middle).  The $\D$-classes of $\P_2$ are ${D_0<D_1<D_2}$, and we have
\[
\Phi(D_2) = \Phi(D_1) = 0 \AND \Phi(D_0) = 1.
\]
Since $0\not\leqJ^M\infty=q$, it follows that $\{0\}\times D_2$ and $\{0\}\times D_1$ are regular, but $\{0\}\times D_0$ is not.  On the other hand, $\{\infty\}\times D_r$ is regular for each $r$.

Next consider the products arising from the canonical twisting on the Brauer monoids $\B_3$ and $\B_4$; see Figure \ref{fig:B34Phi} (middle).  These monoids have $\D$-classes $D_1<D_3$ and $D_0<D_2<D_4$, respectively, and we have
\[
\Phi(D_3)=\Phi(D_1)=0 \COMMA \Phi(D_4)=\Phi(D_2)=0 \AND \Phi(D_0)=1.
\]
Thus, $\{0\}\times D_0$ is a non-regular $\D$-class of $M\times_\Phi^\infty\B_4$, but every other $\D$-class of $M\times_\Phi^\infty\B_4$ and $M\times_\Phi^\infty\B_3$ is regular.  In particular, $M\times_\Phi^\infty\B_3$ is itself regular.  More generally, $M\times_\Phi^\infty\B_n$ is regular if and only if $n$ is odd.

Finally, when $\Phi$ is the rank-based twisting on $\B_n$, we have $\Phi(D_r) = n-r$.  Consequently, the $\D$-class $\{0\}\times D_r$ of $M\times_\Phi^\infty\B_n$ is only regular when $r=n$.  This can be seen in Figure \ref{fig:B34Phi} (right) for $n=3$ and $4$.
\end{rem}

\section{Idempotent-generated submonoids}\label{sect:IG}

The results of Sections \ref{sect:GR}--\ref{sect:reg} apply to arbitrary tight twisted products $T = M\times_\Phi^qS$, and give detailed structural information about them.  On the other hand, there are many general problems one could pose for twisted products, whose solutions are not uniform, but require situation-specific/ad hoc reasoning.  In this final section we address one of these, namely the determination of the idempotent-generated submonoid~$\la E(T)\ra$.  By considering a number of concrete examples, including some in the existing literature, we will see that the nature of this submonoid depends heavily on the structure of $M$ and $S$, and also on the twisting $\Phi$.  Specifically, we consider rigid twistings in Section \ref{subsect:IGPhirm}, and twisted diagram monoids in Section~\ref{subsect:IGPn}.

For convenience, we will write $\E(S) = \la E(S)\ra$ for the idempotent-generated submonoid of~$S$.  It is well known, and easy to check, that $\E(S)\sm\{1\}$ is a subsemigroup of $\E(S)$, and we will denote it by ${\Ef(S) = \E(S)\sm\{1\}}$.  It is also clear that $\Ef(S) = \la E(S)\sm\{1\}\ra$.  We similarly write $\E(T) = \la E(T)\ra$ and $\Ef(T) = \E(T) \sm \{(0,1)\}$.  But note that $\E(M) = E(M)$ because $M$ is commutative.  

Perhaps the most general statement one can make is the following:

\begin{prop}\label{pr:EfT}
Let $T = M\times_\Phi^qS$ be a tight twisted product.  
\ben
\item \label{EfT0} We have $\E(T) \sub (E(M)\times\{1\}) \cup (M\times\Ef(S))$.
\item \label{EfT1} If $M$ is unipotent, then $\Ef(T) \sub M\times\Ef(S)$.
\item \label{EfT2} If $q$ is a unit, then $\E(T)$ is a homomorphic pre-image of $\E(S)$.
\een
\end{prop}

\pf
\firstpfitem{\ref{EfT0}}  Consider a product of idempotents $a = (i_1,e_1)\cdots(i_k,e_k)$ in $T$.  If each $e_j=1$, then each $i_j$ is an idempotent (cf.~Remark \ref{rem:ETunit}), and $a = (i_1+\cdots+i_k,1)$ with $i_1+\cdots+i_k\in E(M)$ by commutativity.  Otherwise, the second coordinate of $a$ belongs to $\Ef(S)$.

\pfitem{\ref{EfT1}}  This follows from part \ref{EfT0}, as $E(M)\times\{1\} = \{(0,1)\}$.

\pfitem{\ref{EfT2}}  The natural map $\E(T)\to \E(S):(i,a)\mt a$ is surjective, as its image contains all of $E(S)$.  This is because any $e\in E(S)$ is the second coordinate of $\ve(0,e) = (-\Phi(e,e)q,e)$.
\epf

\begin{rem}\label{rem:extreme}
With reference to Proposition \ref{pr:EfT}\ref{EfT0} and \ref{EfT2}, we will see specific examples in which:
\bit
\item $\E(T) = (E(M)\times\{1\}) \cup (M\times\Ef(S))$ (Theorems \ref{thm:IGNPn},  \ref{thm:IGZPn}),
\item $q$ is a unit and $\E(T)\cong\E(S)$ (Theorem \ref{thm:Phirm}).
\eit
In Remark \ref{rem:BTLM} we discuss an intermediate example between these two extremes.
\end{rem}

Before we do this, however, it will be convenient to extend the domain of a twisting $\Phi$ to any number of coordinates.  For $k\geq0$ and $a_1,\ldots,a_k\in S$, we inductively define
\[
\Phi(a_1,\ldots,a_k) = \begin{cases}
0 &\text{if $k\leq1$}\\
\Phi(a_1,\ldots,a_{k-1}) + \Phi(a_1\cdots a_{k-1},a_k) &\text{if $k\geq2$.}
\end{cases}
\]
When $k=2$ or $3$, it is easy to see that this agrees with the numbers $\Phi(a,b)$ and $\Phi(a,b,c)$ we have already encountered.  It is also easy to show by induction that in the twisted product~${M\times_\Phi^qS}$, we have
\[
(i_1,a_1)\cdots(i_k,a_k) = (i_1+\cdots+i_k+\Phi(a_1,\ldots,a_k)q,a_1\cdots a_k).
\]

\subsection{Rigid twistings}\label{subsect:IGPhirm}

We first consider rigid twistings, i.e.~twistings of the form $\Phi_{r,m}$ from Lemma \ref{lem:Phirm}.  Here we do not need to assume that such a twisting is tight.

\begin{thm}\label{thm:Phirm}
If $T = M\times_\Phi^qS$ is a twisted product, where $M$ is a group and $\Phi = \Phi_{r,m}$ is rigid, then $\langle E(T)\rangle \cong \langle E(S)\rangle$.
\end{thm}

\pf
For $a\in S$, we define $j_a = r(a)-m$.  It is easy to check, using the definition of $\Phi = \Phi_{r,m}$, that $j_a + j_b + \Phi(a,b) = j_{ab}$ for all $a,b\in S$.  It follows from this that
\[
W = \set{(j_aq,a)}{a\in S}
\]
is a subsemigroup of $T$, and that $S\to W:a\mt(j_aq,a)$ is an isomorphism.  The result now follows from the claim that $W$ contains $E(T)$, and hence $\E(T)$.  To see that this is the case, note that Remark \ref{rem:0e} gives $E(T) = \set{\ve(0,e)}{e\in E(S)}$, and we observe that 
\[
\ve(0,e) = (-\Phi(e,e)q,e) = (j_eq,e) \qquad\text{for $e\in E(S)$,}
\]
as $\Phi(e,e) = m-r(e)-r(e)+r(e^2) = m-r(e) = -j_e$.  
\epf

\subsection{Twisted diagram monoids}\label{subsect:IGPn}

We now consider the idempotent-generated subsemigroup $\E(T)$, where $T$ is either of the (tight) twisted products $\N\times_\Phi^1\P_n$ or $\Z\times_\Phi^1\P_n$ arising from the partition monoid $\P_n$ and its canonical twisting $\Phi$.  To avoid trivialities, we assume that $n\geq2$ throughout this section.  It will also be convenient to write
\[
E_0(\P_n) = \set{e\in \P_n}{\Phi(e,e)=0},
\]
and we note that $E(\N\times_\Phi^1\P_n) = \set{(0,e)}{e\in E_0(\P_n)}$, as in Remark \ref{rem:ETunit}.

It was shown in \cite[Proposition 16]{JEpnsn} that the singular ideal $\Sing(\P_n) = \P_n\sm\S_n$ is generated by the set
\[
\Si = \set{t_{ij}}{1\leq i<j\leq n} \cup \set{t_i}{1\leq i\leq n},
\]
where these partitions are defined by:
\[
t_{ij} = t_{ji} = \custpartn{1,4,5,6,9,10,11,14}{1,4,5,6,9,10,11,14}{\stline11\stline44\stline66\stline99\stline{11}{11}\stline{14}{14}\stline55\stline{10}{10} \uarc5{10}\darc5{10} \udotted14\ldotted14\udotted69\ldotted69\udotted{11}{14}\ldotted{11}{14}\cplab11\cplab5i\cplab{10}j\cplab{14}n} 
\AND
t_i = \custpartn{1,4,5,6,9}{1,4,5,6,9}{\stline11\stline44\stline66\stline99\udotted14\ldotted14\udotted69\ldotted69\cplab11\cplab5i\cplab9n}.
\]
We will also use the partitions $e_{ij}$, defined for distinct $i,j\in\bn$ by 
\[
e_{ij} = 
\begin{cases}
\custpartn{1,4,5,6,9,10,11,14}{1,4,5,6,9,10,11,14}{\stline11\stline44\stline66\stline99\stline{11}{11}\stline{14}{14}\stline55\stline{10}5 \udotted14\ldotted14\udotted69\ldotted69\udotted{11}{14}\ldotted{11}{14}\cplab11\cplab5i\cplab{10}j\cplab{14}n} &\text{if $i<j$}\\[5mm]
\custpartn{1,4,5,6,9,10,11,14}{1,4,5,6,9,10,11,14}{\stline11\stline44\stline66\stline99\stline{11}{11}\stline{14}{14}\stline5{10}\stline{10}{10} \udotted14\ldotted14\udotted69\ldotted69\udotted{11}{14}\ldotted{11}{14}\cplab11\cplab5j\cplab{10}i\cplab{14}n} &\text{if $i>j$.}
\end{cases}
\]
In what follows, we will use (without mention) the fact that each $t_{ij}$, $e_{ij}$ and $e_{ij}^*$ belong to $E_0(\P_n)$.  This follows from the more general fact (which will also be used without mention) that
\[
\Phi(a,b) = 0 \qquad\text{if $\codom(a) \cup \dom(b) = \bn$,}
\]
which includes the case that either $\codom(a)=\bn$ or $\dom(b) = \bn$.

\begin{lemma}\label{lem:SingPn}
For any $a\in\Sing(\P_n)\cup\{1\}$ we have
\[
a = e_1\cdots e_k \qquad\text{for some $e_1,\ldots,e_k\in E_0(\P_n)$ with $\Phi(e_1,\ldots,e_k) = 0$.}
\]
\end{lemma}

\pf
By the above-mentioned result of \cite{JEpnsn}, we have $a = f_1\cdots f_l$ for some $l\geq0$ and ${f_1,\ldots,f_l\in\Si}$.  We now proceed by induction on $l$.  If $l=0$, we take $k=1$ and $e_1=1$.  So suppose $l\geq1$, and let $b = f_1\cdots f_{l-1}$.  By induction, we have 
\[
b = g_1\cdots g_m \qquad\text{for some $g_1,\ldots,g_m\in E_0(\P_n)$ with $\Phi(g_1,\ldots,g_m) = 0$.}
\]
If $f_l = t_{ij}$ for some $1\leq i<j\leq n$, then $a = g_1\cdots g_mt_{ij}$, with
\[
g_1,\ldots,g_m,t_{ij}\in E_0(\P_n) \AND \Phi(g_1,\ldots,g_m,t_{ij}) = \Phi(g_1,\ldots,g_m) + \Phi(g_1\cdots g_m,t_{ij}) = 0.
\]
We are therefore left to consider the case in which $f_l = t_i$ for some $i\in\bn$, and here we have~$a=bt_i$.  

\pfcase1  If $\{i'\}$ is a block of $b$, then $b=bt_i=a$, so $a = g_1\cdots g_m$, and we are done.

\pfcase2  Next suppose $b$ has a block of the form $A\cup\{i'\}$ for some non-empty subset $A\sub\bn$.  Then for any $j\in\bn\sm\{i\}$ we have $t_i = e_{ji}^*e_{ji}$, and so
\[
a = bt_i = g_1\cdots g_me_{ji}^*e_{ji}.
\]
This factorisation has the desired form, since
\begin{align*}
\Phi(g_1,\ldots,g_m,e_{ji}^*,e_{ji}) &= \Phi(g_1,\ldots,g_m,e_{ji}^*) + \Phi(g_1\cdots g_me_{ji}^*,e_{ji}) \\
&= \Phi(g_1,\ldots,g_m) + \Phi(g_1\cdots g_m,e_{ji}^*) + 0 
= \Phi(b,e_{ji}^*) = 0.
\end{align*}

\pfcase3 If we are not in the above cases, then $i$ belongs to a nontrivial $\coker(b)$-class.  If $j\in\bn\sm\{i\}$ also belongs to this class, then
\[
a = bt_i = be_{ji} = g_1\cdots g_me_{ji},
\]
which again has the required form, since 
\[
\Phi(g_1,\ldots,g_m,e_{ji}) = \Phi(g_1,\ldots,g_m) + \Phi(g_1\cdots g_m,e_{ji}) = 0.  \qedhere
\]
\epf

\begin{lemma}\label{lem:SingPn2}
For any $a\in\Sing(\P_n)$ we have $a = agh$ for some $g,h\in E_0(\P_n)$ with $\Phi(a,g,h) = 1$.
\end{lemma}

\pf
If $\coker(a) \not= \De_\bn$, then for any $(i,j)\in\coker(a)$ with $i\not=j$, we have
\[
a = ae_{ij}e^*_{ij} \AND \Phi(a,e_{ij},e^*_{ij}) = \Phi(a,e_{ij}e^*_{ij}) + \Phi(e_{ij},e^*_{ij}) = 0+1 = 1.
\]
Otherwise, since $\rank(a)<n$ (as $a$ is singular), $a$ must contain some lower singleton, $\{i'\}$ say, and then for any $j\in\bn\sm\{i\}$ we have
\[
a = ae^*_{ji}e_{ji} \AND \Phi(a,e^*_{ji},e_{ji}) = \Phi(a,e^*_{ji}e_{ji}) + \Phi(e^*_{ji},e_{ji}) = 1+0 = 1.  \qedhere
\]
\epf

We can now deal with the twisted monoid $\N\times_\Phi^1\P_n$:

\begin{thm}\label{thm:IGNPn}
If $\Phi$ is the canonical twisting on the partition monoid $\P_n$ $n\geq 2$, then
\[
\bigl\langle E(\N\times_\Phi^1\P_n)\bigr\rangle = \bigl\{(0,1)\bigr\}\cup\bigl(\N\times\Sing(\P_n)\bigr).
\]
\end{thm}

\pf
Throughout the proof we write $T = \N\times_\Phi^1\P_n$, and for $e\in E_0(\P_n)$ write ${\ol e = (0,e) \in E(T)}$.  By Proposition \ref{pr:EfT}\ref{EfT1}, we only need to show that $\N\times\Sing(\P_n) \sub \E(T)$.  If $a\in\Sing(\P_n)$, then with $e_1,\ldots,e_k\in E_0(\P_n)$ as in Lemma \ref{lem:SingPn} we have $(0,a) = \ol e_1\cdots\ol e_k\in\E(T)$.  Thus, we can complete the proof by showing that:
\bit
\item For any $a\in\Sing(\P_n)$ and any $m\in\N$, we have $(m+1,a) = (m,a)\ol g\ol h$ for some $g,h\in E_0(\P_n)$.
\eit
But this follows quickly from Lemma \ref{lem:SingPn2}.  Indeed, if $g,h\in E_0(\P_n)$ are as in that lemma, then
\[
(m,a)\ol g\ol h = (m,a)(0,g)(0,h) = (m+\Phi(a,g,h),agh) = (m+1,a).  \qedhere
\]
\epf

Here is the corresponding result for $\Z\times_\Phi^1\P_n$:

\begin{thm}\label{thm:IGZPn}
If $\Phi$ is the canonical twisting on the partition monoid $\P_n$ with $n\geq 2$, then
\[
\bigl\langle E(\Z\times_\Phi^1\P_n)\bigr\rangle = \bigl\{(0,1)\bigr\}\cup\bigl(\Z\times\Sing(\P_n)\bigr).
\]
\end{thm}

\pf
This time we write $T = \Z\times_\Phi^1\P_n$, and $\ol e = (-\Phi(e,e),e)\in E(T)$ for $e\in E(\P_n)$.  Again, the proof boils down to showing that $\Z\times\Sing(\P_n) \sub \E(T)$.  
Since $T = \Z\times_\Phi^1\P_n$ contains $\N\times_\Phi^1\P_n$, it follows from Theorem \ref{thm:IGNPn} that $\N\times\Sing(\P_n)\sub\E(T)$.  Thus, we are left to show that:
\bit
\item For any $a\in\Sing(\P_n)$ and any $m\in\Z$, we have $(m-1,a) = (m,a)\ol g\ol h$ for some $g,h\in E(\P_n)$.
\eit
To prove this, suppose first that $\coker(a)\not=\De_\bn$, and fix some $(i,j)\in\coker(a)$ with $i\not=j$.  Then
\[
(m,a)\ol t_i\ol t_{ij} = (m,a)(-1,t_i)(0,t_{ij}) = (m-1,a),
\]
since we have $at_it_{ij}=a$ and $\Phi(a,t_i,t_{ij})=0$.

Now suppose $\coker(a)=\De_\bn$.  As in the proof of Lemma \ref{lem:SingPn2}, $a$ contains a lower singleton,~$\{i'\}$~say.  This time we take $g=t_{ij}$ and $h=t_i$ for any $j\in\bn\sm\{i\}$.  
\epf

\begin{rem}\label{rem:BTLM}
Results analogous to Theorems \ref{thm:IGNPn} and \ref{thm:IGZPn} hold for the Brauer monoid~$\B_n$.  Specifically, for $n\geq3$ we have
\[
\bigl\langle E(\N\times_\Phi^1\B_n)\bigr\rangle = \bigl\{(0,1)\bigr\}\cup\bigl(\N\times\Sing(\B_n)\bigr)\ANd
\bigl\langle E(\Z\times_\Phi^1\B_n)\bigr\rangle = \bigl\{(0,1)\bigr\}\cup\bigl(\Z\times\Sing(\B_n)\bigr).
\]
Indeed, the first of these is \cite[Theorem 3.23]{DE2018}, while the second can be deduced from the first, using the following:
\bit
\item For any $a\in\Sing(\B_n)$ and any $m\in\Z$, we have $(m-1,a) = (m,a)\ol g\ol h$ for some $g,h\in E(\B_n)$.
\eit
To prove this, note first that $a$ must have some lower block $\{i',j'\}$.  Then for any $k\in\bn\sm\{i,j\}$, and writing $\bn\sm\{i,j,k\} = \{l_1,\ldots,l_{n-3}\}$, we take
\[
g = \begin{partn}{5}
i&j,k&l_1&\cdots&l_{n-3}\\\hhline{~|-|~|~|~}i&j,k&l_1&\cdots&l_{n-3}
\end{partn}
\AND
h = \begin{partn}{5}
j&i,k&l_1&\cdots&l_{n-3}\\\hhline{~|-|~|~|~}k&i,j&l_1&\cdots&l_{n-3}
\end{partn},
\]
noting that $\ol g = (-1,g)$ and $\ol h = (0,h)$.  The claim follows because $a=agh$ and $\Phi(a,g,h)=0$.

The situation for the Temperley--Lieb monoid $\TL_n$ is far more intricate.  A description of $\bigl\langle E(\N\times_\Phi^1\TL_n)\bigr\rangle=\Ef(\N\times_\Phi^1\TL_n)\cup\{ (0,1)\}$ is the main result of \cite{DE2017}.  Specifically, it was shown there that $(i,a)$ belongs to $\Ef(\N\times_\Phi^1\TL_n)$ precisely when $i\geq\chi(a)$ and $\chi(a) \equiv i$ (mod $2$), where here $\chi(a)$ is a certain parameter defined in terms of the so-called \emph{Jones normal form} of~$a$, which was itself introduced in~\cite{BDP2002}.  The definitions of these normal forms, and the associated $\chi$ parameters, are rather involved, so for reasons of space we will not repeat them here.  But we can at least observe that:
\bit
\item $\Ef(\N\times_\Phi^1\TL_n)$ is a proper subset of $\N\times\Sing(\TL_n)$, and
\item $\Ef(\N\times_\Phi^1\TL_n)$ is infinite, and maps onto $\Ef(\TL_n)$.
\eit
Thus, $\Ef(\N\times_\Phi^1\TL_n)$ lies between the two `extremes' discussed in Remark \ref{rem:extreme}.

Things are potentially even more complicated for other diagram monoids, such as the partial Brauer and Motzkin monoids, $\PB_n$ and $\M_n$.  For one thing, their singular parts are not idempotent-generated \cite{DEG2017}.  For another, their twistings are loose, by Proposition \ref{prop:canPhi}.  We believe it would be interesting, and challenging, to determine the idempotent-generated submonoids of the associated twisted products, over both $\N$ and $\Z$, but this is beyond the scope of the current work.
\end{rem}

\footnotesize
\def\bibspacing{-1.1pt}
\bibliography{biblio}
\bibliographystyle{abbrv}

\end{document}